\documentclass{article}

\usepackage[english]{babel}

\usepackage[letterpaper,top=2cm,bottom=2cm,left=3cm,right=3cm,marginparwidth=1.75cm]{geometry}

\usepackage[mathscr]{eucal}
\usepackage{amsmath}
\usepackage{amsthm}
\usepackage{amscd}
\usepackage{amssymb}
\usepackage{latexsym}
\usepackage{graphicx}
\usepackage{color}

\def\vcorrection#1{\advance\voffset#1\relax}
\def\hcorrection#1{\advance\hoffset#1\relax}

\theoremstyle{plain}
\newtheorem{theorem}{Theorem}

\newtheorem*{Proof}{Proof}

\theoremstyle{definition}

\usepackage[colorlinks=true, allcolors=blue]{hyperref}

\begin{document}

\title{ Invariant Differential Operators 
and the Radon Transform 
\\ 
on the Horocycle Spaces   
 }
\author{Satoshi Ishikawa
\thanks{
Polytechnic University of Japan, 2-32-1 Ogawanishi-Cho Kodaira-Shi Tokyo 187-0035, Japan; 
e-mail: ishikawa@uitec.ac.jp
}
}
\date{ }

\maketitle

\begin{abstract}

We investigate the Radon transform for double fibrations 
of the horocycle spaces for the semisimple symmetric spaces 
with respect to the inclusion incidence relations. 
We present the inversion formula, support theorem 
and the range theorem by the invariant differential operators 
or the invariant system of differential operators 
constructed from the left action of the Pfaffian type elements 
in the universal enveloping algebra for the transformations group. 
In order to prove the range theorem, we make 
the explicit calculations of the Pfaffian type elements 
which lead to the calculations for 
the eigenvalue  of the central elements 
of the universal enveloping algebra. 
We deal with the Radon transform on 
the Schwartz space, the compactly supported 
smooth function spaces and the space 
of the sections of the line bundle. 
The range theorem on the space 
of the sections of the line bundle 
yields the range theorem 
for the Radon transform for double fibrations 
of compact homogeneous spaces which are not necessarily symmetric.

\end{abstract}

\section{Introduction}

Let $U$ be a connected unimodular Lie group 
and $L_1,L_2$ be closed unimodular subgroups. 
The Radon Transform for a double fibration 
$$U/L_1\leftarrow U/(L_1\cap L_2)\rightarrow U/L_2$$ 
where $U/L_1$ and $U/L_2$ are homogeneous spaces in duality 
is defined by S.Helgason in [17-20] as follows,  

\begin{equation}\tag{1.0.1}
(R.f)(g L_2)
=\int\nolimits_{L_2/(L_1\cap L_2)}f(g h L_1)dh_{L_1\cap L_2}
\text{ for }g\in G
\end{equation}
for $f\in C^\infty_c(U/L_1)$     
where $dh_{L_1\cap L_2}$ is the $L_2$- invariant measure 
on $L_2/L_1\cap L_2$. 

For $\varphi\in C^\infty_c(U/L_2)$, 
we also define the dual Radon transform $R^*$ of $R$ by 

\begin{equation}\tag{1.0.2}
(R^*.f)(g L_1)
=\int\nolimits_{L_1/(L_1\cap L_2)}\varphi(g h L_2)dh_{L_1\cap L_2}
\text{ for }g\in G. 
\end{equation}

The principal problems for $R$ or $R^*$ which are posed by S.Helgason 
are  
(A) Injectivity of $R$ or $R^*$, 
(B) Inversion formulas of $R$ or $R^*$, 
(C) Characterizations of the range of $R$ or $R^*$, 
(D) Support Theorem of $R$ or $R^*$ etc  
(See Chapter I in [17],
Chapter II in [18], Chapter I in [19] or Chapter II in [20]).

\subsection{Basic Notation} 

Let $G$ be $SO_0(m,n)(m\geq n\geq 1)$, namely the identity component 
of $O(m,n)$, 
where we put 

\begin{equation*} I_{m,n}=\begin{pmatrix} 1_m & 0 \\ 0 & -1_n\end{pmatrix} 
\text{ where } 1_m\text{ or }1_n 
\text{is the identity matrix of }
GL(m,\mathbb R)\text{ or }GL(n,\mathbb R).
\end{equation*} 
and 
\begin{equation*}SO(m,n)=\{g\in GL(m+n;\mathbb R)\Bigm\vert g^t\cdot I_{m,n}\cdot g=I_{m,n},\det g=1\}.\end{equation*} 

We define 
\begin{equation*} K=\{\begin{pmatrix} u_1 & 0 \\ 0 & u_2\end{pmatrix} \in G  \Bigm\vert
u_1\in O(m),u_2\in O(n)\},\end{equation*} 
and for $0\leq r\leq m$, 
\begin{equation*} H^{(r)}=\{\begin{pmatrix} u & 0 \\ 0 & g\end{pmatrix} \in G  \Bigm\vert
u\in O(m-r),g\in O(r,n)\}.\end{equation*} 

Let $\mathfrak g$ be $\mathfrak s\mathfrak o(m,n)$ namely, 
the Lie algebra of $G$. 
Let $\mathfrak k$ be 
$\mathfrak s\mathfrak o(m)\times\mathfrak s\mathfrak o(n)$ 
then $\mathfrak k$ is a maximal compact subalgebra.  
Let $\mathfrak p$ be the orthogonal complement of $\mathfrak k$ 
with respect to the Killing form of $\mathfrak g$. 
Let $\mathfrak h^{(r)}$ be the Lie algebra of $H^{(r)}$. 
Let $\mathfrak q_{\mathfrak h^{(r)}}$ be 
the orthogonal complement of $\mathfrak h^{(r)}$ 
with respect to the Killing form of $\mathfrak g$. 
\par 
Then $H^{(r)} K$ is a totally geodesic submanifold in $G/K$ 
(See \S 7,Chapter IV in [16] and (0) in [24].) 
Here the semisimple symmetric space 
$G/H^{(r)}$ can be identified 
with the space of all the $G$-translations of $H^{(r)} K$ 
(See Proposition 4-1 in [24]).

Assume $0\leq r\leq m-n$. Let $\mathfrak a$ be a maximal abelian subspace 
in $\mathfrak p\cap\mathfrak q_{\mathfrak h^{(r)}}$ 
which is also maximal abelian in $\mathfrak p$. 
Let $M$ be the centralizer $\mathfrak a$ in $K$, 
$\mathfrak m$ be the Lie algebra of $M$, 
$\Sigma$ be the set of all the roots, 
$\Sigma^+$ be a set of positive roots. 
We put $A=\exp(\mathfrak a)$. We also put 
$$\mathfrak g_\alpha=\{X\in\mathfrak g\Bigm\vert 
[A,X]=\alpha(A)X\text{ for }A\in\mathfrak a\},
\mathfrak n=\bigoplus\limits_{\alpha\in\Sigma^+}\mathfrak g_\alpha
\text{ and }N=\exp(\mathfrak n).$$ 
Then $(MAN)H^{(r)}$ is open dense in $G/H^{(r)}$ 
(See \S 3,Proposition 1 in [33],Theorem 13 in [36] and p151 in [26]). 
\par 
Let $\mathcal H_r$ for $0\leq r\leq m-n$ 
be the space of all the $G$-translations 
of the horocycle (horosphere) 
$((M\cap H^{(r)})N) H^{(r)}$ in $G/H^{(r)}$. 
Then the space $\mathcal H_r$ can be identified 
with the homogeneous spaces $G/((M\cap H^{(r)})N)$
(See Proposition B2 in [31]) 
and does not depend on the choice of $A$ and $N$(See Proposition 2-2). 
We call $\mathcal H_r$ the horocycle space 
for the semisimple symmetric space $G/H^{(r)}$. 
\par 
The horocycle space $\mathcal H_0$  
for the Riemannian symmetric space $G/K$ 
can be identified with the homogeneous space $G/(MN)$ or with the space 
of all the orbits $g_1 N g_2 K$ in $G/K$ for $g_1,g_2\in G$
(See ChII,\S 1,Theorem 1.1(ii) in [19])).  
We call an orbit $g N K$ in $G/K$ for $g\in G$ 
a Riemannian horocycle. 
An element $g MN\in G/(MN)$ with $g\in G$ 
can be written $g MN=u a MN$ with $u\in K$ and $a\in A$ 
where $u M\in K/M$ and $a\in A$ are unique. 
The composite distance of the Riemannian horocycle 
$u a N K$ from $e K$ is $\log a\in\mathfrak a$ 
(See Ch II,\S 1,p64 in [19]). 
We call the point $u a K\in G/K$ 
the composite distance point of the Riemannian horocycle $u a N K$. The simultaneous Euclidean imbedding 
of $G/K$ and of $G/(MN)$ with horocycles as plane sections are known in ChII,\S 4,p64 in [19].  
\par 
A element $g(M\cap H^{(r)})N\in G/((M\cap H^{(r)})N)$ with $g\in G$ 
can be written $g(M\cap H^{(r)})N=u a (M\cap H^{(r)})N$ 
with $u\in K$ and $a\in A$ 
where $u (M\cap H^{(r)})\in K/(M\cap H^{(r)})$ and $a\in A$ are unique. 
Then the horocycle space $\mathcal H_r$ is diffeomorphic to $K/(M\cap H)\times A$. 
Hence the the horocycle space $\mathcal H_r$ can be identified with  
the following manifold  
$$\mathcal H(r)
=\{(u a MN,u a H^{(r)})\in G/(MN) \times G/H^{(r)}\bigm\vert u\in K,a\in A\}$$ 
by the following diffeomorphism(See Proposition 2-4), 
\begin{equation}\tag{1.1}
i_r(u a (M\cap H^{(r)})N)=(u a MN,u a H^{(r)})
\text{ for }u\in K,a\in A.
\end{equation} 
Here $u a M N K$ is a Riemannian horocycle in $G/K$ 
and $u a H^{(r)} K$ is a totally geodesic submanifold in $G/K$ which is tangent to the Riemannian horocycle $u a M N K$ at the composite distance point $u a K$
(See Proposition 2-5). 
\par 
We put $S=H^{(r_1)}$ and  $T=H^{(r_2)}$ 
for $n\geq 1; 0\leq r_1,r_2 \leq m-n; r_1\not= r_2$.  
Let $\mathfrak a$ be a maximal abelian subspace 
in $\mathfrak p\cap\mathfrak q_{\mathfrak h^{(r_1)}}\cap \mathfrak q_{\mathfrak h^{(r_2)}}$ 
which is also maximal abelian in $\mathfrak p$. 
\par 
We define the Radon transform $R$ which maps $C^\infty_c(\mathcal H_{r_1})$ 
into $C^\infty(\mathcal H_{r_2})$ and the dual Radon transform $R^*$ which maps 
$C^\infty_c(\mathcal H_{r_2})$ into 
$C^\infty(\mathcal H_{r_1})$ 
as the Radon transform in (1.0.1) 
and its dual in (1.0.2) for the double fibrations 
$$G/((M\cap S)N)\leftarrow G/((M\cap S\cap T)N)
\rightarrow G/((M\cap T)N)$$ 
under the identification of $\mathcal H_{r_1}$ or 
$\mathcal H_{r_2}$ 
with $G/((M\cap H^{r_1})N)$ or $G/((M\cap H^{r_2})N)$. 
Here $R$ or $R^*$ does not depend on the choice of $A$ and $N$ (See Proposition 2-3). 
\par     
Then the two horocycles $u_1 a_1 (M\cap S)N\in G/((N\cap S)N)$ and $u_2 a_2(M\cap T)N\in G/((M\cap T)N)$ 
where $u_1\in K/(M\cap S), a_1\in A$ 
and $u_2\in K/(M\cap T), a_2\in A$ 
are incident if and only if $u_1 a_1 MN=u_2 a_2 MN$ 
and $u_1 a_1 S K\subset u_2 a_2 T K$ if $r_1<r_2$ 
or if and only if  $u_1 a_1 MN=u_2 a_2 MN$ 
and $u_1 a_1 S K\supset u_2 a_2 T K$ if $r_1>r_2$.  
\par 
Namely the two horocycles $u_1 a_1 (M\cap S)N\in G/((N\cap S)N)$ 
and $u_2 a_2 G/((M\cap T)N)$ are incident 
if and only if the two Riemannian horocycles 
$u_1 a_1 M N K$ and $u_2 a_2 M N K$ 
coincides and the two totally geodesic submanifolds 
$u_1 a_1 S K$ and $u_2 a_2 T K$ 
have inclusion incidence relations.  
\par 
Let $(\xi,\gamma_T)\in \mathcal H(r_2)$ 
where $\xi\in\mathcal H_{0}$ and $\gamma_T\in G/T$ 
is a totally geodesic submanifold in $G/K$ 
which is tangent to the Riemannian horocycle $\xi$ 
at the composite distance point of $\xi$. 
Then we have for $f\in C^\infty_c(\mathcal H(r_1))$, 
\begin{equation*}
(R f)(\xi,\gamma_T) 
=\begin{cases}  \int\nolimits_{(\xi,\gamma _S)
\in\mathcal H(r_1),\gamma_S\subset\gamma_T}
         f((\xi,\gamma_S)d\mu_{(\xi,\gamma_T)}((\xi,\gamma_S))
 &  \text{ if } r_1<r_2    \\ 
         \int\nolimits_{(\xi,\gamma _S)\in\mathcal H(r_1),\gamma_T\subset\gamma_S}
         f((\xi,\gamma_S)d\mu_{(\xi,\gamma_T)}((\xi,\gamma_S))
        & \text{ if } r_1>r_2
        \end{cases}
\end{equation*}
where we denote by $d\mu_{(\xi,\gamma_T)}((\xi,\gamma_S))$ 
the canonical measure  on the set 
of all the totally geodesic submanifolds in $G/S$ 
such that $(\xi,\gamma _S)\in\mathcal H(r_1)$ 
and such that $\gamma_S\subset\gamma_T$ if $r_1<r_2$ 
or  $\gamma_T\subset\gamma_S$ if $r_1>r_2$. 
\par 
Moreover, let $(\xi,\gamma_S)\in \mathcal H(r_1)$ 
where $\xi\in\mathcal H_{0}$ and $\gamma_S\in G/S$ 
is a totally geodesic submanifold in $G/K$ 
which is tangent to the Riemannian horocycle $\xi$ 
at the composite distance point of $\xi$. 
Then we have for $f\in C^\infty_c(\mathcal H(r_2))$, 
\begin{equation*}
(R^* f)(\xi,\gamma_S)  
=\begin{cases}  \int\nolimits_{(\xi,\gamma_T)\in\mathcal H(r_2),\gamma_S\subset\gamma_T}
         f((\xi,\gamma_T)d\mu_{(\xi,\gamma_S)}((\xi,\gamma_T))
 &  \text{ if } r_2>r_1    \\ 
         \int\nolimits_{(\xi,\gamma_T)\in\mathcal H(r_2),\gamma_T\subset\gamma_S}
         f((\xi,\gamma_T)d\mu_{(\xi,\gamma_S)}((\xi,\gamma_T))
        & \text{ if } r_2<r_1
\end{cases}
\end{equation*}
where we denote by $d\mu_{(\xi,\gamma_S)}((\xi,\gamma_T))$ 
the canonical measure on the set of all the totally geodesic submanifolds in $G/T$ 
such that $(\xi,\gamma_T)\in\mathcal H(r_2)$ 
and such that $\gamma_T\subset\gamma_S$ if $r_2<r_1$ 
or  $\gamma_S\subset\gamma_T$ if $r_2>r_1$.  
\par 
Let $R_K$ or $R_K^*$ be the Radon transform in (1.0.1) 
or its dual in (1.0.2) 
for double fibrations of compact homogeneous spaces(See Proposition 5-1(2) and Proposition 5-3(1))   
$$K/(M\cap S)\leftarrow K/(M\cap S\cap T)
\rightarrow K/(M\cap T).$$ 
Then the Radon transform $R_K$ 
maps $C^\infty(K/(M\cap S))$ into $C^\infty(K/(M\cap T))$. 
Here $K/(M\cap S)$ can be identified 
with the space of all the horocycles 
whose Riemannian horocycles 
have the composite distance $\log a$ for any fixed $a\in A$. 
In view of the diffeomorphism 
between $G/((M\cap H)N)$ and $K/(M\cap H)\times A$, 
the Radon transform $R_K$ can be identified with 
the restricted Radon transform $R$ 
on $Ka(M\cap H)$ which is diffeomorphic 
to $K/(M\cap H)$.  
Here we have 
$$(R f)(ua(M\cap T)N)=(R_{K}f)(ua(M\cap T))
\text{ for }u\in K,a\in A.$$ 
Therefore the incidence relation 
amounts to the coincidence of the Riemannian horocycles 
and the inclusion incidence relation 
of the two totally geodesic submanifolds 
which are tangent to the coincident Riemannian horocycle 
of the composite distance $\log a$. 
\par 
Let $R_M$ or $R_M^*$ be the Radon transform in (1.0.1) 
or its dual in (1.0.2) 
for double fibrations of compact symmetric spaces(See Proposition 5-3(1)) 
$$M/(M\cap S)\leftarrow M/(M\cap S\cap T)
\rightarrow M/(M\cap T).$$ 
Then the Radon transform $R_M$ 
maps $C^\infty(M/(M\cap S))$ into $C^\infty(M/(M\cap T))$. 
Here $M/(M\cap S)$ or $M/(M\cap T)$ can be identified 
with the space of all the totally geodesic submanifolds 
which are tangent to the Riemannian horocycle 
at its composite distance point. 
\par 
Then the Radon transform $R_M$ is 
the restricted Radon transform $R_K$ on $M/(M\cap S)$. 
Here we have 
$$(R_K(f))(um(M\cap T))=R_{M}(f)(um(M\cap T))
\text{ for }u\in K,m\in M.$$ 
Therefore the incidence relation 
amounts to the inclusion incidence relation 
of the two totally geodesic submanifolds 
which are tangent to the Riemannian horocycle 
at its composite distance point. 
\par 
We put 
$$X_{ij}=e_{ji}-e_{ij}\in\mathfrak s\mathfrak o(m,n) 
\text{ for }1 \leq i<j \leq m,m+1\leq i<j\leq m+n$$ 
and 
$$Y_{ij}=e_{ji}+e_{ij}\in\mathfrak s\mathfrak o(m,n)
\text{ for }1 \leq i \leq m,m+1\leq j\leq m+n$$ 
where we put  
$$e_{ij}=(\delta_{pi}\dot\delta_{qj})_{1\leq p,q\leq m+n}
\text{ for }1\leq i,j\leq m+n.$$  
We put as follows. 
$$\mathfrak a(0)=\bigoplus\limits_{1\leq i\leq n}
\mathbb R Y_{i,m+n+1-i}
\text{ and }A(0)=\exp(\mathfrak a(0)).$$ 
$$Z_{i,j}^\pm=Y_{i,m+n+1-j}\mp Y_{j,m+n+1-i}\pm X_{i,j}+X_{m+n+1-j,m+n+1-i}
\text{ for }1\leq i<j\leq n.$$ 
$$W_{i,j}=Y_{i,j}+X_{i,m+n+1-j}
\text{ for }n+1\leq i\leq m,m+1\leq j\leq m+n.$$ 
For $1\leq l\leq n$ we have 
$$[Y_{l,m+n+1-l},Z_{i,j}^\pm]=(\delta_{li}\pm\delta_{lj})Z_{i,j}^\pm
\text{ for }1\leq i<j\leq n$$
and 
$$[Y_{l,m+n+1-l},W_{i,j}]=\delta_{m+n+1-l,j}W_{i,j}
\text{ for }n+1\leq i\leq m,m+1\leq j\leq m+n.$$ 
For $1\leq i\leq n$ we define $e_i\in\mathfrak a(0)^*$ by 
$$e_i(Y_{j,m+n+1-j})=\delta_{ij}\text{ for }1\leq j\leq n.$$  
We choose a set of positive roots $\Sigma^+$ by 
$$\Sigma^+=\{e_i\pm e_j\bigm\vert 1\leq i<j\leq n\}
\cup\{e_k\bigm\vert 1\leq k\leq n.\}.$$ 
Then we have 
$$\mathfrak g(0)_{e_i\pm e_j}=\mathbb R\dot Z_{i,j}^\pm
\text{ for }1\leq i<j\leq n $$ 
and 
$$\mathfrak g(0)_{e_k}=\bigoplus\limits_{n+1\leq i\leq m}
\mathbb R W_{i,m+n+1-k}
\text{ for }1\leq k\leq n.$$
We put 
$$\mathfrak n(0)=(\bigoplus\limits_{1\leq i<j\leq n}
(\mathfrak g_{e_i\pm e_j}))
\oplus 
(\bigoplus\limits_{1\leq k\leq n}\mathfrak g_{e_k})
\text{ and }
N(0)=\exp(\mathfrak n(0)).$$  
Let $M(0)$ be the centralizer of $\mathfrak a(0)$ in $K$, 
and $\mathfrak m(0)$ be the Lie algebra of $M(0)$.
We define a subgroup $M(0)_0$ of $M(0)$ by  
\begin{equation*}    
M(0)_0=\{\begin{pmatrix} 1_n & 0 & 0 
\\ 
0 & u & 0  
\\ 
0 & 0 & 1_n
\end{pmatrix} 
\bigm\vert 
u\in SO(m-n)\}.
\end{equation*}
We also define an abelian subgroup $M(0)_0^+$ of $M(0)$ by 
\begin{equation*}  M(0)_0^+ 
=\{\begin{pmatrix} diag(a_1,\cdots,a_n)  & 0 & 0 
\\ 
0 &  1_{m-n} & 0  
\\ 
0 & 0 & diag(a_n,\cdots,a_1)
\end{pmatrix}
\bigm\vert 
a_i=\pm 1(1\leq i\leq n),\prod\limits_{1\leq i\leq n}a_i=1\}.
\end{equation*} 
Then we have the following commuting products(See Lemma 2-4 in [26]). 
\begin{equation}\tag{1.2}  
\begin{split}
M(0)=M(0)^+_0\cdot M(0)_0, 
\\
M(0)\cap S=M(0)^+_0\cdot(M(0)_0\cap S),
\\
M(0)\cap T=M(0)^+_0\cdot(M(0)_0\cap T), 
\\
\text{ and }M(0)\cap S\cap T
=M(0)^+_0\cdot(M(0)_0\cap S\cap T).
\end{split}
\end{equation}

Here the compact homogeneous space 
$M(0)/(M(0)\cap S)$ or $M(0)/(M(0)\cap T)$ is 
diffeomorphic to $M(0)_0/(M(0)_0\cap S)$ or $M(0)_0/(M(0)_0\cap T)$ 
and to the compact Grassmann manifold 
$$SO(m-n)/S(O(m-n-r_1)\times O(r_1))
\text{ or }SO(m-n)/S(O(m-n-r_2)\times O(r_2).$$ 
Moreover, the Radon transform $R_{M}$ and its dual $R_{M}$ 
can be identified with the Radon transform 
and its dual on compact Grassmann manifolds in 
[8,13-15,29,34,37,39,40,45] 
(See Lemma 2-5 in [26] and Proposition 2-1). 
\par  
Here the horocycle space $\mathcal H_{r_1}$ or 
$\mathcal H_{r_2}$ 
is the $G$-translation of $(M(0)\cap S)N$ or  $(M(0)\cap T)N$ 
where $\mathcal H_{r_1}$ or $\mathcal H_{r_2}$ 
can be identified with the homogeneous space 
$G/((M(0)\cap H^{(r_1)})N(0))$ or $G/((M(0)\cap H^{(r_2)})N(0))$ (See Proposition 2-2). 
Under this identification the Radon transform $R$ and its dual $R^*$ 
can be identified with the Radon transform in (1.0.1) 
and its dual in (1.0.2) for the double fibrations 
(See Proposition 2-3) 
$$G/((M(0)\cap T)N)\leftarrow G/((M(0)\cap T\cap S)N)
\rightarrow G/((M(0)\cap S)N).$$ 
\par 
The Riemannian symmetric space $G/K$ 
can be identified with the bounded domain 
$$\mathcal D=\{W\in M_{m,n}\bigm\vert I_n-W\dot W^*>>0\},$$ 
where $X>>0$ for $X\in M_{n\times n}$ implies that 
$X$ is a positive definite symmetric $n\times n$ matrix
(See Proposition 13.2.1 in [12]). 
We put 
\begin{equation*}
o=\begin{pmatrix} 0 \\ 0\end{pmatrix}\in M_{m,n}.
\end{equation*}
We have  
$$g W=(X_{11}W+X_{12})(X_{21}W+X_{22})^{-1}
\text{ for }W\in\mathcal D,$$ 
where we put 
\begin{equation*}
g=\begin{pmatrix} X_{11} & X_{12} \\ X_{21} & X_{22}\end{pmatrix}
\in G,
X_{11}\in M_{m\times m},X_{12}\in M_{m\times n},
X_{21}\in M_{n\times m},X_{22}\in M_{n\times n}.
\end{equation*}
We put 
$$E_{ij}=(\delta_{pi}\dot\delta_{qj})_{1\leq p\leq m,1\leq q\leq n}
\text{ for }1\leq i\leq m,1\leq j\leq n.$$  
Then we have 
$$A(0) o=\{\sum\limits_{1\leq i\leq n} t_i E_{i,n+1-i}\in M_{m,n}
\bigm\vert -1<t_i<1\text{ for }1\leq i\leq n\}.$$ 
For $0\leq r\leq m-n$ 
we define the totally geodesic submanifold $D(r)$ of $D$ by 
\begin{equation*}\mathcal D(r)
=\{\begin{pmatrix} 0 \\  W_0\end{pmatrix}\in\mathcal D
\bigm\vert W_0\in M_{r\times n}\}.
\end{equation*} 
The semisimple symmetric space $G/H^{(r)}$ 
can be identified with the all the $G$-translations of $D(r)$.  
The horocycle $u a N(0) H^{(r)} K$ where $u\in K$ and $a\in A(0)$ 
can be identified with $(u a N(0) o,u a \mathcal D(r))$ 
where $u a N(0) o$ is a Riemannian horocycle 
and  $u a \mathcal D(r)$ is a totally geodesic submanifold  $u a \mathcal D(r)$ 
which is tangent to the Riemannian horocycle $u a N(0) o$ 
at its composite distance point $u a o$. 
\par 
Hereafter we assume the setting of the identification 
of $\mathcal H_{r_1}$ or $\mathcal H_{r_2}$ 
with $G/((M\cap H^{(r_1)})N)$ or $G/((M\cap H^{(r_2)})N)$
and the identification of $R$ and $R^*$ 
with the Radon transform in (1.0.1) 
and its dual in (1.0.2) for double fibrations 
$$G/((M\cap H^{(r_1)})N)\leftarrow   
G/((M\cap H^{(r_1)}\cap H^{(r_2)})N)\rightarrow G/((M\cap H^{(r_2)})N).$$
\par 
Let $L$ be a closed subgroup of $M$. 
We define the Schwartz space $\mathcal S(G/(LN))$ 
on $G/(LN))$ as follows. 
\begin{align*} & \mathcal S(G/(LN)))=\{
f\in C^\infty(G/(LN)))
\bigm\vert 
\\& 
\sup\limits_{u\in K,a\in A}\vert (X.f)(ua(LN)) \vert (1+\vert\log a\vert^n)
<\infty\text{ for }X\in U(\mathfrak g),n\in\mathbb Z^+_{\geq 0}\}.
\end{align*} 
We define for $d>0$, 
\begin{equation*} C^\infty_d(G/(LN))
=\{f\in C^\infty(G/((LN)))\bigm\vert 
f(ua((LN)))=0
\text{ for }u\in K,a\in A\text{ such that }
\vert\log A\vert> d\},
\end{equation*} 
where we denote by $\vert *\vert$ 
the norm on $\mathfrak a$ 
induced by the Killing form. 
\par 
We have 
$$C^\infty_c(G/((LN)))
=\bigcup\limits_{d>0}C^\infty_d(G/((LN))).$$ 
We define for $\lambda\in\mathfrak a^*_\mathbb C$, 
\begin{equation*}  C^\infty_\lambda(G/(LN)))
=\{f\in C^\infty(G/N)\bigm\vert 
f(ga((LN)))=f(g((LN)))a^{\lambda}
\text{ for }g\in G,a\in A\}.
\end{equation*} 
Then the space $C^\infty_\lambda(G/((LN))$ 
can be regarded as the space 
of all the $C^\infty$-sections of the line bundle  
of the homogeneous space $G/((LAN)$ 
with respect to $\lambda\in\mathfrak a_{\mathbb C}$. 
In view of the Iwasawa decomposition $G=KAN$, 
$C^\infty_\lambda(G/((LN))$ 
can be identified with $C^\infty(K/L)$ 
by the restriction from $G$ to $K$. 
In view of this identification 
the Radon transform $R$ which maps 
$C^\infty_\lambda(G/((M\cap S)N))$ 
into $C^\infty_\lambda(G/((M\cap T)N))$  
can be identified with the Radon transform $R_K$
for double fibrations of compact homogeneous spaces 
$K/(M\cap S)\leftarrow K/(M\cap S\cap T)
\rightarrow K/(M\cap T)$ 
which maps $C^\infty(K/(M\cap S))$ 
into $C^\infty(K/(M\cap T))$ 
for any fixed $\lambda\in\mathfrak a_{\mathbb C}$
(See Proposition 5-1(2))). 
Then we have the following diagram is commutative 
for $\lambda\in\mathfrak a^*_\mathbb C$. 
\[\begin{CD} 
  C^\infty_\lambda(G/((M\cap S)N)) @>  R >> C^\infty_\lambda(G/((M\cap T)N)) \\
@V\  restriction  VV   @VV\ restriction  V\\
   C^\infty(K/(K\cap S))  @>>  R_K    >  C^\infty(K/(K\cap T))
\end{CD}\tag{1.3}\] 
We put for $1\leq i<j\leq m+n$, 
\begin{equation*}{\overline X}_{ij}   
=\begin{cases}  X_{ij}  &  \text{ if } 1\leq i<j\leq m\text{ or }m+1\leq i<j\leq m+n\\ 
        -\sqrt{-1}\dot Y_{ij}  
        & \text{ if }1\leq i\leq m,m+1\leq j\leq m+n. 
\end{cases}\tag{1.4} 
\end{equation*}
We denote by $T_{s}(N)$ 
the collection of all increasing    
sequences $I=(i_1,\cdots,i_s)$ of length $s$ in 
$\{1,\cdots,N\}$ 
where $N,s\in \mathbb Z^+_{\geq 0}$. 
We denote by $P(I)$ the set of all the permutations of $I$. 
We denote by $\epsilon(\sigma)$ 
the sign of $\sigma\in P(I)$. 
We denote by $P(I)$ the set of all the permutations of $I$. 
We denote by $\epsilon(\sigma)$ 
the sign of $\sigma\in P(I)$. 
For $p\geq 1$ and $I\in T_{2p}(N)$  we put 
\begin{equation*}  
P_+(I) =\{\sigma\in P(I)\bigm\vert 
\sigma(i_{2a-1})<\sigma(i_{2a+1})\text{ for }1\leq a\leq p-1, 
\text{ and }
\sigma(i_{2a-1})<\sigma(i_{2a})\text{ for }1\leq a\leq p \}.
\end{equation*} 
Assume that $p\geq 1,2p\leq m+n$. 
We put for $I=\{i_1,\cdots,i_{2p}\}\in T_{2p}(m+n)$, 
\begin{equation*}W_{I} 
=\sum_{\sigma\in P_+(I)} 
\epsilon(\sigma){\overline X}_{\sigma(i_1)\sigma(i_2)}
\cdots {\overline X}_{ \sigma(i_{2p-1}) \sigma(i_{2p}) }
\in U(\mathfrak g_{\mathbb C})\tag{1.5}
\end{equation*}
and 
\begin{equation*}W_p=\sum\limits_{J\in T_{2p}(m+n)}{W}_{J}^2
\in U(\mathfrak g_{\mathbb C}).\tag{1.6}
\end{equation*}
\subsection{Main Results}
Let $\tilde r_1=\min(m-r_1,r_1+n)$ and 
$\tilde r_2=\min(m-r_2,r_2+n)$. 
Then $\tilde r_1$ or $\tilde r_2$ 
is the split rank 
of the semisimple symmetric space $G/S$ or $G/T$. 
\begin{theorem} 
[Injectivity(Cf.Proposition 5-5)]
Assume that $n\geq 1; 0\leq r_1,r_2 \leq m-n ; r_1\not= r_2$. 
Let $\Gamma=C^\infty,\mathcal S,C^\infty_c,C^\infty_d$ with $d>0$ 
or $C^\infty_\lambda$ 
with $\lambda\in\mathfrak a^*_{\mathbb C}$. 
Then the following conditions are equivalent.  
\par 
$(1)$ $r_1<r_2,r_1+r_2\leq m-n\text{ or }r_1>r_2,r_1+r_2\geq m-n$, 
\par 
$(2)$ $\dim G/((M\cap S)N))\leq \dim G/((M\cap T)N))$ 
\par 
$(3)$ $\tilde{r_1}\leq \tilde{r_2}$.  
\par 
$(4)$ The Radon transform $R$ maps $\Gamma(G/((M\cap S)N))$ 
injectively into $\Gamma(G/((M\cap T)N))$.  
\end{theorem} 
In the setting of Theorem 1
we also obtain the inversion formula of the Radon transform $R$ 
on $C^\infty(G/((M\cap S)N))$(See Proposition 5-7(3)). 
Our inversion formluas are nonlocal(See [5]) 
and expressed as infinite series. 
By the inversion formula we can prove the following support theorem. 
\begin{theorem} [Support Theorem(Cf.Proposition 5-7(4)))]
Assume the setting of Theorem 1.  
Suppose $f\in C^\infty(G/((M\cap S)N))$. 
Then we have $f\in C^\infty_c(G/((M\cap S)N))$ 
if and only if $R(f)\in  C^\infty_c(G/((M\cap T)N))$. 
Moreover, we have $f\in C^\infty_d(G/((M\cap S)N))$ 
if and only if $R(f)\in  C^\infty_d(G/((M\cap T)N))$ 
for $d>0$. 
\end{theorem} 
\begin{theorem}[Surjectivity(Cf.Proposition 6-1)] 
Assume that $n\geq 1; 0\leq  r_1,r_2\leq m-n ; r_1\not= r_2$. 
Let $\Gamma=C^\infty,\mathcal S,C^\infty_c,C^\infty_d$ with $d>0$ 
or $C^\infty_\lambda$ with $\lambda\in\mathfrak a^*_{\mathbb C}$. 
Then the following conditions are equivalent.  
\par 
$(1)$ $r_1<r_2,r_1+r_2\geq m-n\text{ or }r_1>r_2,r_1+r_2\leq m-n$, 
\par 
$(2)$ $\dim G/((M\cap S)N))\geq \dim G/((M\cap T)N))$ 
\par 
$(3)$ $\tilde{r_1}\geq \tilde{r_2}$.  
\par 
$(4)$ The Radon transform $R$ maps $\Gamma(G/((M\cap S)N))$ 
surjectively onto $\Gamma(G/((M\cap T)N))$.  
\end{theorem} 
\begin{theorem}[Bijectivity(Cf.Proposition 5-5,6-1)] 
Assume that $n\geq 1; 0\leq r_1,r_2\leq m-n; r_1\not= r_2$. 
Let $\Gamma=C^\infty,\mathcal S,C^\infty_c,C^\infty_d$ with $d>0$ 
or $C^\infty_\lambda$ with 
$\lambda\in\mathfrak a^*_{\mathbb C}$. 
Then the following conditions are equivalent.  
\par 
$(1)$ $r_1+r_2=m-n$. 
\par 
$(2)$ $\dim G/((M\cap S)N))=\dim G/((M\cap T)N))$ 
\par 
$(3)$ $\tilde{r_1}=\tilde{r_2}$.  
\par 
$(4)$ The Radon transform $R$ maps $\Gamma(G/((M\cap S)N))$ 
bijectively onto $\Gamma(G/((M\cap T)N))$. 
\end{theorem} 
\begin{theorem}[(Cf.Proposition 4-1)] 
Assume that $p\geq 1$ and $2p\leq m+n$. Then we have 
$$Ad(G)(\bigoplus\limits_{I\in \in T_{2p}(m+n)}
\mathbb C\dot W_I)
\subset(\bigoplus\limits_{I\in \in T_{2p}(m+n)}
\mathbb C\dot W_I)
\text{ and }W_p\in\mathcal Z(\mathfrak g).$$
\end{theorem}
Assume that $n\geq 1; 0\leq r_1,r_2\leq m-n; r_1\not= r_2$. 
The condition that $r_1<r_2,r_1+r_2<m-n\text{ or }r_1>r_2,r_1+r_2>m-n$ 
is equivalent to the condition of $\tilde{r_1}<\tilde{r_2}$ 
or to the condition of the overdeterminedness: 
$\dim G/((M\cap S)N)<\dim G/((M\cap T)N)$(See  [38]). 
In this case our inversion formula is overdetermined(See Proposition 5-7(3)). 
Moreover, we have the following range theorems 
by the invariant system of differential operators $(1.5)$ 
or the invariant differential operator $(1.6)$. 
\begin{theorem}
[Main Theorem(Range Theorem(Cf.Theorem 6-5))] 
Assume that $n\geq 1; 0\leq r_1,r_2 \leq m-n ; r_1\not= r_2$. 
Suppose that $r_1<r_2,r_1+r_2<m-n\text{ or }r_1>r_2,r_1+r_2>m-n$ 
which is equivalent to $\dim G/((M\cap S)N)<\dim G/((M\cap T)N)$ 
or to $\tilde{r_1}<\tilde{r_2}$.  
Then the following (I) and (II) hold. 
\par 
(I) The following $(1)$ and $(2)$ hold 
for $\Gamma=\mathcal S,C^\infty_c$ or $C^\infty_d$ with $d>0$. 
\par 
(II) The following $(1)$ holds  
for $\Gamma=C^\infty_\lambda$ with generic $\lambda\in\mathfrak a^*_{\mathbb C}$. 
\begin{equation} R(\Gamma(G/((M\cap S)N)))
=\{f\in \Gamma(G/((M\cap T)N))
\bigm\vert W_{I}.f\equiv 0\text{ for }I\in T_{2(\tilde r_1+1)}(m+n)\}.\tag{1} \end{equation}
\begin{equation}R(\Gamma(G/((M\cap S)N)))
=\{f\in \Gamma(G/((M\cap T)N))\bigm\vert 
W_{\tilde r_1+1}.f\equiv 0\}.\tag{2}\end{equation}
\end{theorem} 
In view of (1.3) the Main Theorem (II) yields the solution 
to the problem (A) and (C) for the Radon transform $R_K$ 
for double fibrations of compact homogeneous spaces 
$K/(M\cap S)\leftarrow K/(M\cap S\cap T)
\rightarrow K/(M\cap T)$. 
\par 
Our approach to prove Theorem 1,2,3,4 and the Main Theorem 
is based on the Fourier transform on $M$ in (4.5) 
and the Fourier transform on $A$ in (6.2) 
and their relations with the Radon transform $R$ 
in Proposition 5-4 and Proposition 6-2. 
\subsection{Historical Notes and Some Comments}
As for the Radon transform in (1.0.1) for double fibrations 
of semisimple symmetric spaces 
$G/S\leftarrow G/(S\cap T)\rightarrow G/T$, 
the solutions to the problems (A) and (D) are obtained 
in [24] and [26] when $r_1<r_2$ 
and the solutions for the problems (A) and (B) 
are obtained in [46] when $r_1=0$. 
In Theorem 2-12(2) in [26] the projection slice theorem 
which relates the Radon transform in (1.0.1) 
for double fibrations $G/S\leftarrow G/(S\cap T)\rightarrow G/T$ 
and the Fourier transform on $G/S$ and $G/T$ in [1-3] 
can be stated in terms of the Radom transform $R$ in [27]. 
\par 
In the setting of $n=1$, 
the range characterizations of the Radon transform 
on the real hyperbolic Grassmann manifolds 
by the invariant system of differential operators (1.5) 
are given in Theorem 5-1 in [25] as the solutions to the problem (C).  
Moreover, the range characterizations 
by the invariant differential operators (1.6) in this setting 
are given in [27] where Proposition 4-2,Proposition 6-4(1),
Theorem 5-1 in [25] and the Fourier transform in [1-3] are used.  
\par 
In the setting of $(n,r_1)=(1,0)$, 
the range characterizations 
of the totally geodesic Radon transform 
on the real hyperbolic spaces 
by the invariant system of differential operators 
(1.5) or the invariant differential operators 
(1.6) for $p=2$ in Theorem A,D,Corollary 11.4 in [23] 
are given as the solutions to the problem (C).  
\par 
The horocycle Radon transform on Riemannian symmetric space has been studied in [19] and played an important role in the study of harmonic analysis 
on the Riemannian symmetric space. 
\par 
The solution to the Problem (A) and (B) 
for the horocycle Radon transform on real hyperbolic space 
with respect to the lower dimensional horocycle(horosphere) has been studied in [4]. 
\par 
The solution to the Problem (A) and (D) 
for the horocycle Radon transform on semismple symmetric space has been studied in [31].  
\par 
The solutions 
to the problem (A) and (B) for the Radon transform 
for double fibrations of compact Grassmann manifolds 
as well as the non-symmetric compact homogeneous spaces 
such as Stiefel manifolds are found in [14,15,29,37,39,40,44,45]. 
\par  
The solutions to the problem (C) 
for the Radon transform for double fibrations 
of compact Grassmann manifolds have been given 
in [7,8,13,29,34,40] etc 
by the invariant system of differential operators 
or the invariant differential operators. 
\par 
In [34] the solutions to the problem (C) 
has been obtained by the identifications 
of the function spaces on the Grassmann manifold 
with the spaces of all the sections of the line bundle 
on the flag manifolds in a similar way to (1.3).    
\par 
The solutions to the problem (C) 
for the Radon transform for double fibrations 
of Euclidean symmetric spaces have been given 
in [6,8-11,28,32,35] etc 
by the invariant system of differential operators 
or the invariant differential operators.

\section{Preliminaries.} 
Throughout this article, 
we denote by $d(*)$ the bi-invariant measure on $L$ 
and denote by $d(*)_{L_0}$ 
the left $L$-invariant measure on $L/L_0$ 
where $L$ is a unimodular Lie group 
and $L_0$ is a closed unimodular subgroup in $L$. 
We normalize the invariant measure $d(*)$ or $d(*)_{L_0}$ 
so that the total volume on $L$ or $L/L_0$ 
with respect to $d(*)$ or $d(*)_{L_0}$ is equal to $1$ 
whenever $L$ or $L/L_0$ is compact respectively. 
We normalize the bi-invariant measure $dg$ on $L$, 
the bi-invariant measure $dh$ on $H$, 
and the invariant measure $dg_{L_0}$ on $L/L_0$ 
(See Chapter II,\S 2 in [20]) so that   
\begin{equation}\int\nolimits_{L}F(g)dg
=\int\nolimits_{L/L_0}
(\int\nolimits_{L_0}
F(g h)dh) dg_{ L_0 }\text{ for }F\in C^\infty_c(L)
\tag{2.0}\end{equation}  
Here we have $K=H^{(0)}$. 
The homogeneous space $G/K$ is a Riemannian symmetric space of rank $n$ 
and $G/H^{(r)}$ is a semisimple 
symmetric space of split rank $n$ if $m-n\geq r$ 
and of split rank $m-r$ if $m-n<r$(See p.113 in [41]). 
\par 
Let $U$ be a real Lie group and $H$ be a closed subgroup of $U$. 
Let $\mathfrak u$ be the Lie algebra of $U$.  
We define the left action of $U(\mathfrak u)$ for $f\in C^\infty(U/H)$ by 
\begin{equation}((X_1\dot X_2\cdots X_l).f)(u H)
=\frac{\partial^l}{\partial t_1\partial t_2\cdots
\partial t_l}\bigm\vert_{t_1=\cdots=t_l=0}\cdot 
f(e^{-t_l\cdot X_l}\cdots e^{-t\cdot X_2}\cdot 
e^{-t_1\cdot X_1}\cdot u H)
\tag{2.1}\end{equation} 
for $u\in U$ and $X_1,X_2,\cdots,X_l\in\mathfrak u$. 
\par 
Let $\mathfrak l$ is a Lie algebra of $L$. 
Then we have for $f\in  C^\infty(U/L)$,$u\in U$, $X\in U(\mathfrak u)$ 
and $X_0\in\mathfrak l$, 
\begin{equation}((Ad(u)X_0)X.f)(u L )=0
\text{ for }u\in U.\tag{2.2}\end{equation} 
\par 
\textbf{Proposition 2-1}. 
\par 
Assume $n\geq 1; 0\leq r_1,r_2\leq m-n$. 
Let $S=H^{(r_1)}$ and $T=H^{(r_2)}$. 
There exists some $u_M\in K\cap S\cap T$ such that 
$$A=u_M A(0) u_M^{-1},M=u_M M(0) u_M^{-1}
\text{ and }N=u_M N(0) u_M^{-1}.$$
\begin{Proof}
We may assume that $0\leq r_1\leq r_2\leq m-n$  
without loss of generality. 
Let $G_1$ be the reductive Lie subgroup of G as follows. 
\begin{equation*}G_1
=\{\begin{pmatrix} g_1 & 0  & g_2 
\\ 
0 & u & 0  
\\ 
g_3 & 0 & g_4
\end{pmatrix}\in G 
\bigm\vert 
\begin{pmatrix} g_1 & g_2 
\\ 
g_3 & g_4 
\end{pmatrix}\in O(m-r_2,n),u\in (O(r_2-r_1)\times O(r_1))\}. 
\end{equation*}
Let $\mathfrak g_1$ be the Lie algebra of $G_1$. 
Then we have $\mathfrak p\cap\mathfrak q_{\mathfrak h^{(r_2)}}=\mathfrak p\cap\mathfrak g_1$. 
It follows from Theorem 6.51. in [30]
that the maximal abelian subspace $\mathfrak a$ 
and $\mathfrak a(0)$ 
in $\mathfrak p\cap\mathfrak q_{\mathfrak h^{(r_2)}}$ 
are conjugate by some element in $G_1\cap K$. 
Since $G_1\cap K\subset K\cap S\cap T$,   
there exists some $u_0\in K\cap S\cap T$ such that $\mathfrak a(0)=Ad(u_0)\mathfrak a$.
Hence we have $A(0)=u_0 A u_0^{-1}$. 
Since $M(0)$ or $M$ is the centralizer of $\mathfrak a(0)$ or $\mathfrak a$, 
we have $M(0)=u_0 M u_0^{-1}$. 
Let $W$ be the factor group $N_K(\mathfrak a(0))/Z_K(\mathfrak a(0))$ 
and $W_1$ be the factor group 
$N_{K\cap G_1}(\mathfrak a(0))/Z_{K\cap G_1}(\mathfrak a(0))$. 
Then we have $W=W_1$, since we have $K_1 \subset G_1\cap K$, where we put 
\begin{equation*}   
K_1=\{\begin{pmatrix} u & 0  & 0 & 0
\\ 
0 & \det(u)^{-1} & 0 & 0 
\\ 
0 & 0 & 1_{m-n-1} & 0
\\ 
0 & 0 & 0 & u'
\end{pmatrix}\in SO(m)\times SO(n) 
\bigm\vert 
u\in O(n),u'\in SO(n)\}.
\end{equation*} 
Here we can easily verify that 
there exists some $u(i)\in K_1(1\leq i\leq n)$ 
and $u_{ij}\in K_1(1\leq i<j\leq n)$ such that the followings hold.
\begin{align*} & 
Ad(u(i))Y_{i,m+n+1-i}=-Y_{i,m+n+1-i}
\text{ and }Ad(u(i))Y_{j,m+n+1-i}=Y_{j,m+n+1-j}
\text{ for }1\leq j\leq n,j\not=i. 
\\& 
Ad(u_{ij})Y_{i,m+n+1-i}=Y_{j,m+n+1-j},
Ad(u_{ij})Y_{j,m+n+1-j}=Y_{i,m+n+1-i}
\\& 
\text{ and }Ad(u_{ij})Y_{k,m+n+1-k}=Y_{k,m+n+1-k}
\text{ for }1\leq k\leq n,k\not=i,j.
\end{align*} 
Then it follows from Corollary 6.55. in $[30]$  
that there exists some 
$u_1\in N_{K\cap G_1}(\mathfrak a(0))
\subset K\cap S\cap T$ 
such that 
$$N(0)=u_1(u_0 N u_0^{-1})u_1^{-1}=(u_1 u_0)N(u_1 u_0)^{-1}.$$ 
We put $u_M=(u_0 u_1)^{-1}\in K\cap S\cap T$. 
Then the assertion follows. \qed 
\end{Proof}
We define 
$$M_0=u_M\cdot M(0)_0\cdot u_M^{-1}
\text{ and }
M_0^+=u_M\cdot M(0)_0^+\cdot u_M^{-1}.$$ 
Then $M_0^+$ is abelian. 
By (1.2) we have the following commuting products. 
\begin{equation}  
\begin{split}
M=M^+_0\cdot M_0, 
\\
M\cap S=M^+_0\cdot(M_0\cap S),
\\
M\cap T=M^+_0\cdot(M_0\cap T), 
\\
\text{ and }M\cap S\cap T
=M^+_0\cdot(M_0\cap S\cap T).
\end{split}
\tag{2.3}\end{equation}
In view of (1.2) and (2.3) and Proposition 2-1,
the compact homogeneous space $K/(M\cap S)$ 
or $K/(M\cap T)$ which is not necessarily symmetric 
is diffeomorphic to a quotient of  
$(O(m)/(O(m-n-r_i)\times O(r_i))\times O(n)$ 
for $i=1$ or $i=2$ respectively. 
\par 
\textbf{Proposition 2-2}  
\par 
Let $0\leq r\leq m-n$. Let $\mathcal H_r$ be 
the $G$-translation of $((M\cap H^{(r)})N) H^{(r)}$ in $G/H^{(r)}$.  
Then $\mathcal H_r$ does not depend on the choice of $A$ and $N$
\begin{Proof}
In the setting of Propsition 2-1 we have 
\begin{equation}((M\cap H^{(r)})N) H^{(r)}=u_M (M(0)\cap H^{(r)})N(0)H^{(r)}.\tag{2.2.1}
\end{equation}
Then it follows that $\mathcal H_r$ 
is the $G$-translation of $(M(0)\cap H^{(r)})N(0)H^{(r)}$. 
This yields the assertion. \qed 
\end{Proof}
\textbf{Proposition 2-3} 
\par 
Assume $n\geq 1; 0\leq r_1,r_2\leq m-n; r_1\not= r_2$. 
Let $S=H^{(r_1)}$ and $T=H^{(r_2)}$. 
We define the Radon transform $R$ which maps $C^\infty_c(\mathcal H_{r_1})$ 
into $C^\infty(\mathcal H_{r_2})$ and the dual Radon transform $R^*$ which maps 
$C^\infty_c(\mathcal H_{r_2})$ into $C^\infty(\mathcal H_{r_1})$ 
as the Radon transform in (1.0.1) 
and its dual in (1.0.2) for the double fibrations 
$$G/((M\cap S)N)\leftarrow G/((M\cap S\cap T)N)
\rightarrow G/((M\cap T)N)$$ 
under the identification of $\mathcal H_{r_1}$ or $\mathcal H_{r_2}$ 
with $G/((M\cap H^{r_1})N)$ or $G/((M\cap H^{r_2})N)$. 
Then $R$ or $R^*$ does not depend on the choice of $A$ and $N$.  
\begin{Proof}
We have only to prove the assertion for $R$, since the assertion for $R^*$ 
can be proved in a similar way to the assertion for $R$. 
Let $R(0)$ be the Radon transform in (1.0.1) 
and its dual $R(0)^*$ in (1.0.2) for the double fibrations 
$$G/((M(0)\cap S)N(0))\leftarrow G/((M(0)\cap S\cap T)N(0))
\rightarrow G/((M(0)\cap T)N(0)).$$ 
Let $\xi_T\in \mathcal H_{r_2}$ be $\xi_T=(M\cap T)N$  
and $\xi_S\in \mathcal H_{r_1}$ be $\xi_S=(M\cap S)N$. 
By the identification 
of $\mathcal H_{r_1}=G/((M\cap S)N)$ and $\mathcal H_{r_2}=G/((M\cap T)N)$
and (2.2.1), we have 
$$\xi_T=(M\cap T)N=u_M(M(0)\cap T)N(0)  
\text{ and }\xi_S=(M\cap S)N=u_M(M(0)\cap S)N(0).$$ 
We have for $f\in C^\infty_c(G/(M\cap S)N)$ and $g\in G$,  
\begin{align*} & (R f)(g \xi_T))=(R f)(g ((M\cap T)N))
\\& 
=\int\nolimits_{((M\cap T)N)/((M\cap S)N)}
f(g\cdot h ((M\cap S)N))dh_{((M\cap S)N)\cap((M\cap T)N)}
\\& 
=\int\nolimits_{(u_M((M(0)\cap T)N(0))u_M^{-1})/u_M((M(0)\cap S)N(0))u_M^{-1}}
f(g\cdot h \xi_S)
dh_{u_M((M(0)\cap S)N(0))\cap((M(0)\cap T)N(0))u_M^{-1}}.
\end{align*} 
Here we have for $h\in(M\cap T)N$, 
$$g\cdot h \xi_S=g\cdot h\cdot u_M ((M(0)\cap S)N(0))
=g u_M u_M^{-1} h u_M ((M(0)\cap S)N(0)).$$ 
We put $h(0)=u_M h u_M^{-1}\in((M(0)\cap T)N(0))$. 
Then we have $g\cdot h \xi_S=g\cdot u_M h(0)((M(0)\cap S)N(0))$. 
Then we have 
\begin{align*} &  
(R f)(g \xi_T))
\\& 
=\int\nolimits_{((M(0)\cap T)N(0))/((M(0)\cap S)N(0))}
f(g u_M h(0)((M(0)\cap S)N(0))) 
dh(0)_{((M(0)\cap S)N(0))\cap((M(0)\cap T)N(0))}. 
\end{align*} 
Then we have for $g\in G$, 
$$(R f)(g \xi_T)=(R(0) f)(g u_M(M(0)\cap T)N))=(R(0) f)(g \xi_T).$$
This yields the assertion \qed 
\end{Proof}
\textbf{Proposition 2-4}  
\par 
Let $0\leq r\leq m-n$. Assume the setting of $(1.1)$. 
The map $i_r$ from $\mathcal H_r$ into $\mathcal H(r)$ 
is a diffeomorphism.  
\begin{Proof}
First we prove that the injectivity of $i_r$. 
Let $(u_0,a_0)\in K\times A$ and $(u_1,a_1)\in K\times A$ satisfy 
$$(u_0 a_0 M N, u_0 a_0 H^{(r)})=(u_1 a_1 M N,u_1 a_1 H^{(r)}).$$
Then by the uniqueness of the Iwasawa decomposition 
it follows from $u_0 a_0 M N=u_1 a_a M N$ that $a_0=a_1$ 
and that there exists some $m\in M$ such that $u_0 m= u_1$. 
Then it follows from $u_0 a_0 H^{(r)}=u_1 a_1 H^{(r)}$ that 
$$u_0 a_0  H^{(r)}=u_1 a_1 H^{(r)}=u_0 m a_0 H^{(r)}=u_0 a_0 m H^{(r)}.$$ 
Then we have $m H^{(r)}=H^{(r)}$. Hence we have $m\in M\cap H^{(r)}$.  
Therefore we have 
$$u_0 a_0 (M\cap H^{(r)})N= u_1 a_1 (M\cap H^{(r)}).$$ 
In view of the diffeomorphism between $\mathcal H_r$ and $K/(M\cap H)\times A$, 
this yields the injectivity of $i_r$. 
The surjectivity of $i_r$ follows from the definition of $i_r$. 
Then $\mathcal H(r)$ is diffeomorphic to $K/(M\cap H^{(r)})\times A$. 
Thus we have proved the assertion. \qed  
\end{Proof}
\textbf{Proposition 2-5}  
\par 
Let $0\leq r\leq m-n$. 
The totally geodesic submanifold $H^{(r)} K$ in $G/K$ 
is tangent to the Riemannian horocycle $NK$ at $e K$. 
\begin{Proof}
By Proposition 2-1 we have only to prove the assertion for $N=N(0)$. 
We can prove the assertion in a similar way to 
Lemma 1.2 in [19], Ch VI Exercise B.2 in $[16]$ 
and the proof in p571 in $[16]$. 
Let $\pi$ be the natural projection from $G$ onto $G/K$ 
and $d\pi$ be the differential map at $g=e$. 
The tangent space of $N(0)K$ at $eK$ is $d\pi(\mathfrak n(0))$. 
We have 
$$d\pi(\bigoplus\limits_{n+1\leq i\leq m,m+1\leq j\leq m+n}
\mathbb R Y_{i,j})
\subset 
d\pi(\bigoplus\limits_{n+1\leq i\leq m,m+1\leq j\leq m+n}
\mathbb R W_{i,j})
\subset d\pi(\mathfrak n(0)).$$
Here the tangent space of $H^{(r)}K$ at $eK$ 
is $d\pi(\mathfrak h^{(r)}\cap\mathfrak p)$. 
It follows from $0\leq r\leq m-n$ that 
$$\mathfrak h^{(r)}\cap\mathfrak p
\subset\bigoplus\limits_{n+1\leq i\leq m,m+1\leq j\leq m+n}
\mathbb R Y_{i,j}.$$  
Then we have 
$$d\pi(\mathfrak h^{(r)}\cap\mathfrak p)\subset 
d\pi(\mathfrak n(0)).$$ 
Thus we have proved the assertion. \qed 
\end{Proof}
\section{Pfaffian type elements 
for the Orthogonal Lie algebra.} 
In this section 
we investigate the properties of the Pfaffian type elements for $\mathfrak s\mathfrak o(N,\mathbb F)$ 
by the explicit calculations of the projection of these elements 
which lead to the calculations 
of the eigenvalues of the central elements of 
$U(\mathfrak s\mathfrak o(N,\mathbb F))$
(See Proposition 5.42. of [30]). 
\par 
Let $\mathbb F$ be $\mathbb R$ or $\mathbb C$ and $N$ 
be a positive integer such that $N\geq 3$. 
We put as follows. 
$$SO(N,\mathbb F)=\{u\in GL(N,\mathbb F) 
\bigm\vert u^{t}\cdot u=I_N,\det u=1\}.$$
$$\mathfrak s\mathfrak o(N,\mathbb F)=\{X\in \mathfrak g \mathfrak l(N,\mathbb F)
\bigm\vert X=-X^t,\text{tr }X=0\}.$$ 
$${\tilde E}_{ij}
=(\delta_{pi}\cdot\delta_{qj})_{1\leq p,q \leq N}
\in \mathfrak g\mathfrak l(N,\mathbb F) 
\text{ for }1\leq i,j\leq N.$$ 
$${\tilde X}_{ij}
={\tilde E}_{ji}-{\tilde E}_{ij}
\in \mathfrak s\mathfrak o(N,\mathbb F)
\text{ for }1\leq i<j\leq N.$$ 
Let $\tilde{\mathfrak g}
=\mathfrak s\mathfrak o(N,\mathbb C)$. 
We put 
$$\tilde{\mathfrak h}
=\bigoplus\limits_{1\leq i\leq [N/2]}
\mathbb C\tilde{X}_{i,N+1-i}.$$ 
Then $\tilde{\mathfrak h}$ is a maximal abelian subalgebra of $\tilde{\mathfrak g}$. 
Let $\tilde{\Sigma}$ be the set of all the roots. 
We put for $\alpha\in\tilde{\Sigma}$, 
$$\tilde{\mathfrak g}_\alpha=\{X\in\tilde{\mathfrak g}
\bigm\vert [H,X]=\alpha(H)
\text{ for }H\in\tilde{\mathfrak g}.\}.$$ 
We put as follows. 
$$\tilde{Z}_{i,j}^\pm
=\sqrt{-1}\tilde{X}_{i,N+1-j}
\mp\sqrt{-1}\tilde{X}_{j,N+1-i}
\pm\tilde{X}_{i,j}+\tilde X_{N+1-j,N+1-i}
\text{ for }1\leq i<j\leq N.$$ 
$$\tilde W_{i,j}=\sqrt{-1}\tilde X_{i,j}+\tilde X_{i,N+1-j}
\text{ for }1\leq i<j\leq N.$$ 
For $1\leq l\leq [N/2]$ we have 
$$[\sqrt{-1}\tilde X_{l,N+1-l},\tilde Z_{i,j}^\pm]
=(\delta_{li}\pm\delta_{lj})\tilde Z_{i,j}^\pm
\text{ for }1\leq i<j\leq N$$
$$\text{ and }[\sqrt{-1}\tilde X_{l,N+1-l},\tilde W_{i,j}]
=\delta_{N+1-l,j}\tilde W_{i,j}
\text{ for }1\leq i<j\leq N.$$ 
For $1\leq i\leq [N/2]$ 
we define $e_i\in\tilde{\mathfrak h}^*$ by 
$$\tilde e_i(\sqrt{-1}\tilde X_{j,N+1-j})
=\delta_{ij}\text{ for }1\leq j\leq [N/2].$$   
We choose a set of positive roots $\tilde\Sigma^+$ by 
$$\tilde{\Sigma}^+   
=\begin{cases} 
\{\tilde e_i\pm\tilde e_j\bigm\vert 1\leq i<j\leq [N/2]\}
\cup\{\tilde e_{[N/2]}\} 
        & \text{ if }N\text{ is odd} \\
\{\tilde e_i\pm\tilde e_j\bigm\vert 1\leq i<j\leq [N/2]\}
        & \text{ if }N\text{ is even}, 
\end{cases}$$ 
where we have 
$\tilde{\Sigma}^+=\{\tilde e_1\}$ if $N=3$. 
Then we have 
$$\tilde{\mathfrak g}_{\tilde e_i\pm\tilde e_j}
=\mathbb C\cdot\tilde Z_{i,j}^\pm
\text{ for }1\leq i<j\leq [N/2]$$ 
and if $N$ is odd we have 
$$\tilde{\mathfrak g}_{\tilde e_{[N/2]}}
=\mathbb C\tilde W_{[N/2]+1,N+1-[N/2]}.$$   
We define the Lie subalgebra 
$\tilde{\mathfrak n}(0)$ of $\tilde{\mathfrak g}$ by 
$$\tilde{\mathfrak n}(0)=
\bigoplus\limits_{\alpha\in\Sigma^+}
\tilde{\mathfrak g}_\alpha.$$  
If $N$ is odd, we put for $1\leq q\leq[N/2]$, 
$$\tilde{\mathfrak n}(0)_q
=(\bigoplus\limits_{q< l\leq [N/2]}
\mathbb C(\tilde Z^+_{q,l}-\tilde Z^-_{q,l}))
\oplus(\mathbb C\tilde W_{[N/2]+1,N+1-q})
\oplus 
(\bigoplus\limits_{q< l\leq [N/2]}
\mathbb C(\tilde Z^+_{q,l}+\tilde Z^-_{q,l})),$$ 
where we have 
$\tilde{\mathfrak n}(0)_{[N/2]}=\mathbb C\tilde  W_{[N/2]+1,N+1-[N/2]}$. 
\par 
If $N$ is even, we put for $1\leq q<[N/2]$, 
$$\tilde{\mathfrak n}(0)_q
=(\bigoplus\limits_{q< l\leq [N/2]}
\mathbb C(\tilde Z^+_{q,l}-\tilde Z^-_{q,l}))
\oplus 
(\bigoplus\limits_{q< l\leq [N/2]}
\mathbb C(\tilde Z^+_{q,l}+\tilde Z^-_{q,l})).$$ 
Then we have 
$$\tilde{\mathfrak n}(0)
=\begin{cases} 
\bigoplus\limits_{1\leq q\leq [N/2]}\tilde{\mathfrak n}(0)_q 
& \text{ if }n\text{ is odd} \\
\bigoplus\limits_{1\leq q<[N/2]}
\tilde{\mathfrak n}(0)_q 
& \text{ if }n\text{ is even}. 
\end{cases}$$

If $N$ is odd, we put for $1\leq q\leq[N/2]$, 
$$\tilde N^q_l   
=\begin{cases} (\tilde Z_{q,l}^--\tilde Z_{q,l}^+)/2 & 
\text{ for }q< l\leq [N/2] \\
        \tilde W_{[N/2]+1,N+1-q} & \text{ for }l=[N/2]+1 \\ 
(\tilde Z_{q,N+1-l}^++\tilde Z_{q,N+1-l}^-)/2  & 
\text{ for }[N/2]+1< l< N+1-q. 
\end{cases}$$ 

If $N$ is even, we put for $1\leq q<[N/2]$, 
$$\tilde N^q_l   
=\begin{cases} (\tilde Z_{q,l}^--\tilde Z_{q,l}^+)/2 & 
\text{ for }q< l\leq [N/2] \\
        (\tilde Z_{q,N+1-l}^++\tilde Z_{q,N+1-l}^-)/2  & 
\text{ for }[N/2]< l< N+1-q. 
\end{cases}$$ 

Then we have 
\begin{equation*}
{\tilde N}^q_{l}   
=\begin{cases}  \tilde X_{l,N+1-q}+\sqrt{-1}\cdot\tilde X_{q,l} 
          &  \text{ for }1\leq q\leq [N/2],q< l< N+1-q 
          \text{ if }N\text{ is odd }\\ 
          \tilde X_{l,N+1-q}+\sqrt{-1}\cdot\tilde X_{q,l}
        & \text{ for }1\leq q<[N/2],q< l< N+1-q\text{ if }N\text{ is even }. 
\end{cases}
\tag{3.1} 
\end{equation*}

Then we have 

$$\tilde{\mathfrak n(0)}_q   
=\begin{cases}  
\bigoplus\limits_{q< l\leq N+1-q}\mathbb C\tilde N^q_l
        &  \text{ for }1\leq q\leq [N/2] 
          \text{ if }N\text{ is odd }    \\ 
          \bigoplus\limits_{q< l\leq N+1-q}\mathbb C\tilde N^q_l
        & \text{ for }1\leq q<[N/2]\text{ if }N\text{ is even }. 
\end{cases}$$ 
We have 
\begin{equation}
[\sqrt{-1} X_{q,N+1-q},N^q_l]=N^q_l\tag{3.2}
\end{equation}
for $1\leq q\leq [N/2]$ and $q< l<N+1-q$  
if $N$ is odd 
or for $1\leq q<[N/2]$ and $q< l<N+1-q$ if $N$ is even.
\par 
We denote by $T_{s}(N)$ 
the collection of all increasing    
sequences $I=(i_1,\cdots,i_s)$ of length $s$ in $\{1,\cdots,N\}$ 
for $N,s\in\mathbb Z^{+}_{\geq 0}$. 
We denote by $P(I)$ the set of all the permutations of $I$. 
We denote by $\epsilon(\sigma)$ the sign of $\sigma$. 
We have 
$$\epsilon(\sigma)
=\prod_{1\leq a<b \leq s}
\frac{\sigma(i_b)-\sigma(i_a)}{i_b-i_a}
=\prod_{1\leq a<b\leq s}
\frac{\sigma(i_b)-\sigma(i_a)}{\vert\sigma(i_b)-\sigma(i_a)\vert }
\text{ for }\sigma\in P(I).$$ 
Assume that $p\geq 1$ and $I\in T_{2p}(N)$. We put 
\begin{align*}  
P_+(I) & =\{\sigma\in P(I)\bigm\vert 
\sigma(i_{2a-1})<\sigma(i_{2a+1})
\text{ for }1\leq a\leq p-1  
\text{ and }
\sigma(i_{2a-1})<\sigma(i_{2a})\text{ for }1\leq a\leq p 
\}
\\& 
=\{\sigma\in P(I) \bigm\vert 
\sigma(i_1)<\sigma(i_3)<\cdots<\sigma(i_{2p-1})
\text{ and } 
\sigma(i_1)<\sigma(i_2),\cdots,\sigma(i_{2p-1})<\sigma(i_{2p})\}.
\end{align*} 
We put as follows. 
$${\tilde W}_{I} 
=\sum_{\sigma\in P_+(I)} 
\epsilon(\sigma){\tilde X}_{\sigma(i_1)\sigma(i_2)}
\cdots {\tilde X}_{ \sigma(i_{2p-1}) \sigma(i_{2p}) }
\in U(\mathfrak s\mathfrak o(N,\mathbb F)).$$ 
$$\tilde W_p=\sum\limits_{J\in T_{2p}(N)}{\tilde W}_{J}^2
\in U(\mathfrak s\mathfrak o(N,\mathbb F)).$$ 
We put $\tilde W_\emptyset =1$. 
\par 
\textbf{Fact 3-1(Pfaffian type elements 
for $U(\tilde{\mathfrak g})$).}  
\par 
Let $N$ and $p$ be positive integers such that $2p\leq N$. 
Then the following $(1)$ and $(2)$ hold.
\par 
$(1)$(See Lemma 2.2 in $[8]$) 
For $u=(u_{ij})\in O(N,\mathbb F)$ 
and $I=(i_1,\cdots,i_{2p})\in T_{2p}(N)$, 
we have 
$$Ad(u){\tilde W}_{I}=\sum\limits_{J=(j_1,\cdots,j_{2p})\in T_{2p}(n)}
\begin{vmatrix}
   u_{i_1,j_1}& \cdots & u_{i_1,j_{2p}} \\
   \vdots & \ddots & \vdots \\
   u_{i_{2p},j_1}& \cdots & u_{i_{2p},j_{2p}} 
\end{vmatrix}\cdot{\tilde W}_{J}.$$ 
\par 
$(2)$(See Theorem 2.3 in $[8]$) 
We have 
$\tilde W_p\
in\mathcal Z(\mathfrak s\mathfrak o(N,\mathbb F))$. 
\par 
\textbf{Proposition 3-2(Laplace Expansions 
for $U(\tilde{\mathfrak g})$)}
\par  
Assume that $p\geq 0,2p+2\leq N$ and $I\in T_{2p+2}(N)$. 
\par 
$(1)$ Let $1\leq q\leq 2p+2$. We have 
$$\tilde W_{I}=\sum\limits_{1\leq l\leq 2p+2,l\not=q}
\frac{q-l}{\vert q-l\vert} (-1)^{l+q-1}
\tilde X_{i_l,i_q}\tilde W_{I\setminus\{i_l,i_q\}}.$$ 
\par 
$(2)$ We have 
$$\tilde W_{I}
=\sum\limits_{1\leq l\leq 2p+1}
(-1)^{l-1}
\tilde X_{i_l,i_{2p+2}}\tilde W_{I\setminus\{i_l,i_{2p+2}\}}.$$ 
\begin{Proof}
$(1)$ The assertion of Proposition 3-2(1) and its proof 
follows from Lemma 2.3. in [21] and in Proposition 2.3.in [22]. 
\par 
$(2)$ The assertion follows from (1) in the setting of $q=2p+2$. \qed 
\end{Proof}
\textbf{Lemma 3-3}
\par 
Let $N$ and $p$ be positive integers such that $2p+2\leq N$. 
Let $I=(i_1,\cdots,i_{2p+2})\in T_{2p+2}(N)$. 
Let $1\leq s,t\leq N$.
Then the following $(1)$ and $(2)$ hold.
\par 
$(1)$ Assume that $s,t\in I$ or $s,t\notin I$. Then we have 
$$\tilde W_{I}\tilde X_{st}=\tilde X_{st} {\tilde W}_{I}.$$
\par 
$(2)$ Assume that  $s\in I,t\notin I$. 
Let $s=i_q$ and $t=j_r$ for $1\leq q,r\leq 2p+2$ 
where we put 
$$J=I\cup\{ t\}\setminus\{ s\}\in T_{2p+2}(N).$$ 
Then we have 
$${\tilde W}_{I}\tilde X_{st}=\tilde X_{st} {\tilde W}_{I}
-(-1)^{q-r}\tilde W_{I\cup\{t\}\setminus\{s \}}.$$ 
\begin{Proof}
$(1)$ It follows directly from Fact 3-1(1) that 
$$Ad(\exp(c\tilde{X}_{st}))\tilde{W}
=\tilde{W}\text{ for }c\in\mathbb R.$$ 
This yields the assertion. 
\par 
$(2)$ By Proposition 3-2(1) we have 
\begin{equation} \tilde W_{I}\tilde X_{st}
=\sum\limits_{1\leq l\leq 2p+2,l\not=q}
\frac{q-l}{\vert q-l\vert} (-1)^{l+q-1}
\tilde X_{i_l,i_q}
\tilde W_{I\setminus\{i_l,i_q\}}\tilde X_{i_q t} .
\tag{ 3.3.1}
\end{equation}
By Lemma 3-3(1) we have 
$$\tilde W_{I\setminus\{i_l,i_q\}}\tilde X_{i_q,t}
=\tilde X_{i_q,t}\tilde W_{I\setminus\{i_l,i_q\}}.$$ 
Hence we have 
\begin{equation}
\tilde W_{I}\tilde X_{st}
=\sum\limits_{1\leq l\leq 2p+2,l\not=q}
\frac{q-l}{\vert q-l\vert} (-1)^{l+q-1}
\tilde X_{i_l,i_q}\tilde X_{i_q t}
\tilde W_{I\setminus\{i_l,i_q\}}.
\tag{ 3.3.2}
\end{equation}

Here we have 
\begin{equation} \tilde X_{i_l,i_q}\tilde X_{i_q,t} 
=\tilde X_{i_q,t} \tilde X_{i_l,i_q}-\tilde X_{i_l,t}.
\tag{ 3.3.3}
\end{equation} 

By (3.3.1),(3.3.2),(3.3.3) and Proposition 3-2(1), we have 
\begin{equation}{\tilde W}_{I}\tilde X_{rs} 
=\tilde X_{st} {\tilde W}_{I}
-\sum\limits_{1\leq l\leq 2p+2,l\not=q}
\frac{q-l}{\vert q-l\vert} (-1)^{l+q-1}
\tilde X_{i_l,t}\tilde W_{I\setminus\{i_l,t\}}.
\tag{ 3.3.4}
\end{equation}
Here we have by Proposition 3-2(1),  
\begin{equation}
\tilde W_{I\cup\{ t\}\setminus\{ s\}}=\tilde W_{J}
=\sum\limits_{1\leq l\leq 2p+2,l\not=r}
\frac{r-l}{\vert r-l\vert} (-1)^{l+r-1}
\tilde X_{j_l,j_r}\tilde W_{I\setminus\{i_{l},t\}}.
\tag{ 3.3.5}
\end{equation}
We assume that $s<t$. 
Then we have 
\begin{equation}
j_l=\begin{cases} i_l  &  \text{ for }1\leq l<q \\ 
         i_{l+1}  & \text{ for } q\leq l <r   \\ 
         t  & \text{ for } l=r \\ 
         i_{l} & \text{ for }r<l\leq 2p+2.   
\end{cases}\tag{ 3.3.6}
\end{equation} 
We assume that $t<s$. 
Then we have 
\begin{equation}j_l
=\begin{cases} i_l  &  \text{ for }1\leq l<r \\ 
         t  & \text{ for } l=r   \\ 
         i_{l-1}  & \text{ for } r<l\leq q \\ 
         i_{l} & \text{ for }q<l\leq 2p+2.   
\end{cases}\tag{ 3.3.7}
\end{equation}
In view of (3.3.4) it follows from (3.3.5),(3.3.6) and (3.3.7) 
that we have  
\begin{equation}
\sum\limits_{1\leq l\leq 2p+2,l\not=q}
\frac{q-l}{\vert q-l\vert} (-1)^{l+q-1}
\tilde X_{i_l,t}\tilde W_{I\setminus\{i_l,i_q\}}
=(-1)^{q-r}\tilde W_{I\cup\{ t\}\setminus\{ s\}}.
\tag{ 3.3.8}
\end{equation}
Then the assertion follows from (3.3.4) and (3.3.8). \qed 
\end{Proof}
\textbf{Remark} 
The assertion of Lemma 3-3(2) can be deduced 
from the assertion for $q=r$ in Lemma 3-3(2), 
Fact 3-1(1) and (2) in [22]. 
In the proof of Lemma 3-3 we provide a direct explicit proof. 
\par 
\textbf{Lemma 3-4}
\par 
Assume that $p\geq 0,2p+4\leq N$ and $I\in T_{2p+4}(N)$. 
\par 
$(1)$ Let $1\leq q,r\leq 2p+4,q\not=r$. We have  
$$\sum\limits_{1\leq l\leq 2p+4,l\not=q,r}
\frac{q-l}{\vert q-l\vert}\frac{r-l}{\vert r-l\vert}(-1)^{l-1}
\tilde X_{i_q,i_l}
\tilde W_{I\setminus\{i_l,i_r\}}
=(-1)^{q-1}(p+1)\tilde W_{I\setminus\{i_q,i_r\}}.$$ 
\par 
$(2)$ We have 
$$\sum\limits_{2\leq l\leq 2p+3}(-1)^{l}
\tilde X_{i_1,i_l}
\tilde W_{I\setminus\{i_l,i_{2p+4}\}}
=(p+1)\tilde W_{I\setminus\{i_1,i_{2p+4}\}}.$$ 
\begin{Proof}
$(1)$ We put 
$$\mathcal I=\sum\limits_{1\leq l\leq 2p+4,l\not=q,r}
\frac{q-l}{\vert q-l\vert}\frac{r-l}
{\vert r-l\vert}(-1)^{l-1}\tilde X_{i_q,i_l}
\tilde W_{I\setminus\{i_l,i_r\}}.$$ 
By Proposition 3-2(1) we have for $1\leq l\leq 2p+4,l\not=q,r$, 
\begin{align*} 
& 
\tilde W_{I\setminus\{i_l,i_r\}}
=\sum\limits_{1\leq k\leq 2p+4,k\not=q,l,r} 
\frac{q-k}{\vert q-k\vert}
\frac{r-k}{\vert r-k\vert}
\frac{l-k}{\vert l-k\vert}
(-1)^{k+q-1}\tilde X_{i_k,i_q}
\tilde W_{I\setminus\{i_q,i_l,i_k,i_r\}} 
\\&
=\sum\limits_{1\leq k\leq 2p+4,k\not=q,l,r} 
\frac{q-k}{\vert q-k\vert}
\frac{r-k}{\vert r-k\vert}
\frac{l-k}{\vert l-k\vert}
(-1)^{k+q}\tilde X_{i_q,i_k}
\tilde W_{I\setminus\{i_q,i_l,i_k,i_r\}}.
\end{align*}
Then we have 
\begin{align*} & \mathcal I
=\sum\limits_{1\leq l\leq 2p+4,l\not=q,r}
\frac{q-l}{\vert q-l\vert}\frac{r-l}{\vert r-l\vert}
(-1)^{l-1}(\tilde X_{i_q,i_l}  
\\& 
\times(\sum\limits_{1\leq k\leq 2p+4,k\not=l,q,r}
\frac{q-k}{\vert q-k\vert}
\frac{r-k}{\vert r-k\vert}
\frac{l-k}{\vert l-k\vert}
(-1)^{k+q}
\tilde X_{i_q,i_k}W_{I\setminus{\{i_q,i_k,i_l,i_r\}}}) 
\\& 
=\sum\limits_{1\leq l\leq 2p+4,l\not=q,r}
\frac{q-l}{\vert q-l\vert}\frac{r-l}{\vert r-l\vert} 
\sum\limits_{1\leq k\leq 2p+4,k\not=l,q,r}
\frac{q-k}{\vert q-k\vert}\frac{r-k}{\vert r-k\vert}
\frac{l-k}{\vert l-k\vert}
(-1)^{k+l+q-1}\tilde X_{i_q,i_l}\tilde X_{i_q,i_k}
\tilde W_{I\setminus\{i_q,i_k,i_l,i_r\}}
\\& 
=\sum\limits_{1\leq k,l\leq 2p+4,k\not=l,q,r,l\not=q,r}
\frac{q-l}{\vert q-l\vert}\frac{r-l}{\vert r-l\vert}
\frac{q-k}{\vert q-k\vert}\frac{r-k}{\vert r-k\vert} 
\frac{l-k}{\vert l-k\vert}
(-1)^{k+l+q-1}\tilde X_{i_q,i_l}\tilde X_{i_q,i_k}
\tilde W_{I\setminus\{i_q,i_k,i_l,i_r\}} 
\\& 
=\sum\limits_{1\leq k<l\leq 2p+4,k,l\not=q,r}
\frac{q-l}{\vert q-l\vert}\frac{r-l}{\vert r-l\vert}
\frac{q-k}{\vert q-k\vert}\frac{r-k}{\vert r-k\vert}
(-1)^{k+l+q-1}\tilde X_{i_q,i_l}\tilde X_{i_q,i_k}
\tilde W_{I\setminus\{i_q,i_k,i_l,i_r\}}
\\&
-\sum\limits_{1\leq l<k\leq 2p+4,k,l\not=q,r}
\frac{q-l}{\vert q-l\vert}\frac{r-l}{\vert r-l\vert}
\frac{q-k}{\vert q-k\vert}\frac{r-k}{\vert r-k\vert}
(-1)^{k+l+q-1}\tilde X_{i_q,i_l}\tilde X_{i_q,i_k}
\tilde W_{I\setminus\{i_q,i_k,i_l,i_r\}}.
\end{align*}
We have 
$$\tilde X_{i_q,i_l}\tilde X_{i_q,i_k}
-\tilde X_{i_q,i_k}\tilde X_{i_q,i_l}
=\tilde X_{i_l,i_k}=-\tilde X_{i_k,i_l}.$$
Then we have 
\begin{align*} & 
\mathcal I=\sum\limits_{1\leq k<l\leq 2p+4,k,l\not=q,r}
\frac{q-l}{\vert q-l\vert}\frac{r-l}{\vert r-l\vert}
\frac{q-k}{\vert q-k\vert}\frac{r-k}{\vert r-k\vert}
(-1)^{k+l+q-1}(\tilde X_{i_q,i_k}\tilde X_{i_q,i_l}
-\tilde X_{i_k,i_l})
\tilde W_{I\setminus\{i_q,i_k,i_l,i_r\}} 
\\& 
-\sum\limits_{1\leq l<k\leq 2p+4,k,l\not=q,r}
\frac{q-l}{\vert q-l\vert}\frac{r-l}{\vert r-l\vert}
\frac{q-k}{\vert q-k\vert}\frac{r-k}{\vert r-k\vert}
(-1)^{k+l+q-1}\tilde X_{i_q,i_l}\tilde X_{i_q,i_k} 
\tilde W_{I\setminus\{i_q,i_k,i_l,i_r\}} 
\\& 
=\sum\limits_{1\leq k<l\leq 2p+4,k,l\not=q,r}
\frac{q-l}{\vert q-l\vert}\frac{r-l}{\vert r-l\vert}
\frac{q-k}{\vert q-k\vert}\frac{r-k}{\vert r-k\vert}
(-1)^{k+l+q}\tilde X_{i_k,i_l}\tilde W_{I\setminus\{i_q,i_k,i_l,i_r\}}.
\end{align*}
First we assume $p=0$. We can easily verify that 
$$\frac{q-l}{\vert q-l\vert}\frac{r-l}{\vert r-l\vert}
\frac{q-k}{\vert q-k\vert}\frac{r-k}{\vert r-k\vert}(-1)^{k+l}=-1.$$ 
This directly yields the assertion in view of 
$\tilde W_\emptyset=1$.
\par 
Next we assume that $p\geq 1$. We have 
$$\tilde W_{I\setminus\{i_q,i_k,i_l,i_r\}}
=\sum\limits_{\tau\in P_+(I\setminus\{i_q,i_k,i_l,i_r\})}
\epsilon(\tau)\prod\limits_{1\leq a\leq p} 
\tilde X_{\tau(j^{k,l}_{2a-1})\tau(j^{k,l}_{2a})},$$ 
where we put 
$$I\setminus\{i_q,i_k,i_l,i_r\}
=(j^{k,l}_1,\ldots,j^{k,l}_s,\ldots,j^{k,l}_{2p}).$$ 
Hence we have 
\begin{align*} & 
\mathcal I=(-1)^{q-1}\sum\limits_{1\leq k<l\leq 2p+4,k,l\not=q,r}
\sum\limits_{\tau\in P_+(I\setminus\{i_q,i_k,i_l,i_r\})}
\frac{q-l}{\vert q-l\vert}\frac{r-l}{\vert r-l\vert}
\\& 
\times\frac{q-k}{\vert q-k\vert}\frac{r-k}{\vert r-k\vert}
(-1)^{k+l-1}\epsilon(\tau)\tilde X_{i_k,i_l}
\prod\limits_{1\leq a\leq p} 
\tilde X_{\tau(j^{k,l}_{2a-1})\tau(j^{k,l}_{2a})}. 
\end{align*}  
We put 
$$I\setminus\{i_q,i_r\}
=(j_1,\ldots,j_s,\ldots,j_{2p+2}).$$
For any $\sigma\in P_+(I\setminus\{i_q,i_r\})$, we define 
\begin{align*} & I(\sigma)
=\{((k,l,\tau)
\bigm\vert 
1\leq k<l\leq 2p+4,k,l\not=q,r,
\tau\in P_+(I\setminus\{i_q,i_k,i_l,i_r\})
\\& 
\tilde X_{i_k,i_l}\prod\limits_{1\leq a\leq p}
\tilde X_{\tau(j^{k,l}_{2a-1})\tau(j^{k,l}_{2a})}
=\prod\limits_{1\leq a\leq p+1}
\tilde X_{\sigma(j_{2a-1})\sigma(j_{2a})}.\}
\end{align*} 
Then we have 
\begin{align*} & 
\mathcal I=(-1)^{q-1}\sum\limits_{\sigma\in P_+(I\setminus\{ i_q,i_r\})}
\sum\limits_{(k,l,\tau)\in I(\sigma)}
\frac{q-l}{\vert q-l\vert}\frac{r-l}{\vert r-l\vert}
\\&
\times\frac{q-k}{\vert q-k\vert}\frac{r-k}{\vert r-k\vert}
(-1)^{k+l-1}\epsilon(\tau)\prod\limits_{1\leq a\leq p+1}
\tilde X_{\sigma(j_{2a-1})\sigma(j_{2a})}. 
\tag{ 3.4.1}
\end{align*}  
We define $\rho\in P(I\setminus\{i_q,i_r\})$ by 
$$(\rho(j_1),\rho(j_2),\cdots,\rho(j_{2p+2}))
=(i_k,i_l,\tau(j^{k,l}_{1}),\cdots,j^{k,l}_{2p})).$$ 
Since we have $k<l$, we have 
\begin{align*} & \epsilon(\rho)
=\prod_{1\leq a<b\leq 2p+2}
\frac{\rho(j_b)-\rho(j_a)}{\vert\rho(j_b)-\rho(j_a)\vert}
\\& 
=\prod_{3\leq b\leq 2p+2}
\frac{\rho(j_b)-\rho(j_1)}{\vert\rho(j_b)-\rho(j_1)\vert}
\times 
\prod_{3\leq b\leq 2p+2}
\frac{\rho(j_b)-\rho(j_2)}{\vert\rho(j_b)-\rho(j_2)\vert}
\times
\prod_{3\leq a<b\leq 2p+2}
\frac{\rho(j^{k,l}_b)-\rho(j^{k,l}_a)}
{\vert\rho(j^{k,l}_b)-\rho(j^{k,l}_a)\vert} 
\\& 
=\prod_{3\leq b\leq 2p+2}
\frac{\rho(j_b)-\rho(j_1)}{\vert\rho(j_b)-\rho(j_1)\vert}
\times 
\prod_{3\leq b\leq 2p+2}
\frac{\rho(j_b)-\rho(j_2)}{\vert\rho(j_b)-\rho(j_2)\vert}
\times
\prod_{1\leq a<b\leq 2p}
\frac{\tau(j^{k,l}_b)-\tau(j^{k,l}_a)}
{\vert\tau(j^{k,l}_b)-\tau(j^{k,l}_a)\vert} 
\\& 
=\prod_{3\leq b\leq 2p+2}
\frac{\rho(j_b)-i_k}{\vert\rho(j_b)-i_k\vert}
\times 
\prod_{3\leq b\leq 2p+2}
\frac{\rho(j_b)-i_l}{\vert\rho(j_b)-i_l\vert}
\times \epsilon(\tau). \tag{ 3.4.2}
\end{align*} 
We have  
$$\{\rho(j_b)\bigm\vert \rho(j_b)<i_k(3\leq b\leq 2p+2)\}
=\{i_c\bigm\vert 1\leq c<k\}\setminus\{q,r\}$$ 
and 
$$\{\rho(j_b)\bigm\vert \rho(j_b)<i_l(3\leq b\leq 2p+2)\}
=\{i_c\bigm\vert 1\leq c<l, c\not=k\}\setminus\{q,r\}.$$ 
Then we have 
$$\prod_{3\leq b\leq 2p+2}
\frac{\rho(j_b)-i_k}{\vert\rho(j_b)-i_k\vert}
=\frac{q-k}{\vert q-k\vert}\frac{r-k}{\vert r-k\vert}(-1)^{k-1}$$
and 
$$\prod_{3\leq b\leq 2p+2}
\frac{\rho(j_b)-i_l}{\vert\rho(j_b)-i_l\vert}
=\frac{q-l}{\vert q-l\vert}\frac{r-l}{\vert r-l\vert}(-1)^{l-2}.$$
Hence we have by (3.4.2), 
$$\epsilon(\rho)
=\frac{q-l}{\vert q-l\vert}\frac{r-l}{\vert r-l\vert}
\frac{q-k}{\vert q-k\vert}\frac{r-k}{\vert r-k\vert}
(-1)^{k+l-1}\epsilon(\tau).$$ 
By the repeated application of (4.2.3.21) in [25] we have 
$\epsilon(\rho)=\epsilon(\sigma)$. 
Hence we have for $(k,l,\tau)\in I(\sigma)$, 
$$\epsilon(\sigma)= 
\frac{q-l}{\vert q-l\vert}\frac{r-l}{\vert r-l\vert}
\frac{q-k}{\vert q-k\vert}\frac{r-k}{\vert r-k\vert}
(-1)^{k+l-1}\epsilon(\tau).$$ 
We have $\# I(\sigma)=p+1$ 
for any $\sigma\in P_+(I\setminus\{i_q,i_r\})$. 
Then we have by (3.4.1), 
\begin{align*} & 
\mathcal I=(-1)^{q-1}\sum\limits_{\sigma\in P_+(I\setminus\{i_q,i_r\})}
\sum\limits_{(k,l,\tau)\in I(\sigma)}
\epsilon(\sigma)
\prod\limits_{1\leq a\leq p+1}\tilde X_{\sigma(j_{2a-1})\sigma(j_{2a})}
\\& 
=(-1)^{q-1}(p+1)\sum\limits_{\sigma\in P_+(I\setminus\{i_q,i_r\})}
\epsilon(\sigma)
\prod\limits_{1\leq a\leq p+1}\tilde X_{\sigma(j_{2a-1})\sigma(j_{2a})}
\end{align*}
Here we have by definition 
$$\tilde W_{I\setminus\{i_q,i_r\}}
=\sum\limits_{\sigma\in P_+(I\setminus\{i_q,i_r\})}
\epsilon(\sigma)
\prod\limits_{1\leq a\leq p+1}\tilde X_{\sigma(j_{2a-1})\sigma(j_{2a})}.$$
Then we have $\mathcal I=(-1)^{q-1}(p+1)\tilde W_{I\setminus\{i_q,i_r\}}$. 
Thus we have proved the assertion. 
\par 
$(2)$ The assertion follows from (1) in the setting of $(q,r)=(1,2p+4)$. \qed 
\end{Proof}
\textbf{Remark} The assertion of Lemma 3-4(1) can be deduced 
from the assertion of Lemma 3-4(2), Fact 3-1(1) and (2) in [22]. 
In the proof of Lemma 3-4 we provide a direct explicit proof. 
\par 
Assume that $p\geq 0,2p+4\leq N,1\leq s\leq [N/2]$. 
We define subsets $(J)^s_-,(J)^s_0$ and $(J)^s_+$ 
of $J\in T_{2p+4}(N)$ by the followings.  
$$(J)^s_{-}=\{j_l\in J\bigm\vert 1\leq j_l\leq s\}.$$
$$(J)^s_0=\{j_l\in J\bigm\vert s+1\leq j_l\leq N-s\}.$$ 
$$(J)^s_+=\{j_l\in J\bigm\vert N-s+1\leq j_l\leq N\}.$$ 
We also define a subset $(J)^s_1$ of $J$ by 
$$(J)^s_1=\{j_l\in (J)^s_+\bigm\vert N+1-j^s_{l}\in (J)^s_{-}\}.$$ 
We define the Lie subalgebra $\tilde{\mathfrak n}(0)^s$ 
of $\tilde{\mathfrak n}(0)$ for $1\leq s\leq [N/2]$  by 
$$\tilde{\mathfrak n}(0)^s=\bigoplus\limits_{1\leq i\leq s}
\tilde{\mathfrak n}(0)_i.$$
\par 
\textbf{Lemma 3-5}
\par 
Assume that $p\geq 0,2p+4\leq N,1\leq r\leq [N/2]$ and $I\in T_{2p+4}(N)$. 
\par 
$(1)$ Suppose that $(I)^s_1\not=\emptyset$.
Let $i_r\in(I)^s_1$ and $i_q\in(I)^s_-$ for $1\leq q<r\leq 2p+4$
such that $i_q=N+1-i_r$. Then we have 
$$\tilde W_I\in 
\sum\limits_{1\leq l\leq q,r<l\leq 2p+4}
\tilde W_{I\setminus\{i_l,i_{r}\}}U(\tilde{\mathfrak g})
+\tilde{\mathfrak n}(0)_{i_q}U(\tilde{\mathfrak g}).$$  
\par 
$(2)$ We have 
$$\tilde W_I\in
\sum\limits_{J\subset I,
\#(J)^s_-+\#(J)^s_+=(\#(I)^s_-+\#(I)^s_+)-2\#(I)^s_1,
\#(J)^s_0=\#(I)^s_0,\#(J)^s_1=0}
\tilde W_J U(\tilde{\mathfrak g})
+(\tilde{\mathfrak n}(0)^s)U(\tilde{\mathfrak g}).$$ 
\begin{Proof} 
$(1)$ We have by Proposition 3-2(1), 
\begin{align*}  & 
\tilde W_{I} =\sum\limits_{1\leq l\leq 2p+4,l\not=r}
\frac{r-l}{\vert r-l\vert}(-1)^{l+r-1}
\tilde X_{i_l,i_{r}} \tilde W_{I\setminus\{i_l,i_{r}\}}
\\& 
=\sum\limits_{1\leq l\leq q,r<l\leq 2p+4}
\frac{r-l}{\vert r-l\vert}(-1)^{l+r-1}
\tilde X_{i_l,i_{r}} \tilde W_{I\setminus\{i_l,i_{r}\}}
+\sum\limits_{q<l<r}
\frac{r-l}{\vert r-l\vert}(-1)^{l+r-1}
\tilde X_{i_l,i_{r}} \tilde W_{I\setminus\{i_l,i_{r}\}}
\\& 
\in\sum\limits_{1\leq l\leq q,r<l\leq 2p+4}
\frac{r-l}{\vert r-l\vert}(-1)^{l+r-1}
\tilde X_{i_l,i_{r}} \tilde W_{I\setminus\{i_l,i_{r}\}}
\\& 
-(\sqrt{-1})\sum\limits_{q<l<r}
\frac{r-l}{\vert r-l\vert}(-1)^{l+r-1}
\tilde X_{i_q,i_l} \tilde W_{I\setminus\{i_l,i_{r}\}}
+\tilde{\mathfrak n}(0)_q U(\tilde{\mathfrak g}),
\tag{ 3.5.1}
\end{align*}
where we have by (3.1),
$$\tilde N^{i_q}_{i_l}
=\tilde X_{i_l,i_{r}}+\sqrt{-1}\cdot
\tilde X_{i_q,i_l}\in\tilde{\mathfrak n}(0)_{i_q}
\text{ for }q<l<r.$$
We have 
\begin{align*} & \sum\limits_{q<l<r}
\frac{r-l}{\vert r-l\vert}(-1)^{l+r-1}
\tilde X_{i_q,i_l}\tilde W_{I\setminus\{i_l,i_{r}\}}
=-\sum\limits_{q<l<r}
\frac{q-l}{\vert q-l\vert}\frac{r-l}{\vert r-l\vert}(-1)^{l+r-1}
\tilde X_{i_q,i_l} \tilde W_{I\setminus\{i_l,i_{r}\}}
\\&
=\sum\limits_{1\leq l<q,r<l\leq 2p+4}
\frac{q-l}{\vert q-l\vert}\frac{r-l}{\vert r-l\vert}(-1)^{l+r-1}
\tilde X_{i_q,i_l} \tilde W_{I\setminus\{i_l,i_{r}\}}
\\& 
-\sum\limits_{1\leq l\leq 2p+4,l\not=q,l\not=r}
\frac{q-l}{\vert q-l\vert}\frac{r-l}{\vert r-l\vert}(-1)^{l+r-1}
\tilde X_{i_q,i_l} \tilde W_{I\setminus\{i_l,i_{r}\}}.
\end{align*}
By Lemma 3-4(1) we have 
\begin{align*} & \sum\limits_{q<l<r}
\frac{r-l}{\vert r-l\vert}(-1)^{l+r-1}
\tilde X_{i_q,i_l} \tilde W_{I\setminus\{i_l,i_{r}\}}
\\& 
=\sum\limits_{1\leq l<q,r<l\leq 2p+4}
\frac{q-l}{\vert q-l\vert}\frac{r-l}{\vert r-l\vert}(-1)^{l+r-1}
\tilde X_{i_q,i_l} \tilde W_{I\setminus\{i_l,i_{r}\}}
-(-1)^{q-1}(p+1)\tilde W_{I\setminus\{i_q,i_{r}\}}.
\tag{3.5.2}
\end{align*}
Then the assertion follows from (3.5.1) and (3.5.2). 
\par 
$(2)$  Suppose that $(I)^s_1=\emptyset$. 
Then we have $\#(I)^s_1=0$. 
Hence the assertion follows if we just put $J=I$. 
\par 
Suppose that $(I)^s_1\not=\emptyset$. 
We assume that $i_r$ is the maximum element in $(I)^s_1$ 
in the setting of Lemma 3-5(1). 
Then we have by Lemma 3-5(1), 
\begin{equation*}\tilde W_I\in
\sum\limits_{J\subset I,\#(J)^s_-+\#(J)^s_+=(\#(I)^s_-+\#(I)^s_+)-2,
\#(J)^s_0=\#(I)^s_0,
\#(J)^s_1=\#(I)^s_1-1}
\tilde W_J U(\tilde{\mathfrak g})
+\tilde{\mathfrak n}(0)^s U(\tilde{\mathfrak g}).
\tag{ 3.5.3 }
\end{equation*}
The repeated application of Lemma 3-5(1) to $\tilde W_J$ 
in the right-hand side of (3.5.3) yields the assertion. \qed 
\end{Proof}
\textbf{Lemma 3-6}
\par 
Assume that $p\geq 0,2p+4\leq N$ and $I\in T_{2p+4}(N)$. 
\par 
$(1)$ Assume that $1\in I$ and $N\in I$. 
Then we have for $2\leq q\leq 2p+3$, 
$$\tilde W_{I\setminus\{1,N\}}\tilde N^1_{i_q}
=\tilde N^1_{i_q}\tilde W_{I\setminus\{1,N\}}
+\sum\limits_{2\leq k\leq 2p+3,k\not= q}
\frac{q-k}{\vert q-k\vert}(-1)^{k+q}
tilde N^1_{i_k}\tilde W_{I\setminus\{1,i_q,i_k,N\}}.$$ 
\par 
$(2)$ Assume $1\in I$ and $N\not\in I$. 
Then we have for $2\leq q\leq 2p+4$,  
$$\tilde W_{I}\tilde N^1_{i_q}=\tilde N^1_{i_q}\tilde W_{I}
+(-1)^{q-1} \tilde W_{I\cup\{N\}\setminus\{i_q\}}.$$ 
\begin{Proof}
$(1)$ Here we have $i_1=1$and $i_{2p+4}=N$. 
Here we put 
$$J=(j_1,j_2,\ldots,j_{2p+2})=(i_2,i_3,\ldots,i_{2p+3})
=I\setminus\{1,N\}.$$ 
Then we have $i_q=j_{q-1}$. 
By Proposition 3-2(1), we have for $2\leq q\leq 2p+3$, 
$$\tilde W_{J}=\sum\limits_{1\leq l\leq 2p+2,l\not=q-1}
\frac{(q-1)-l}{\vert (q-1)-l\vert} (-1)^{l+(q-1)-1}
\tilde X_{j_l,j_{q-1}}\tilde W_{J\setminus\{j_l,j_{q-1}\}}.$$ 
Since we have $i_k=j_{l}$ for $k=l+1$, 
we have for $2\leq q\leq 2p+3$, 
\begin{equation*}\tilde W_{I\setminus\{1,N\}}
=\sum\limits_{2\leq k\leq 2p+3,k\not=q}
\frac{q-k}{\vert q-k\vert}(-1)^{k+q-1}
\tilde X_{i_k,i_q}\tilde W_{I\setminus\{1,i_k,i_q,N\}}.
\tag{ 3.6.1} 
\end{equation*}
We have by (3.1), 
\begin{equation*}\tilde N^1_{i_q}
=\tilde X_{i_q,N}+\sqrt{-1}\cdot \tilde X_{1,i_q}.
\tag{3.6.2} 
\end{equation*}
Then we have by (3.6.1),(3.6.2) and Lemma 3-3(1), 
\begin{align*} & 
\tilde W_{I\setminus\{1,N\}}\tilde N^1_{i_q} 
=\sum\limits_{2\leq k\leq 2p+3,k\not=q}
\frac{q-k}{\vert q-k\vert}(-1)^{k+q-1}
\tilde X_{i_k,i_q}\tilde W_{I\setminus\{1,i_k,i_q,N\}}
\tilde N^1_{i_q}
\\& 
=\sum\limits_{2\leq k\leq 2p+3,k\not=l}
\frac{q-k}{\vert q-k\vert}(-1)^{k+q-1}
\tilde X_{i_k,i_q} \tilde N^1_{i_q} 
\tilde W_{I\setminus\{1,i_k,i_q,N\}}.
\tag{]3.6.3}
\end{align*}
Since we have 
$$\tilde X_{i_k,i_q} \tilde X_{1,i_q}
=\tilde X_{1,i_q}\tilde X_{i_k,i_q}-\tilde X_{1,i_k}
\text{ and }
\tilde X_{i_k,i_q} \tilde X_{i_q,N}
=\tilde X_{i_q,m+n}\tilde X_{i_k,i_q}-\tilde X_{i_q,N},$$ 
we have by (3.6.2), 
\begin{equation*}
\tilde X_{i_k,i_q}\tilde N^1_{i_q}
=\tilde N^1_{i_q}\tilde X_{i_k,i_q}
-\tilde N^1_{i_k}.\tag{ 3.6.4}
\end{equation*}  
It follows from (3.6.3) that we have by (3.6.4) and (3.6.1), 
\begin{align*} & \tilde W_{I\setminus\{1,N\}}\tilde N^1_{i_q} 
=\sum\limits_{2\leq k\leq 2p+3,k\not=q}
\frac{q-k}{\vert q-k\vert}(-1)^{k+q-1}
(\tilde N^1_{i_q}\tilde X_{i_k,i_q}-\tilde N^1_{i_k})
\tilde W_{I\setminus\{1,i_k,i_q,N\}} 
\\& 
=\tilde N^1_{i_q}\sum\limits_{2\leq k\leq 2p+3,k\not=q}
\frac{q-k}{\vert q-k\vert}(-1)^{k+q-1}
\tilde X_{i_k,i_q}\tilde W_{I\setminus\{1,i_k,i_q,N\}} 
\\& 
+\sum\limits_{2\leq k\leq 2p+3,k\not=q}
\frac{q-k}{\vert q-k\vert}(-1)^{k+q}
\tilde N^1_{i_k}\tilde W_{I\setminus\{1,i_q,i_l,N\}} 
\\&  
=\tilde N^1_{i_q}\tilde W_{I\setminus\{1,N\}}
+\sum\limits_{2\leq k\leq 2p+3,k\not=q}
\frac{q-k}{\vert q-k\vert}(-1)^{k+q}
\tilde N^1_{i_k}\tilde W_{I\setminus\{1,i_k,i_q,N\}}.
\end{align*}
This yields the assertion. 
\par 
$(2)$ Here we have $i_1=1$ and $i_{2p+4}<N$.  
By Lemma 3-3(1) we have 
\begin{equation*}
\tilde W_I \tilde X_{1,i_q}=\tilde X_{1,i_q}\tilde W_I.
\tag{3.6.5}
\end{equation*}
In the setting of Lemma 3-3(2) we put $s=i_q$ and $t=N$. 
Then we have $r=2p+4$. Then we have $(-1)^{q-r}=(-1)^{q}$.   
Hence we have by Lemma 3-3(2) 
\begin{equation*}\tilde W_I \tilde X_{i_q,N}
=\tilde X_{i_q,N}\tilde W_I
-(-1)^{q}\tilde W_{I\cup\{ N\}\setminus\{ i_q\}}.
\tag{3.6.6}
\end{equation*}
Then we have by (3.6.2),(3.6.5) and (3.6.6), 
\begin{align*}  & \tilde W_I \tilde N^1_{i_q}
=\sqrt{-1}\cdot\tilde W_I\tilde X_{1,i_q}
+\tilde W_I \tilde X_{i_q,N}
\\& 
=\sqrt{-1}\cdot \tilde X_{1,i_q}\tilde W_I
+\tilde X_{i_q,N}\tilde W_I
-(-1)^{q}\tilde W_{I\cup\{N\}\setminus\{i_q\}}
\\& 
=\tilde N^1_{i_q}\tilde W_I
+(-1)^{q-1}\tilde W_{I\cup\{N\}\setminus\{i_q\}}
\end{align*} 
This yields the assertion. \qed 
\end{Proof}
\par 
\textbf{Proposition 3-7}
\par 
Assume that $p\geq 0,2p+4\leq N$ 
and $I\in T_{2p+4}(N)$. 
\par 
$(1)$ Assume that $1\in I$ and $N\in I$. Then we have 
$$\tilde W_I
\in(-\sqrt{-1})(\sqrt{-1} \tilde X_{1,N}-(p+1))\cdot 
\tilde W_{I\setminus\{1,N\}}
+\tilde{\mathfrak n}(0)_1 \cdot U(\tilde{\mathfrak g})$$
and 
$$\tilde W_I^2\in
-(\sqrt{-1} \tilde X_{1,N}-(p+1))^2
(\tilde W_{I\setminus\{1,N\}})^2
+\tilde{\mathfrak n}(0)_1\cdot U(\tilde{\mathfrak g}).$$
\par 
$(2)$  Assume that $1\not\in I$ and $N\in I$. Then we have 
$$\tilde W_I\in(-\sqrt{-1})\tilde W_{I\cup\{1\}\setminus\{N \}}
+\tilde{\mathfrak n}(0)_1\cdot U(\tilde{\mathfrak g})$$
and 
$$\tilde W_I^2\in - \tilde W_{I\cup\{1\}\setminus\{N \}}^2
+(\sqrt{-1} \tilde X_{1,N}-(p+1))\sum\limits_{1\leq l\leq 2p+3}
(\tilde W_{I\setminus\{i_l,N\}})^2
+\tilde{\mathfrak n}(0)_1\cdot U(\tilde{\mathfrak g}).$$ 
$(3)$ We have 
\begin{align*} & \tilde W_{p+2}\in 
-((\sqrt{-1} \tilde X_{1,N}-(p+1))(\sqrt{-1} \tilde X_{1,N}-(N-3-p))
\sum\limits_{1,N\in I}\tilde W_{I\setminus\{1,N\}}^2
\\& 
+\sum\limits_{1\not\in I,N\not\in I}\tilde W_I^2
+\tilde{\mathfrak n}(0)_1\cdot U(\tilde{\mathfrak g}).
\end{align*} 
\begin{Proof} 
$(1)$ Here we have $i_1=1$ and $i_{2p+4}=N$. 
By Proposition 3-2(2) and (3.1) we have 
\begin{align*}  & 
\tilde W_{I}=\sum\limits_{1\leq l\leq 2p+3}
(-1)^{l-1}\tilde X_{i_l,N}\tilde W_{I\setminus\{i_l,N\}}
\\& 
=\tilde X_{1,N}\tilde W_{I\setminus\{1,N\}}
+\sum\limits_{2\leq l\leq 2p+3}
(-1)^{l-1}\tilde X_{i_l,N}\tilde W_{I\setminus\{i_l,N\}}
\\& 
=(-\sqrt{-1})(\sqrt{-1}\tilde X_{1,N})\tilde W_{I\setminus\{1,N\}}
\\& 
+\sum\limits_{2\leq l\leq 2p+3}(-1)^{l-1}
\tilde N^1_{i_l}\tilde W_{I\setminus\{i_l,N\}}
-\sqrt{-1}(\sum\limits_{2\leq l\leq 2p+3}(-1)^{l-1}
(\tilde X_{1,i_l}\tilde W_{I\setminus\{i_l,N\}})
\\& 
=(-\sqrt{-1})(\sqrt{-1} \tilde X_{1,N})\tilde W_{I\setminus\{1,N\}}
\\& 
-(-\sqrt{-1})(\sum\limits_{2\leq l\leq 2p+3}(-1)^{l}
(\tilde X_{1,i_l}\tilde W_{I\setminus\{i_l,N\}})
+\sum\limits_{2\leq l\leq 2p+3}(-1)^{l-1}
\tilde N^1_{i_l}\tilde W_{I\setminus\{i_l,N\}}.
\end{align*} 
By Lemma 3-4(2) we have 
\begin{equation*}\tilde W_{I} 
=(-\sqrt{-1})(\sqrt{-1}\tilde X_{1,N}-(p+1))\tilde W_{I\setminus\{1,N\}}
+\sum\limits_{2\leq l\leq 2p+3}
(-1)^{l-1}(\tilde N^1_{i_l}\tilde W_{I\setminus\{i_l,N\}}).
\tag{3.7.1}
\end{equation*}
This yields the first assertion. Then we have by (3.7.1), 
\begin{align*}  
\tilde W_{I}^2 & =-(\sqrt{-1}\tilde X_{1,N}-(p+1))^2 
\tilde W_{I\setminus\{1,N\}}^2
\\& 
+(-\sqrt{-1})(\sqrt{-1}\tilde X_{1,N}-(p+1))
\tilde W_{I\setminus\{1,N\}}
\sum\limits_{2\leq l\leq 2p+3}(-1)^{l-1}
(\tilde N^1_{i_l}\tilde W_{I\setminus\{i_l,N\}})
\\& 
+\sum\limits_{2\leq l\leq 2p+3}(-1)^{l-1}
(\tilde N^1_{i_l}\tilde W_{I\setminus\{i_l,N\}})\tilde W_{I}. 
\tag{3.7.2}
\end{align*} 
By Lemma 3-6(1) for $q=l$ where $2\leq l\leq 2p+3$, we have 
\begin{equation*}
\tilde W_{I\setminus\{1,N\}}\tilde N^1_{i_l}
=\tilde N^1_{i_l}\tilde W_{I\setminus\{1,N\}}
+\sum\limits_{2\leq k\leq 2p+3,k\not= l}
\frac{l-k}{\vert l-k\vert}(-1)^{k+l}
\tilde N^1_{i_k}\tilde W_{I\setminus\{1,i_l,i_k,N\}}.
\tag{3.7.3}
\end{equation*}
Since we have 
$\sqrt{-1}\tilde X_{1,N}\tilde N^1_{i_s}=\tilde N^1_{i_s}
(\sqrt{-1}\tilde X_{1,N}+1)$ 
for $2\leq s\leq 2p+3$ by (3.2), we have by (3.7.3), 
\begin{align*} & \sqrt{-1}\tilde X_{1,N}\cdot 
\tilde W_{I\setminus\{1,N\}}\tilde N^1_{i_l}
=\tilde N^1_{i_l}(\sqrt{-1} X_{1,N}+1)\tilde W_{I\setminus\{1,N\}}
\\& 
+\sum\limits_{2\leq k\leq 2p+3,k\not= l}
\frac{l-k}{\vert l-k\vert}(-1)^{k+l}
\tilde N^1_{i_k}(\sqrt{-1} X_{1,N}+1)\tilde W_{I\setminus\{1,i_l,i_k,N\}}.
\tag{3.7.4}
\end{align*}
By (3.7.2),(3.7.3) and (3.7.4) we have 
\begin{align*}  
& \tilde W_{I}^2 =-(\sqrt{-1}\tilde X_{1,N}-(p+1))^2
\tilde W_{I\setminus\{1,N\}}^2
\\& 
+(-\sqrt{-1})\sum\limits_{2\leq l\leq 2p+3}
(-1)^{l-1}(\tilde N^1_{i_l}
((\sqrt{-1} \tilde X_{1,N}-p)\tilde W_{I\setminus\{1,N\}}
\tilde W_{I\setminus\{i_l,N\}})
\\& 
+(-\sqrt{-1})\sum\limits_{2\leq l\leq 2p+3}
(-1)^{l-1}
\\&
\times\sum\limits_{2\leq k\leq 2p+3,k\not= l}
\frac{l-k}{\vert l-k\vert}(-1)^{k+l}
\tilde N^1_{i_k}(\sqrt{-1} \tilde X_{1,N}-p)
\tilde W_{I\setminus\{1,i_l,i_k,N\}}
\tilde W_{I\setminus\{i_l,N\}})
\\& 
+\sum\limits_{2\leq l\leq 2p+3}(-1)^{l-1}(\tilde N^1_{i_l}
\tilde W_{I\setminus\{i_l,N\}})\tilde W_I. 
\end{align*} 
This yields the second assertion. 
\par 
$(2)$ Here we have $i_1>1$ and $i_{2p+4}=N$. 
By Proposition 3-2(2) and (3.1) we have 
\begin{align*}  & 
\tilde W_{I} =\sum\limits_{1\leq l\leq 2p+3}
(-1)^{l-1}\tilde X_{i_l,N}\tilde W_{I\setminus\{i_l,N\}}
\\& 
=\sum\limits_{1\leq l\leq 2p+3}(-1)^{l-1}
\tilde N^1_{i_l}\tilde W_{I\setminus\{i_l,N\}}
-\sqrt{-1}(\sum\limits_{1\leq l\leq 2p+3}(-1)^{l-1}\tilde X_{1,i_l}
\tilde W_{I\setminus\{i_l,N\}}).
\tag{3.7.5}
\end{align*} 
Here we put 
$$J=(j_1,j_2,\ldots,j_{2p+4})=(1,i_1,\ldots,i_{2p+3})
=I\cup\{1\}\setminus\{N\}.$$ 
By Proposition 3-2(1) for $q=1$ we have 
$$\tilde W_{J}
=\sum\limits_{2\leq k\leq 2p+4}
(-1)^{k}
\tilde X_{1,j_k}\tilde W_{J\setminus\{1,j_k\}}.$$ 
Hence we have 
$$\tilde W_{I\cup\{1\}\setminus\{N\}}
=\sum\limits_{1\leq l\leq 2p+3}
(-1)^{l-1}
\tilde X_{1,i_l}\tilde W_{I\setminus\{i_l,N\}}.$$ 
By (3.7.5) we have 
\begin{equation*}\tilde W_{I}=\sum\limits_{1\leq l\leq 2p+3}(-1)^{l-1}
(\tilde N^1_{i_l}\tilde W_{I\setminus\{i_l,N\}})
-(\sqrt{-1})\tilde W_{I\cup\{1\}\setminus\{N\}}
\tag{3.7.6}
\end{equation*}
This yields the first assertion. 
Then we have by (3.7.6) 
\begin{align*}  
\tilde W_{I}^2 & =-(\tilde W_{I\cup\{1\}\setminus\{N\}})^2
-(\sqrt{-1})\tilde W_{I\cup\{1\}\setminus\{N\}}
\sum\limits_{1\leq l\leq 2p+3}
(-1)^{l-1}(\tilde N^1_{i_l}\tilde W_{I\setminus\{i_l,N\}}) 
\\& 
+\sum\limits_{1\leq l\leq 2p+3}(-1)^{l-1}
(\tilde N^1_{i_l}\tilde W_{I\setminus\{i_l,N\}}))\tilde W_I. 
\tag{3.7.7}
\end{align*} 
By Lemma 3-6(2) for $q=l+1$ where $1\leq l\leq 2p+3$, we have 
\begin{align*} & 
\tilde W_{I\cup\{1\}\setminus\{N\}}\cdot\tilde N^1_{i_l}
=\tilde W_J\cdot \tilde N^1_{j_{l+1}}
=\tilde N^1_{j_{l+1}}\cdot \tilde W_J+(-1)^l
\tilde W_{J\cup\{N\}\setminus\{j_{l+1}\}}
\\& 
=\tilde N^1_{i_l}\cdot \tilde W_{I\cup\{1\}\setminus\{N\}}
+(-1)^l\cdot \tilde W_{I\cup\{1\}\setminus\{i_l\}}.
\end{align*}
Hence we have for $1\leq l\leq 2p+3$,
\begin{equation*}
\tilde W_{I\cup\{1\}\setminus\{N\}}\cdot\tilde N^1_{i_l}
=\tilde N^1_{i_l}\cdot\tilde W_{I\cup\{1\}\setminus\{N\}}
+(-1)^l\cdot\tilde W_{I\cup\{1\}\setminus\{i_l\}}.
\tag{3.7.8}
\end{equation*}
By (3.7.1) we have 
\begin{align*} & \tilde W_{I\cup\{1\}\setminus\{i_l\}}
=(-\sqrt{-1})(\sqrt{-1}\tilde X_{1,N}-(p+1))
\tilde W_{I\setminus\{i_l,N\}} 
\\& 
+\sum\limits_{1\leq k\leq 2p+3,k\not=l}
\frac{l-k}{\vert l-k\vert}(-1)^{k}
(\tilde N^1_{i_k}\tilde W_{I\cup\{1\}\setminus\{i_l,i_k,N\}}).
\tag{3.7.9}
\end{align*} 
By (3.7.8) and (3.7.9) we have 
\begin{align*} & 
\tilde W_{I\cup\{1\}\setminus\{N\}}\cdot\tilde N^1_{i_l}
=\tilde N^1_{i_l}\cdot\tilde W_{I\cup\{1\}\setminus\{N\}}
\\&
+(-1)^l(-\sqrt{-1})(\sqrt{-1}\tilde X_{1,N}-(p+1))
\tilde W_{I\setminus\{i_l,N\}} 
\\& 
+(-1)^l\sum\limits_{1\leq k\leq 2p+3,k\not=l}
\frac{l-k}{\vert l-k\vert}(-1)^{k}
(\tilde N^1_{i_k}\tilde W_{I\cup\{1\}\setminus\{i_l,i_k,N\}}).
\tag{3.7.10}
\end{align*} 
We substitute (3.7.10) into the second term in the right-hand side in (3.7.7) 
where we have 
$$(-\sqrt{-1})(-1)^{l-1}(-1)^l(-\sqrt{-1})
=1\text{ and }(-\sqrt{-1})(-1)^{l-1}(-1)^l=\sqrt{-1}.$$
Then we have 
\begin{align*} & 
\tilde W_{I}^2 
=-(\tilde W_{I\cup\{1\}\setminus\{N\}})^2
-\sqrt{-1}\sum\limits_{1\leq l\leq 2p+3}(-1)^{l-1}
(\tilde N^1_{i_l}\tilde W_{I\cup\{1\}\setminus\{N\}}
\tilde W_{I\setminus\{i_l,N\}})
\\&  
+(\sqrt{-1}\tilde X_{1,N}-(p+1))\sum\limits_{1\leq l\leq 2p+3}
(\tilde W_{I\setminus\{i_l,N\}})^2  
\\&
+\sqrt{-1}\sum\limits_{1\leq l\leq 2p+3}
\sum\limits_{1\leq k\leq 2p+3,k\not=l}
\frac{l-k}{\vert l-k\vert}(-1)^{k}
(\tilde N^1_{i_k}\tilde W_{I\cup\{1\}\setminus\{i_l,i_k,N\}})
\tilde W_{I\setminus\{i_l,N\}}
\\& 
+\sum\limits_{1\leq l\leq 2p+3}(-1)^{l-1}
(\tilde N^1_{i_l}\tilde W_{I\setminus\{i_l,N\}}))\tilde W_I.
\end{align*} 
This yields the second assertion. 
\par 
$(3)$ We have 
\begin{equation*}\tilde W_{p+1}
=\mathcal W_1+(\sum\limits_{1\in I,N\not\in I}\tilde W_I^2
+\sum\limits_{1\not\in I,N\not\in I}\tilde W_I^2)
\tag{3.7.11}
\end{equation*}
where we put   
$$\mathcal W_1=\sum\limits_{1,N\in I}\tilde W_I^2
+\sum\limits_{1\not\in I,N\in I}\tilde W_I^2.$$ 
By Proposition 3-7(1)(2) we have 
\begin{align*} & \mathcal W_1\in-(\sqrt{-1}\tilde X_{1,N}-(p+1))^2
\sum\limits_{1,N\in I} (\tilde W_{I\setminus\{1,N\}})^2
\\& 
+(\sqrt{-1}\tilde X_{1,N}-(p+1))\sum\limits_{1\not\in I,N\in I}
\sum\limits_{l=1}^{2p+3} (\tilde W_{I\setminus\{i_l,N\}})^2
-\sum\limits_{1\not\in I,N\in I} 
(\tilde W_{I\cup\{1\}\setminus\{N\}})^2
+\tilde{\mathfrak n}(0)_1\cdot U(\tilde{\mathfrak g}).
\end{align*}
Here we have 
$$\sum\limits_{1\not\in I,N\in I}
\sum\limits_{l=1}^{2p+3} (\tilde W_{I\setminus\{i_l,N\}})^2
=(N-2p-4)\sum\limits_{1\in I,N\in I} (\tilde W_{I\setminus\{1,N\}})^2.$$ 
Hence we have 
\begin{align*} & \mathcal W_1\in(-(\sqrt{-1}\tilde X_{1,N}-(p+1))^2)
\sum\limits_{1\in I,N\in I} (\tilde W_{I\setminus\{1,N\}})^2
\\& 
+(N-2p-4)(\sqrt{-1}\tilde X_{1,N}-(p+1))
\sum\limits_{1,N\in I} (\tilde W_{I\setminus\{1,N\}})^2)
-\sum\limits_{1\not\in I,N\in I} (\tilde W_{I\cup\{1\}\setminus\{N\}})^2
+\tilde{\mathfrak n}(0)_1\cdot U(\tilde{\mathfrak g}).
\end{align*}
We have 
\begin{align*} & -(\sqrt{-1}\tilde X_{1,N}-(p+1))^2
+(N-2p-4)(\sqrt{-1}\tilde X_{1,N}-(p+1))
\\& 
=-(\sqrt{-1}\tilde X_{1,N}-(p+1))(\sqrt{-1}\tilde X_{1,N}-(N-p-3)).
\end{align*}
Then we have 
\begin{align*} & \mathcal W_1\in(-(\sqrt{-1}\tilde X_{1,N}-(p+1))(\sqrt{-1}\tilde X_{1,N}-(N-p-3))
\sum\limits_{1\in I,Nin I} (\tilde W_{I\setminus\{1,N\}})^2
-\sum\limits_{1\not\in I,N\in I} (\tilde W_{I\cup\{1\}\setminus\{N\}})^2.
\end{align*}
Here we have  
$$\sum\limits_{1\not\in I,N\in I} (\tilde W_{I\cup\{1\}\setminus\{N\}})^2
=\sum\limits_{1\in I,N\not\in I}\tilde W_I^2.$$
Then we have by (3.7.11), 
\begin{align*} & 
\tilde W_{p+1}\in(-(\sqrt{-1}\tilde X_{1,N}-(p+1))
(\sqrt{-1}\tilde X_{1,N}-(N-p-3))
\sum\limits_{1\in I,N\in I}(\tilde W_{I\setminus\{1,N\}})^2
\\& 
+\sum\limits_{1\not\in I,N\not\in I}\tilde W_I^2 
+\tilde{\mathfrak n}(0)_1\cdot U(\tilde{\mathfrak g}).
\end{align*}
This yields the assertion. \qed
\end{Proof}
\par 
\textbf{Remark}
In view of $(1.4)$ the assertions of 
Proposition 3-7(1)(2)(3) in the setting of $n=1$ and $p=0$ 
follow from Proposition 2-3 (2.3.1)(2.3.2)(2.3.3) in $[23]$. 
\par 
\textbf{Proposition 3-8(The Calculations of Pfaffians 
of $\tilde{\mathfrak g}$)}
\par 
$(1)$ We have 
\begin{align*} & \tilde W_{1}=\sum\limits_{1\leq i<j\leq N}
\tilde X_{i,j}^2
\in(-\sum\limits_{1\leq l\leq [N/2]}
(\sqrt{-1}\tilde X_{l,N+1-l})
((\sqrt{-1}\tilde X_{l,N+1-l})-(N-2l))))
+\tilde{\mathfrak n}(0)\cdot U(\tilde{\mathfrak g}).
\end{align*} 
\par 
$(2)$ Assume that $p\geq 0,2p+4\leq N$.  
We have 
\begin{align*} & \tilde W_{p+2}
\in(-1)^{p}(\sum\limits_{J\in T_{p+1}([N/2]-1)}
\\& 
\{(\prod\limits_{1\leq l\leq p+1}
((\sqrt{-1} X_{j_l,N+1-j_l}-(p+2-j_l))
(\sqrt{-1} X_{j_l,N+1-j_l}-(N-4-p+j_l))))
\\& 
\times
(\sum\limits_{1\leq l\leq [(N-2j_{p+1})/2] }
((\sqrt{-1}\tilde X_{j_{p+1}+l,N+1-j_{p+1}-l})
(\sqrt{-1}\tilde X_{j_{p+1}+l,N+1-j_{p+1}-l}-(N-2j_{p+1}-2l)))\}
\\& 
+\tilde{\mathfrak n}(0)\cdot U(\tilde{\mathfrak g}).
\end{align*} 
\begin{Proof} 
$(1)$ By (3.1) we have  
\begin{align*} & \tilde W_{1}=\sum\limits_{1\leq i<j\leq N}
\tilde X_{i,j}^2.
=\sum\limits_{1\leq k\leq N-1}\tilde X_{1,N}^2
-\sqrt{-1}\sum\limits_{2\leq l\leq N-1}
\tilde X_{1,l}\tilde X_{l,N}
+\sum\limits_{2\leq l\leq N-1}
\tilde N^1_{l}\tilde X_{l,N}
+\sum\limits_{2\leq i<j\leq N}\tilde X_{i,j}^2.  
\end{align*} 
Here we have  
$$\tilde X_{1,l}\tilde X_{l,N} 
=\tilde X_{l,N} \tilde X_{1,l}-\tilde X_{1,N}.$$ 
Hence we have by (3.1), 
\begin{align*} & 
\tilde W_1=\sum\limits_{1\leq k\leq N-1}\tilde X_{1,N}^2
-\sqrt{-1}(\sum\limits_{2\leq l\leq N-1}
\tilde X_{l,N} \tilde X_{1,l})
+\sqrt{-1} (N-2)\tilde X_{1,N}
+\sum\limits_{2\leq i<j\leq N}\tilde X_{i,j}^2
+\sum\limits_{2\leq l\leq N-1}\tilde N^1_{l}\tilde X_{l,N}
\\& 
=\sum\limits_{1\leq k\leq N-1}\tilde X_{1,N}^2
-\sqrt{-1}(\sum\limits_{2\leq l\leq N-1}
\tilde N^1_l\tilde X_{1,l}
-\sum\limits_{2\leq l\leq N-1}\sqrt{-1}\tilde X_{1,l}^2)
+\sqrt{-1} (N-2)\tilde X_{1,N}
\\& 
+\sum\limits_{2\leq i<j\leq N}\tilde X_{i,j}^2
+\sum\limits_{2\leq l\leq N-1}\tilde N^1_{l}\tilde X_{l,N}
\\& 
=\sum\limits_{1\leq k\leq N-1}\tilde X_{1,N}^2
-\sum\limits_{2\leq l\leq N-1}\tilde X_{1,l}^2 
+\sqrt{-1} (N-2)\tilde X_{1,N}
+\sum\limits_{2\leq i<j\leq N}\tilde X_{i,j}^2
\\& 
+\sum\limits_{2\leq l\leq N-1}\tilde N^1_{l}\tilde X_{l,N}
-\sqrt{-1}(\sum\limits_{2\leq l\leq N-1}
\tilde N^1_l\tilde X_{1,l})
\\& 
=X_{1,N}^2+\sqrt{-1} (N-2)\tilde X_{1,N}
+\sum\limits_{2\leq i<j\leq N}\tilde X_{i,j}^2
+\sum\limits_{2\leq l\leq N-1}\tilde N^1_{l}\tilde X_{l,N}
-\sqrt{-1}(\sum\limits_{2\leq l\leq N-1}
\tilde N^1_l\tilde X_{1,l}). 
\end{align*} 
Therefore we have 
$$\tilde W_{1}\in 
-(\sqrt{-1}\tilde X_{1,N})(\sqrt{-1}\tilde X_{1,N})-(N-2))
+\sum\limits_{2\leq i<j\leq N}\tilde X_{i,j}^2
+\tilde{\mathfrak n}(0)\cdot U(\tilde{\mathfrak g}).$$ 
By induction we have 
$$\tilde W_{1}\in -\sum\limits_{1\leq l\leq [N/2]}
(\sqrt{-1}\tilde X_{l,N+1-l}))
((\sqrt{-1}\tilde X_{l,N+1-l})-(N-2l))
+\tilde{\mathfrak n}(0)\cdot U(\tilde{\mathfrak g}).$$ 
This yields the assertion. 
\par  
$(2)$ We have only to prove the following (3.8.1), 
since the assertion  follows from (3.8.1) 
and the assertion of Proposition 3-8(1)
\begin{align*} & \tilde W_{p+2}\in 
(-1)^{p+1}\sum\limits_{J\in T_{p+1}([N/2]-1)}
\\&
(\prod\limits_{1\leq l\leq p+1}
((\sqrt{-1} X_{j_l,N+1-j_l}-(p+2-j_l))
(\sqrt{-1} X_{j_l,N+1-j_l}-(N-4-p+j_l))))
\\& 
\times
\sum\limits_{j_{p+1}+1\leq l_1<l_2\leq N-j_{p+1}}
\tilde X_{l_1,l_2}^2)
+\tilde{\mathfrak n}(0)\cdot U(\tilde{\mathfrak g}).
\tag{3.8.1}
\end{align*} 
We prove the assertion of (3.8.1) by the induction on $N$. 
The assertion of (3.8.1) when $N\leq 5$ follows directly 
from Proposition 3-7(3) and Proposition 3-8(1). 
\par 
Assume that $N\geq 6$. In the setting of Proposition 3-7(3), 
we have the following (3.8.2) and (3.8.3) 
by the inductive hypothesis. 
\begin{align*} 
& \sum\limits_{1,N\in I}\tilde W_{I\setminus\{1,N\}}^2
\in(-1)^p\sum\limits_{J\in T_{p+1}([N/2]-2),j_1=1}
\\&
(\prod\limits_{2\leq l\leq p+1}
((\sqrt{-1} X_{j_l,N+1-j_l}-(p+2-j_l))
(\sqrt{-1} X_{j_l,N+1-j_l}-(N-4-p+j_l))))
\\& 
\times
\sum\limits_{j_{p+1}+1\leq l_1<l_2\leq N-j_{p+1}}
\tilde X_{l_1,l_2}^2)
+\tilde{\mathfrak n}(0)\cdot U(\tilde{\mathfrak g}).
\tag{3.8.2}
\end{align*} 

\begin{align*} & 
\sum\limits_{1\not\in I,N\not\in I}\tilde W_I^2
\in(-1)^{p+1}\sum\limits_{J\in T_{p+1}([N/2]-2),j_1\not=1}
\\&
(\prod\limits_{1\leq l\leq p+1}
((\sqrt{-1} X_{j_l,N+1-j_l}-(p+2-j_l))
(\sqrt{-1} X_{j_l,N+1-j_l}-(N-4-p+j_l))))
\\& 
\times
\sum\limits_{j_{p+1}+1\leq l_1<l_2\leq N-j_{p+1}}
\tilde X_{l_1,l_2}^2)
+\tilde{\mathfrak n}(0)\cdot U(\mathfrak g).
\tag{3.8.3}
\end{align*} 
We substitute (3.8.2) and (3.8.3) into Proposition 3-7(3).
Here we have 
$$U(\tilde{\mathfrak h})\tilde{\mathfrak n}(0)
\subset(\tilde{\mathfrak n}(0))U(\tilde{\mathfrak h}).$$
Then the assertion (3.8.1) follows immediately. \qed 
\end{Proof}
\par 
Assume that $\lambda\in\tilde{\mathfrak h}^*$ satisfies that 
$$\lambda(\sqrt{-1} X_{l,N+1-l})
\in\mathbb Z\text{ for }1\leq l\leq [N/2].$$
If $N$ is odd, we assume that 
\begin{equation*}
\lambda(\sqrt{-1}\tilde X_{1,N}))
\geq\cdots
\geq\lambda(\sqrt{-1}\tilde X_{l,N+1-1}))\cdots
\geq\lambda(\sqrt{-1}\tilde X_{[N/2],N+1-[N/2]}))\geq 0.
\tag{3.9.0.1}
\end{equation*}
If $N$ is even, we assume that 
\begin{equation*}
\lambda(\sqrt{-1}\tilde X_{1,N}))
\geq\cdots
\geq\lambda(\sqrt{-1}\tilde X_{l,N+1-1}))\cdots
\geq\vert\lambda(\sqrt{-1}\tilde X_{[N/2],[N/2]+1}))
\vert\geq 0.
\tag{3.9.0.2}
\end{equation*}
Then it follows from Theorem 5.110 in [30] and 
18 and 20 in IV 9 in [30] that 
uniquely determined up to equivalence classes, 
there exists a finite dimensional irreducible unitary representation 
$(\pi_\lambda, V_\lambda)$ of $SO(N,\mathbb R)$ 
with the highest weight $\lambda\in\mathfrak h^*$. 
Moreover, any finite dimensional irreducible unitary representation of $SO(N,\mathbb R)$ which is isomorphic to 
$(\pi_\lambda,V_\lambda)$ 
for $\lambda\in\tilde{\mathfrak h}^*$  
satisfies (3.9.0.1) and (3.9.0.2). 
Let $v_\lambda\in V_\lambda$ be a nonzero highest weight vector 
in $V_\lambda$. Then we have 
\begin{equation*} 
d\pi(X).v_\lambda=0
\text{ for }X\in\tilde{\mathfrak n}(0),
\tag{3.9.0.3}
\end{equation*} 
where $d\pi$ is the differential representation of $\pi$. 
Let $\langle*,*\rangle_\lambda$ be a nonzero 
$SO(N,\mathbb R)$-invariant inner product on $V_\lambda$. 
For $v\in V_\lambda$ we define a function 
on $SO(N,\mathbb R)$ by 
$$F^\lambda_v(g)=\langle g.v_\lambda,v\rangle_\lambda
\text{ for }g\in SO(N,\mathbb R).$$
\par 
\textbf{Proposition 3-9(The the eigenvalue 
of the central elements of 
$\tilde{\mathfrak g}$)}
\par 
Assume that $\lambda\in\tilde{\mathfrak h}^*$ satisfies 
$(3.9.0.1)$ or $(3.9.0.2)$. 
\par 
$(1)$ Assume that $1\leq r\leq [N/2]$. 
Then there exists some $D^N_r(\lambda)\in\mathbb C$ such that 
$$\tilde W_r.F^\lambda_v(g)=D^N_r(\lambda)F^\lambda_v(g)
\text{ for }g\in SO(N,\mathbb R),v\in V_\lambda.$$
\par 
$(2)$ Assume that $p\geq 0$ and $2p+4\leq N$. 
The following assertions (i) and (ii) hold. 
\par 
(i) $D^N_1(\lambda)=0$  
if and only if $\lambda=0$. 
\par 
(ii) $D^N_{p+2}(\lambda)=0$ 
if and only if  
$$\lambda(\sqrt{-1}\tilde X_{l,N+1-l}))=0
\text{ for }p+2\leq l\leq [N/2].$$
\begin{Proof}
Since we have 
$\tilde W_r\in\mathcal Z(\mathfrak s\mathfrak o(N;\mathbb R))$, 
we have for $g\in SO(N,\mathbb R)$ and $v\in V_\lambda$, 
\begin{equation*}
\tilde W_r.F^\lambda_v(g)=(Ad(g)\tilde W_r).F^\lambda_v(g)
=(Ad(g)\tilde W_r).\langle\pi_\lambda(g).v_\lambda,v\rangle.
\tag{3.9.1}
\end{equation*}
Assume that $r=1$. 
It follows from (3.9.1), Proposition 3-8(1) 
and (3.9.0.3) that 
we have for $g\in SO(N,\mathbb R)$ and $v\in V_\lambda$, 
$$\tilde W_1.(F^\lambda_v(g))
=D^N_1(\lambda)F^\lambda_v(g),$$ 
where we put 
\begin{align*} & D^N_1(\lambda)
=-\sum\limits_{1\leq l\leq [N/2]}
(-\lambda(\sqrt{-1}\tilde X_{l,N+1-l}))
((-\lambda(\sqrt{-1}\tilde X_{l,N+1-l}))-(N-2l)))
\\& 
=-\sum\limits_{1\leq l\leq [N/2]}
(\lambda(\sqrt{-1}\tilde X_{l,N+1-l}))
((\lambda(\sqrt{-1}\tilde X_{l,N+1-l}))+(N-2l)). 
\end{align*} 
In view of (3.9.0.1) or (3.9.0.2) we have 
$D^N_1(\lambda)=0$ if and only if 
$$\lambda(\sqrt{-1}\tilde X_{l,N+1-l}))=0
\text{ for }1\leq l\leq [N/2],$$
which is equivalent to $\lambda=0$. 
This yields the assertion of (1) for $r=1$ and (2)(i).  
\par 
Assume that $p\geq 0,r=p+2$ and $p\geq 0$ and $2p+4\leq N$. 
It follows from (3.9.1), Proposition 3-8(2) 
and (3.9.0.3) that 
we have for $g\in SO(N,\mathbb R)$ and $v\in V$, 
$$\tilde W_{p+2}.(F^\lambda_v(g))=D^N_{p+2}(\lambda)F^\lambda_v(g).$$  
where we put 
\begin{align*} & D^N_{p+2}(\lambda)
=(-1)^{p}(\sum\limits_{J\in T_{p+1}([N/2]-1)}
(\prod\limits_{1\leq l\leq p+1}
\\& 
((-\lambda(\sqrt{-1}\tilde X_{j_l,N+1-j_l})-(p+2-j_l))
(-\lambda(\sqrt{-1} \tilde X_{j_l,N+1-j_l})-(N-4-p+j_l)))
\\& 
\times
(\sum\limits_{1\leq l\leq [(N-2j_{p+1})/2] }
((-\lambda(\sqrt{-1}\tilde X_{j_{p+1}+l,N+1-j_{p+1}-l}))
(-\lambda(\sqrt{-1}\tilde X_{j_{p+1}+l,N+1-j_{p+1}-l})-(N-2j_{p+1}-2l)))
\\& 
=(-1)^{p}(\sum\limits_{J\in T_{p+1}([N/2]-1)}
(\prod\limits_{1\leq l\leq p+1}
\\& 
((\lambda(\sqrt{-1}\tilde X_{j_l,N+1-j_l})+(p+2-j_l))
(\lambda(\sqrt{-1}\tilde X_{j_l,N+1-j_l})+(N-4-p+j_l))
\\& 
\times
(\sum\limits_{1\leq l\leq [(N-2j_{p+1})/2] }
((\lambda(\sqrt{-1}\tilde X_{j_{p+1}+l,N+1-j_{p+1}-l}))
(\lambda(\sqrt{-1}\tilde X_{j_{p+1}+l,N+1-j_{p+1}-l})+(N-2j_{p+1}-2l))).
\end{align*}
In view of (3.9.0.1) or (3.9.0.2) we have $D^N_{p+2}(\lambda)=0$ 
if and only if 
if and only if 
$$\lambda(\sqrt{-1}\tilde X_{l,N+1-l})=0
\text{ for }p+2\leq l\leq [N/2].$$ 
This yields the assertion of (1) for $r\geq 2$ and (2)(ii). 
\qed 
\end{Proof}
\par 
\textbf{Remark} 
Let $1\leq r\leq [N/2]$. 
We have to calculate the projection of $\tilde W_{r}$ 
mod $U(\tilde{\mathfrak g})\tilde{\mathfrak n}(0)$ 
in order to directly calculate 
$C^N_r(\lambda)\in\mathbb C$ such that 
$d\pi(\tilde W_{r}).v_\lambda=C^N_r(\lambda)\cdot v_\lambda$. 
However it turns out that $C^N_r(\lambda)=D^N_r(\lambda)$. 
The explicit calculation of $C^N_1=D^N_1$ can be found in $[7]$. 
\par 
Let $(\pi_\lambda,V_\lambda)$ be an irreducible finite dimensional 
unitary representation of $SO(N,\mathbb R)$ 
with the highest weight $\lambda\in\tilde{\mathfrak h}^*$. 
We define a subgroup $H^N_k$ of $SO(N,\mathbb R)$ 
for $1\leq k<N$ by 
\begin{equation*} H^N_k
=\{ u\in SO(N,\mathbb R) 
\bigm\vert 
u=\begin{pmatrix} 
u_1 & 0 
\\ 
0 & u_2 \end{pmatrix} ,
u_1\in O(N-k),u_2\in O(k) \}.\end{equation*} 

Then we have 
$$\min(k_1,N-k_1)<\min(k_2,N-k_2)  
\text{ if and only if }
\dim(SO(N,\mathbb R)/H^N_{k_1})
<\dim(SO(N,\mathbb R)/H^N_{k_2}).$$ 
\par 
\textbf{Proposition 3-10}
\par 
Assume that $1\leq k_1,k_2<N$ 
and $\min(k_1,N-k_1)<\min(k_2,N-k_2)$. 
Let of $(\pi_\lambda,V_\lambda)$ be an irreducible finite dimensional unitary representation 
of $SO(N,\mathbb R)$. 
\par 
$(1)$ If there exists a nonzero $H^N_{k_1}$-fixed element in $V_\lambda$, 
then there exists a nonzero $H^N_{k_2}$-fixed element in $V_\lambda$. 
\par 
$(2)$ If there exists a a nonzero $H^N_{k_2}$-fixed element in $V_\lambda$ 
and $D^N_{\min(k_1,N-k_1)+1}(\lambda)=0$, 
then there exists a nonzero $H^N_{k_1}$-fixed element in $V_\lambda$. 
\begin{Proof}
Let $1\leq k<N$. It follows from  \S 6 in [42] that 
there exists a nonzero $H^N_k$-fixed element in $V_\lambda$ 
if and only if 
$\lambda\in\mathfrak h^*$ satisfies (3.9.0.1),(3.9.0.2) and 
\begin{align*} & \lambda(\sqrt{-1}\tilde X_{l,N+1-l})
\in 2\mathbb Z\text{ for }1\leq l\leq N
\text{ and }
\lambda(\sqrt{-1} X_{l,N+1-l})=0
\text{ for }\min(k,N-k)+1\leq l\leq N.
\tag{3.10}
\end{align*} 
The assertion (1) directly follows from (3.10).   
The assertion (2) follows from (3.10) and Proposition 3-9. 
\qed 
\end{Proof}
\section{ The Properties of 
the Invariant Differential Operators.} 
In this section 
we investigate the properties of the Pfaffian type elements 
for $\mathfrak s\mathfrak o(m,n)$  
with respect to the Fourier transform on $M$. 
\par 
We put for $1\leq q\leq n$, 
$$\mathfrak n(0)_q
=(\bigoplus\limits_{q< i\leq n}\mathbb R(Z^+_{q,i}-Z^-_{q,i}))
\oplus
(\bigoplus\limits_{n+1\leq i\leq m}W_{i,m+n+1-q})
\oplus 
(\bigoplus\limits_{q< i\leq n}\mathbb R(Z^+_{q,i}+Z^-_{q,i})).$$ 
Then we have 
$$\mathfrak n(0)=\bigoplus\limits_{1\leq q\leq n}\mathfrak n(0)_q.$$
We put for $1\leq q\leq n$,
$$N^q_l   
=\begin{cases} (Z_{q,l}^--Z_{q,l}^+)/2 & \text{ for }q<l\leq n \\
        W_{l,m+n+1-q} & \text{ for } n+1\leq l\leq m \\ 
(Z_{q,m+n+1-l}^++Z_{q,m+n+1-l}^-)/2  & 
\text{ for }m+1\leq l< m+n+1-q. 
\end{cases}$$ 
Then we have for $1\leq q\leq n$, 
\begin{equation*}N^q_l   
=\begin{cases} Y_{l,m+n+1-q}+X_{l,q} & \text{ for }q< l\leq m \\
        X_{l,m+n+1-q}+Y_{l,q}  & 
\text{ for }m+1\leq l< m+n+1-q. 
\end{cases}
\end{equation*}
We have for $1\leq q\leq n$, 
$$\mathfrak n(0)_q
=\bigoplus\limits_{q< l< m+n+1-q}\mathbb R N^q_l.$$
and 
\begin{equation*}
[Y_{q,m+n+1-q},N^q_l]=N^q_l\text{ for }q<l< m+n+1-q.
\tag{4.1}
\end{equation*}
We put for $1\leq q\leq n$, 
\begin{equation*}
N^q_{l}   
=\begin{cases}  (-\sqrt{-1})\cdot N^q_{l}  &  
 \text{ for } q< l\leq m \\ 
          N^q_{l}  
        & \text{ for }m+1\leq l< m+n+1-q. 
\end{cases}
\end{equation*}
We put for $1\leq i<j\leq m+n$, 
\begin{equation*}{\overline X}_{ij}   
=\begin{cases}  X_{ij}  &  \text{ if } 1\leq i<j\leq m\text{ or }m+1\leq i<j\leq m+n\\ 
        -\sqrt{-1}\cdot Y_{ij}  
        & \text{ if }1\leq i\leq m,m+1\leq j\leq m+n. 
\end{cases} 
\end{equation*}
We put for $1\leq j<i\leq m+n$, 
$${\overline X}_{ij}=-{\overline X}_{ji}.$$   
Then we have for $1\leq q\leq n$, 
\begin{equation*}\overline N^q_l=\overline X_{l,m+n+1-q}
+\sqrt{-1}\cdot\overline X_{q,l}
\in(\mathfrak n(0)_q)_{\mathbb C} 
\text{ for }q< l < m+n+1-q.
\tag{4.2}
\end{equation*}
Assume that $m\geq n\geq 1$ and $p\geq 1,2p\leq m+n$. 
We put for $I=(i_1,\cdots,i_{2p})\in T_{2p}(m+n)$, 
Assume that $m\geq n\geq 1$ and $p\geq 1,2p\leq m+n$. 
We put for $I=(i_1,\cdots,i_{2p})\in T_{2p}(m+n)$, 
$${W}_{I} 
=\sum_{\sigma\in P_+(I)} 
\epsilon(\sigma){\overline X}_{\sigma(i_1)\sigma(i_2)}
\cdots {\overline X}_{\sigma(i_{2p-1}) \sigma(i_{2p}) }$$ 
and 
$${W}_p=\sum\limits_{I\in T_{2p}(m+n)}{W}_{I}^2.$$ 
We put $W_\emptyset=1$.   
\par 
\textbf{Proposition 4-1(Differential Operators on the Symmetric Spaces)}
\par 
Assume that $m\geq n\geq 1$. Suppose that $p\geq 1,2p\leq m+n$. 
\par 
$(1)$ There exists some $a_{I'}^{I}(g)\in\mathbb C$ 
for all $I,I'\in T_{2p}(m+n)$ and $g\in G$ such that 
we have for $I\in T_{2p}(m+n)$, 
$$Ad(g){W}_{I} 
=\sum\limits_{I'\in T_{2p}(m+n)}a_{I}^{I'}(g)\cdot {W}_{I'}
\text{ for }g\in G.$$ 
\par 
$(2)$ We have ${W}_p\in \mathcal Z(\mathfrak g)$. 
\par
$(3$ Assume $m\geq r\geq 0,p+1>\min(m-r,r+n)$ 
and $I\in T_{2p+2}(m+n)$. 
Then the elements ${W}_{I}$ and ${W}_{p+1}$ are contained in both 
$$Ann_{U(\mathfrak g_{\mathbb C})}(C^\infty(G/H^{(r)}))
\text{ and }Ann_{U(\mathfrak g_{\mathbb C})}(C^\infty(G/H^{(m-n-r)})).$$
\begin{Proof}
$(1)$ We assume the setting of Fact 3-1 
for $N=m+n$ and $\mathbb F=\mathbb C$. 
Let $E:\mathfrak s\mathfrak o(m,n)\rightarrow
\mathfrak s\mathfrak o(m+n;\mathbb C)$ 
be the linear map defined by $E(X_{ij})={\tilde X}_{ij}$ 
for $1\leq i<j \leq m,m+1\leq i<j\leq m+n$ 
and $E(Y_{ij})=\sqrt{-1}\cdot {\tilde X}_{ij}$ 
for $1\leq i \leq m,m+1\leq j\leq m+n$.
Then $E$ is an injective Lie algebra homomorphism. 
We also denote by $E$ the induced algebra homomorphism 
from $U(\mathfrak s\mathfrak o(m,n))\otimes_{\mathbb R}\mathbb C$ 
into $U(\mathfrak s\mathfrak o(m+n,\mathbb C))
\otimes_{\mathbb R}\mathbb C$. 
Let $\tilde E:Ad(SO_0(m,n))\rightarrow Ad(SO(m+n;\mathbb C))$ 
be the Lie group homomorphism naturally induced by $E$ 
such that $E(g.X)={\tilde E}(g).E(X)$ 
for all $g\in Ad(SO_0(m,n))$ and all $X\in U(\mathfrak g)$. 
For all $g\in SO_0(m,n)$, we can choose some 
$\tilde g\in SO(m+n;\mathbb C)$ 
such that $Ad(\tilde g)=\tilde E(Ad(g))$. 
For all $g\in SO_0(m,n)$, we can choose some 
$\tilde g\in SO(m+n;\mathbb C)$ 
such that $Ad(\tilde g)=\tilde E(Ad(g))$. 
We have 
\begin{equation*}E(\overline X_{ij})=\tilde X_{ij}\text{ for }
1\leq i,j\leq m+m,i\not=j.
\tag{4.1.1}
\end{equation*}
Then we have  
\begin{equation*}
E(W_I)=\tilde W_I\text{ for }I\in T_{2p}(m+n) 
\tag{4.1.2]}
\end{equation*}
Then we have by Fact 3-1(1) and (4.1.2), 
\begin{align*}
& E(Ad(g) W_{I})=\tilde E(Ad(g))E(W_{I})
=Ad(\tilde g)\tilde W_{I}
\\& 
=\sum_{I'\in T_{2p}(m+n)}a_{I}^{I'}(g)\tilde W_{I'}
=\sum_{I'\in T_{2p}(m+n)}a_{I}^{I'}(g)E(W_{I'})
=E(\sum_{I'\in T_{2p}(m+n)}a_{I}^{I'}(g)W_{I'})
\end{align*} 
where we put for $I'=({i'}_1,\cdots,{i'}_{2p})\in T_{2p}(m+n)$, 
$$\tilde g=(\tilde g_{i,j})\text{ and }
a_{I}^{I'}(g)=\begin{vmatrix} 
\tilde g_{i_1,{i'}_1} & \cdots & \tilde g_{i_{2p},{i'}_{1}} \\ 
\vdots & \ddots & \vdots \\
\tilde g_{i_{2p},{i'}_{1}} & \cdots & \tilde g_{i_{2p},{i'}_{2p}} 
\end{vmatrix}.$$ 
Thus the assertion follows from the injectivity of $E$. 
\par 
$(2)$ We  have 
$$E({W}_p)=E(\sum\limits_{I\in T_{2p}(m+n)}W_{I}^2)
=\sum\limits_{I\in T_{2p}(m+n)}\tilde W_{I}^2.$$ 
The right-hand side belongs to 
$\mathcal Z(\mathfrak s\mathfrak o(m+n;\mathbb C))$ 
by Fact 3-1(2) for $N=m+n$ and $\mathbb F=\mathbb C$. 
Then it also belongs to $\mathcal Z(E(\mathfrak g))$. 
Thus we have ${W}_p\in\mathcal Z(\mathfrak g)$. 
This yields the assertion. 
\par
$(3)$ By Proposition 4-1(1) we have for $I\in T_{2p+2}(m+n)$ and 
$f\in C^\infty(G/H^{(r)})$, 
\begin{align*} 
(& W_{I}.f)(g H^{(r)})
=((Ad(g)Ad(g^{-1}) W_{I}).f)(g H^{(r)})
\\
&=\sum_{I'\in T_{2p+2}(m+n)}
a_{I}^{I'}(g^{-1})((Ad(g) W_{I'}).f)(g H^{(r)})=0.
\end{align*} 
Then by (4.1.1) 
we have $W_{I}\in Ann_{U(\mathfrak g_\mathbb C)}C^{\infty}(G/H^{(r)})$, since at least one of the $(p+1)$ factors 
in the summands in (1.5) 
is an element of the Lie algebra of the isotropy group $H^{(r)}$ 
in view of the rank of $G/H^{(r)}$ 
which is equal to $\min(m-r,n+r)$. Similarly, we have 
$W_{I}\in Ann_{U(\mathfrak g_\mathbb C)} C^{\infty}(G/H^{(m-n-r)})$ 
in view of the rank of $G/H^{(m-n-r)}$ which is equal to $\min(m-r,n+r)$. 
Then $W_{I}$ is contained in both 
$Ann_{U(\mathfrak g_\mathbb C)}C^{\infty}(G/H^{(r)})$ 
and $Ann_{U(\mathfrak g_\mathbb C)}C^{\infty}(G/H^{(m-n-r)})$. 
\qed
\end{Proof}
\par 
\textbf{Remark}
The assertions in Proposition 4-1 have been obtained in 
Proposition 2-2 in $[25]$
in the setting of $n=1$ 
in Proposition 2-2 in $[23]$
in the setting of $n=1,r=0$. 
\par 
\textbf{Proposition 4-2(Differential Operators 
on the Horocycle Spaces))}
\par 
\par 
Assume that $n\geq 1; 0\leq r_1,r_2 \leq m-n ; r_1\not= r_2$. 
and $p\geq 0,2p+4\leq m+n$. 
Suppose that $p+2>\min(m-r,r+n)$ and $I\in T_{2p+4}(m+n)$. 
Then the elements ${W}_{I}$ and ${W}_{p+2}$ are contained in both 
$$Ann_{U(\mathfrak g_\mathbb C)}(C^\infty(G/((M\cap H^{(r)}) N))
\text{ and }
Ann_{U(\mathfrak g_\mathbb C)}(C^\infty(G/((M\cap H^{(m-n-r)})N)).$$
\begin{Proof}
By Proposition 2-1 and Proposition 4-1(1) 
we have only to prove the assertions 
when $M=M(0),A=A(0)$ and $N=N(0)$. 
We have 
\begin{equation*}\#(I)^n_-+\#(I)^n_+\leq n+\#(I)^n_1
\text{ and }
\#(I)^n_0=2p+4-(\#(I)^n_-+\#(I)^n_+).
\tag{4.2.1}
\end{equation*}
By Lemma 3-5(2) with (4.1.1) and $N=m+n$, we have 
\begin{align*} & W_I\in
\sum\limits_{J\subset I,
\#(J)^n_-+\#(J)^n_+=(\#(I)^n_-+\#(I)^n_+)-2\#(I)^n_1,
\#(J)^n_0=\#(I)^n_0,\#(J)^n_1=0}
W_J U(\mathfrak g_\mathbb C)+\mathfrak n(0)_\mathbb C U(\mathfrak g_\mathbb C).
\end{align*} 
Here we have 
\begin{align*} & \#(J)^n_-+\#(J)^n_+
=(\#(I)^n_-+\#(I)^n_+)-2\#(I)^n_1 
\\& 
\text{ and } 
\#(J)^n_0=\#(I)^n_0=2p+4-(\#(I)^n_-+\#(I)^n_+).
\tag{4.2.2}
\end{align*}  
By (4.2.1) and (4.2.2) we have  
\begin{align*} & \#(J)^n_0-(\#(J)^n_-+\#(J)^n_+)
=2p+4-2(\#(I)^n_-+\#(I)^n_+)+2\#(I)^n_1
\\& 
\geq(2p+4)-2(n+\#(I)^n_1)+2\#(I)^n_1
=2(p+2-n)>2\min(m-n-r,r).
\end{align*} 
Hence we have 
\begin{equation*}W_I\in
\sum\limits_{J\subset I,\#(J)^n_0-(\#(J)^n_-+\#(J)^n_+)
>2\min(m-n-r,r)}
W_J U(\mathfrak g_\mathbb C)
+\mathfrak n(0)_\mathbb C U(\mathfrak g_\mathbb C).
\tag{4.2.3} 
\end{equation*}
In the setting of (4.2.3), we have by (2.2), 
\begin{align*} & ((Ad(g)W_J).f_0)(g((M(0)\cap H^{(r)})N(0)))=0
\\& 
\text{ and }
((Ad(g)W_J).f_1)(g((M(0)\cap H^{(m-n-r)})N(0)))=0
\text{ for }g\in G,
\end{align*}
where $f_0\in C^\infty(G/((M(0)\cap H^{(r)}) N(0)))$ 
and $f_1\in C^\infty(G/((M(0)\cap H^{(m-n-r)}) N(0)))$, 
since at least one of the factors 
in the summands of $W_J$ in (1.5) 
is an element of the Lie algebra 
of $M(0)\cap H^{(r)}$ or $M(0)\cap H^{(m-n-r)}$ 
in view of the rank of $M(0)/(M(0)\cap H^{(r)})$ 
or $M(0)/(M(0)\cap H^{(m-n-r)})$ 
which is equal to $\min(m-n-r,r)$. 
Hence by (4.2.3) we have for any $I\in T_{2p+4}(m+n)$ and $g\in G$, 
\begin{align*} & ((Ad(g)W_I).f_0)(g((M(0)\cap H^{(r)})N(0)))=0
\\& 
\text{ and }((Ad(g)W_I).f_1)(g((M(0)\cap H^{(m-n-r)})N(0)))=0
\end{align*}
where $f_0\in C^\infty(G/((M(0)\cap H^{(r)}) N(0)))$ 
and $f_1\in C^\infty(G/((M(0)\cap H^{(m-n-r)}) N(0)))$. 
It follows from Proposition 4-1(1) that 
the element $W_I$ is contained in both 
$$Ann_{U(\mathfrak g_\mathbb C)}
(C^\infty(G/((M\cap H^{(r)}) N))
\text{ and }
Ann_{U(\mathfrak g_\mathbb C)}(C^\infty(G/((M\cap H^{(m-n-r)})N)).$$
This yields the assertion for the element $W_I$. 
The assertion for the element $W_{p+2}$ follows directly 
from the assertion for the element $W_I$. \qed 
\end{Proof}
\par 
Let $m-n\geq 2,n\geq 1$ and $p'\geq 1,2p'\leq m+n$. 
We define 
$$W^{n}_{p'}  
=\begin{cases}  W^{n}_{p'}=\sum\limits_{I\in T_{2p'}(m+n),
I\subset(n+1,\ldots,m)}W_{I}^2 & \text{ if } 2p'\leq m-n \\ 
        0  & \text{ if }  2p'>m-n. 
\end{cases}$$ 
By Fact 3-1(2) we have $W^n_{p'}\in \mathcal Z(\mathfrak m(0))$. 
Here we have $W^n_{1}=\sum\limits_{n+1\leq i<j\leq m}X_{ij}^2$. 
Since $W^n_{p'}\in \mathcal Z(\mathfrak m(0))$, 
it follows from Proposition 2-1 that 
\begin{equation*}(Ad(u_M)W^n_{p'})\in\mathcal Z(\mathfrak m).
\tag{4.3}
\end{equation*}
\par 
\textbf{Proposition 4-3
(The Calculations of the Pfaffians 
of $\mathfrak s\mathfrak o(m,n)$)}
\par 
Assume that $m-n\geq 2$ and $n\geq 1$.  
Suppose that $p\geq 0,p+2>n$ and $2p+4\leq m+n$. 
\par 
$(1)$ Let $I\in T_{2p+4-2n}(m+n)$ such that $I\subset(n+1,\ldots,m)$. 
Then we have 
\begin{align*} & 
(W_{I\cup\{1,\ldots,n\}\cup\{m+1,\ldots,m+n\}})^2
=((\prod\limits_{1\leq l\leq n}
(-(Y_{l,m+n+1-l}-(p+2-l))^2)(W_{I}^2))+
\mathfrak n(0)_{\mathbb C}U(\mathfrak g_{\mathbb C}). 
\end{align*}  
$(2)$ We have 
\begin{align*} & W_{p+2}
\\& 
\in 
(\sum\limits_{1\leq k\leq n}
\sum\limits_{J\in T_{k}(n)}
(\prod\limits_{1\leq l\leq k}
(-((Y_{j_l,m+n+1-j_l}-(p+2-j_l))(Y_{l,m+n+1-j_l}-(m+n-4-p+j_l))))
W^n_{p+2-k}))
\\& 
+W^n_{p+2}
+\mathfrak n(0)_{\mathbb C}\cdot U(\mathfrak g_\mathbb C). 
\end{align*} 
\par 
\begin{Proof}
$(1)$ It follows from the latter assertion in Proposition 3-7(1) 
with (4.1.1) and $N=m+n$ that 
\begin{equation*}
(W_{I\cup\{1,\ldots,n\}\cup\{m+1,\ldots,m+n\}})^2
=((-(Y_{1,m+n}-(p+1))^2))
(W_{I'\setminus\{1,m+n\}}^2))+\mathfrak n(0)_{\mathbb C}
U(\mathfrak g_{\mathbb C}), 
\tag{4.3.1}
\end{equation*}
where we put 
$$I'=I\cup\{1,\ldots,n\}\cup\{m+1,\ldots,m+n\}.$$ 
We prove the assertion by the induction on $n$. 
The assertion for $n=1$ follows directly from (4.3.1). 
\par 
Assume that $n\geq 2$. By the inductive hypothesis we have 
\begin{align*} & W_{I'\setminus\{1,m+n\}}^2
=W_{I\cup\{2,\ldots,n\}\cup\{m+1,\ldots,m+n-1\}})^2
\\& 
=((\prod\limits_{2\leq l\leq n}
(-(Y_{l,m+n+1-l}-(p+2-l))^2)(W_{I}^2)))
+\mathfrak n(0)_{\mathbb C}U(\mathfrak g_{\mathbb C}). 
\tag{4.3.2}
\end{align*} 
We substitute (4.3.2) into (4.3.1). 
Here we have 
$$U(\mathfrak a(0)_{\mathbb C})\mathfrak n(0)_{\mathbb C}
\subset(\mathfrak n(0)_{\mathbb C})
U(\mathfrak a(0)_{\mathbb C}).$$
Then the assertion follows immediately.  
\par 
$(2)$ By Proposition 3-7(3) 
with (4.1.1) and $N=m+n$, 
we have 
\begin{align*} & 
W_{p+2}
\in((-(Y_{1,m+n}-(p+1))(Y_{1,m+n}-(m+n-3-p))
\sum\limits_{1\in I,m+n\in I}W_{I\setminus\{1,m+n\}}^2
\\& 
+\sum\limits_{1\not\in I,m+n\not\in I}W_I^2
+\mathfrak n(0)_{\mathbb C}U(\mathfrak g_{\mathbb C}).
\tag{4.3.3}
\end{align*} 
We prove the assertion by the induction on $n$.  
The assertion for $n=1$ follows directly from (4.3.3). 
Assume that $n\geq 2$. In the setting of (4.3.3), 
we have the following (4.3.4) and (4.3.5) by the inductive hypothesis.
\begin{align*} & \sum\limits_{1\in I,m+n\in I}W_{I\setminus\{1,m+n\}}^2
\in(\sum\limits_{2\leq k\leq n}
\sum\limits_{J\in T_{k}(n),j_1=1}
\\& 
(\prod\limits_{2\leq l\leq k}
(-((Y_{j_l,m+n+1-j_l}-(p+2-j_l)))(Y_{l,m+n+1-j_l}-(m+n-4-p+j_l))))
W^n_{p+2-k})
\\& 
+W^n_{p+1} 
+\mathfrak n(0)_{\mathbb C}\cdot U(\mathfrak g_\mathbb C). 
\tag{4.3.4}
\end{align*} 
\begin{align*} & 
\sum\limits_{1\not\in I,m+n\not\in I}W_I^2
\in(\sum\limits_{1\leq k\leq n}
\sum\limits_{J\in T_{k}(n),j_1\not=1}
\\&
(\prod\limits_{1\leq l\leq k}
(-((Y_{j_l,m+n+1-j_l}-(p+2-j_l))(Y_{l,m+n+1-j_l}-(m+n-4-p+j_l))))
W^n_{p+2-k})
\\& 
+W^n_{p+2}
+\mathfrak n(0)_{\mathbb C}\cdot U(\mathfrak g_\mathbb C). 
\tag{4.3.5}
\end{align*} 
We substitute (4.3.4) and (4.3.5) into (4.3.3). 
We have 
$$U(\mathfrak a(0)_{\mathbb C})\mathfrak n(0)_{\mathbb C}
\subset(\mathfrak n(0)_{\mathbb C})U(\mathfrak a(0)_{\mathbb C}).$$ 
Then the assertion follows immediately. 
\qed 
\end{Proof}
\par 
Suppose $H=H^{(r)}$ for $0\leq r\leq m-n$. 
Let $(\tau,V_\tau)$ be 
a finite dimensional irreducible 
unitary representation of $M$ 
with a nonzero $(M\cap H)$-fixed element in $V_\tau$. 
Let $\hat{M}_H$ be the set of equivalence classes 
of such representations. 
Let $\langle *,*\rangle_\tau$ be a nonzero $M$-invariant 
inner product on $V_\tau$. 
Let $\Vert v\Vert_\tau=(\langle v,v\rangle_\tau)^{1/2}$ 
for any $v\in V_\tau$.
Let $v_H$ be a nonzero $(M\cap H)$-fixed element in $V_\tau$ 
such that $\Vert v_H\Vert_\tau=1$. 
Let $v_\tau$ be a highest weight vector in $V_\tau$ 
in the setting of Proposition 3-9 and Proposition 3-10 
where $N=m-n$. 
By the proof of Theorem 4-1 on p536 in [17], 
we may assume that 
\begin{equation*} 
v_H=\int\nolimits_{M\cap H}\tau(h).v_\tau dh.
\tag{4.4.0}
\end{equation*}
In view of (1.2) and (2.3) it follows from the proof 
of Lemma 2-7 in [26]
that the restriction of the representation $\tau$ 
from $M$ to $M_0$ on $V_\tau$ is irreducible 
with a unit $(M_0\cap H)$-fixed element $v_H\in V_\tau$. 
Then the restriction of the representation $\tau$ 
from $M$ to $M_0^+$ on $V_\tau$ is trivial. 
\par 
Then $\hat{M}_H$ can be identified with the set 
of equivalence classes of all the irreducible unitary representations 
of $M_0$ on finite dimensional vector spaces  
with nonzero $(M_0\cap H)$-fixed elements. 
\par 
\textbf{Proposition 4-4}
\par 
Assume that $\dim G/((M\cap S)N)<\dim G/((M\cap T)N)$. 
For any $\tau\in\hat{M}_T$ and $p'\geq 1,2p'\leq m+n$, 
there exists some $D_{p'}(\tau)\in\mathbb C$ such that 
$$(Ad(u_M)W^n_{p'}).
\langle\tau(m).v_T,v\rangle_\tau=D_{p'}(\tau)\langle\tau(m).v_T,v\rangle_\tau
\text{ for }m\in M,v\in V_\tau$$ 
Moreover we have $D_{\tilde{r_1}+1-n}(\tau)=0$ 
if and only if $\tau\in\hat{M}_S$. 
\par 
\begin{Proof}
\par 
Assume that $\dim G/((M\cap S)N)<\dim G/((M\cap T)N)$. 
It follows from the assumption that 
$\dim (M\cap S)/(M\cap S\cap T)
<\dim (M\cap T)/(M\cap S\cap T)$. 
Since $\hat{M}_S$ or $\hat{M}_T$ can be identified 
with the set of equivalence classes 
of all irreducible unitary representations of $M_0$ 
on a finite dimensional vector space 
with a nonzero $M_0\cap S$ or $M_0\cap T$-fixed element, 
it follows from Proposition 3-9(1),Proposition 3-10 
with $N=m-n$ and Proposition 2-1 with (4.3) that 
there exists some $D_{p'}(\tau)\in\mathbb C$ for $\tau\in\hat{M}_T$ 
such that we have for $m\in M\cap T$ and $v\in V_\tau$, 
\begin{equation*}
(Ad(u_M)W^n_{p'}).\langle \tau(m).v_\tau,v\rangle_\tau
=D_{p'}(\tau)\langle v_\tau, v\rangle_\tau 
\tag{4.4.1}
\end{equation*}
and such that $D_{\tilde{r_1}+1-n}(\tau)=0$ 
if and only if $\tau\in\hat{M}_S$. 
By (4.4.0) and (4.4.1) we have for $m\in M$ and $v\in V_\tau$,
\begin{align*} & 
(Ad(u_M)W^n_{p'}).(\langle\tau(m).v_T,v\rangle_\tau)
\\& 
=(Ad(u_M)W^n_{p'}).
\int\nolimits_{M\cap T}(\langle\tau(m).\tau(t).v_\tau,v\rangle_\tau)dt 
=\int\nolimits_{M\cap T}
(Ad(u_M)W^n_{p'}).(\langle\tau(m\cdot t).v_\tau,v\rangle_\tau)dt.  
\end{align*}  
By (4.4.1) and (4.4.0), 
we have for $m\in M$ and $v\in V_\tau$, 
\begin{align*} & (Ad(u_M)W^n_{p'}).(\langle\tau(m).v_T,v\rangle_\tau)
=\int\nolimits_{M\cap T}D_{p'}(\tau)
\langle\tau(m\cdot t).v_\tau, v \rangle_\tau dt 
\\&
=D_{p'}(\tau)
\langle\int\nolimits_{M\cap T}\tau(m).\tau(t).v_\tau dt, v \rangle_\tau
=D_{p'}(\tau)
\langle\tau(m).\int\nolimits_{M\cap T}
\tau(t).v_\tau dt, v \rangle_\tau 
=D_{p'}(\tau)\langle\tau(m).v_T,v \rangle_\tau.
\end{align*}
Thus we have proved the former assertion. 
The latter assertion follows from Proposition 3-10. \qed 
\end{Proof}
\par 
Let $H=H^{(r)}$ for $0\leq r\leq m-n$. 
Let $F\in C^\infty(G/(M\cap H)N)$. 
For any irreducible unitary representation $\tau$ of $M$ 
on a finite dimensional vector space $V_\tau$ 
with a unit $M\cap H$-fixed element $v_H$, 
we define a function $(F)^\tau_v\in C^\infty(G/N)$ by   
\begin{equation*}
(F)^\tau_v(gN)
=\int\nolimits_{M}
F(g\cdot m(M\cap H)N)
\overline{\langle \tau(m).v_{H},v\rangle_\tau} dm 
\text{ for }g\in G,v\in V_\tau.
\tag{4.5}
\end{equation*}
We call $(F)^\tau_v$ the Fourier transform 
of $F\in C^\infty(G/((M\cap H)N))$ on $M$. 
If $F\in\Gamma(G/((M\cap H)N))$ then 
$(F)^\tau_v\in \Gamma(G/(M\cap H)N)$ 
for $\Gamma=C^\infty,\mathcal S,C^\infty_c,C^\infty_d$ with $d>0$ 
and $C^\infty_\lambda$ 
with $\lambda\in\mathfrak a^*_{\mathbb C}$. 
\par 
\textbf{Proposition 4-5(The Fourier Transform on $M$)}
\par 
Let $H=H^{(r)}$ for $0\leq r\leq m-n$. 
\par 
$(1)$(The injectivity) 
Let $F\in C^\infty(G/((M\cap H)N))$. Suppose 
$$(F)^\tau_v(gN)=0
\text{ for any }g\in G,\tau\in\hat{M}_{H},v\in V_\tau.$$ 
Then we have 
$$F(gN)=0\text{ for }g\in G.$$ 
\par 
$(2)$(Fourier Inversion Formula) 
We have for $g\in G$ and $m\in M$, 
$$F(g mN)=\sum\limits_{\tau\in\hat{M}_H}
(\dim V_\tau)
\sum\limits_{1\leq i\leq\dim V_\tau}
F^\tau_{v^\tau_i}(g N)\langle\tau(m).v_H,v^\tau_i\rangle_\tau,$$ 
where $\{v^\tau_i\bigm\vert 1\leq i\leq\dim V_\tau \}$ 
is an orthonormal basis of $V_\tau$. 
\par 
$(3)$ Let $0\leq k\leq n$. For $\tau\in\hat{M}_H$ 
and $v\in V_\tau$, we have for $g\in G$ and $m\in M$, 
$$(((Ad(g)(Ad(u_M)W^{n}_{\tilde{r_1}+1-k}))).F)^\tau_v(g N)
=D_{\tilde{r_1}+1-k}(\tau)F^\tau_v(gN).$$
\begin{Proof}
$(1)$ In view of Lemma 2-4(3) in [26] and Lemma 2-5(3)(4) in [26] it follows 
from (8) on p391 in [17] and 
Theorem 3.5 on p533 in [17] that 
the functions on $M/(M\cap H)$ 
$$m(M\cap H)\to\langle \tau(m).v_H,v\rangle_\tau$$
where $\tau\in\hat{M}_H,v\in V_\tau$ 
span $L^2(M/(M\cap H))$. 
This yields the assertion. 
\par 
$(2)$ It follows from (8) on p391 in [17] and 
Theorem 3.5 on p533 in [17] that 
the functions on $M/(M\cap H)$ 
$$m(M\cap H)\to(\dim V_\tau)^{1/2}\langle \tau(m).v_H,v^\tau_i\rangle_\tau$$
where $\tau\in\hat{M}_H$ and $1\leq i\leq\dim V_\tau$ 
constitute an orthonormal basis of $L^2(M/(M\cap H))$.  
This yields the assertion of Proposition 4-5(2). 
\par 
$(3)$ By Proposition 4-5(2) and Proposition 4-4 we have 
\begin{align*} 
& ((Ad(g)(Ad(u_M)W^{n}_{\tilde{r_1}+1-k}))).F)(g m N)
\\& 
=\sum\limits_{\tau_0\in\hat{M}_H}
(\dim V_{\tau_0})
\sum\limits_{1\leq i\leq\dim V_{\tau_0}}
(F)^{\tau_0}_{v^{\tau_0}_i}(g N)
(Ad(u_M)W^{n}_{\tilde{r_1}+1-k}).\langle\tau_0(m).v_H,v^\tau_i\rangle_{\tau_0}
\\& 
=\sum\limits_{\tau_0\in\hat{M}_H}
(\dim V_{\tau_0})
\sum\limits_{1\leq i\leq\dim V_{\tau_0}}
(F)^{\tau_0}_{v^{\tau_0}_i}(g N)
D_{\tilde{r_1}+1-k}(\tau)\langle\tau_0(m).v_H,v^{\tau_0}_i\rangle_{\tau_0}. 
\tag{4.5.1}
\end{align*} 
We put $m=e$. Then we have 
\begin{align*} & ((Ad(g)(Ad(u_M)W^{n}_{\tilde{r_1}+1-k}))).F)(g N)
\\& 
=\sum\limits_{\tau_0\in\hat{M}_H}
(\dim V_{\tau_0})
\sum\limits_{1\leq i\leq\dim V_{\tau_0}}
(F)^{\tau_0}_{v^{\tau_0}_i}(g N)
D_{\tilde{r_1}+1-k}(\tau)\langle v_H,v^{\tau_0}_i\rangle_{\tau_0}. 
\tag{4.5.2}
\end{align*} 
For $\tau\in\hat{M}_H$ and $v\in V_\tau$, we have, 
\begin{align*} & (((Ad(g)(Ad(u_M)W^{n}_{\tilde{r_1}+1-k}))).F)^\tau_v(g N)
=\sum\limits_{\tau_0\in\hat{M}_H}
(\dim V_{\tau_0})
\\&
\times\sum\limits_{1\leq i\leq\dim V_{\tau_0}}
(F)^{\tau_0}_{v^{\tau_0}_i}(g N)
D_{\tilde{r_1}+1-k}(\tau_0)
\int\nolimits_{M}
\langle\tau_0(m).v_H,v^{\tau_0}_i\rangle_{\tau_0}
\overline{\langle\tau(m).v_H,v\rangle_\tau}dm
\\& 
=\sum\limits_{1\leq i\leq\dim V_{\tau}}
(F)^{\tau}_{v^\tau_i}(g N)
D_{\tilde{r_1}+1-k}(\tau)
\int\nolimits_{M}
\langle\tau(m).v_H,v^{\tau}_i\rangle_{\tau}
\overline{\langle\tau(m).v_H,v\rangle_\tau}dm,  
\end{align*}
where we use the Schur orthogonality relations in (8) in p391 in [17]. By the Schur orthogonality relations in (8) in p391 in [17], 
we have 
$$\int\nolimits_{M}
\langle\tau(m).v_H,v^{\tau}_i\rangle_{\tau}
\overline{\langle\tau(m).v_H,v\rangle_\tau}dm
=\langle v_H,v_H\rangle_{\tau}
\overline{\langle v^{\tau}_i,v\rangle_{\tau}}
=\overline{\overline{\langle v,v^{\tau}_i\rangle_{\tau}}}.$$ 
Then we have 
\begin{align*} & (((Ad(g)(Ad(u_M)W^{n}_{\tilde{r_1}+1-k}))).F)^\tau_v(g N)
=\sum\limits_{1\leq i\leq\dim V_{\tau}}
(F)^\tau_{v_i}(g N)
D_{\tilde{r_1}+1-k}(\tau)
\overline{\overline{\langle v,v^{\tau}_i\rangle_{\tau}}} 
\\& 
=D_{\tilde{r_1}+1-k}(\tau)\sum\limits_{1\leq i\leq\dim V_{\tau}}
\int\nolimits_{M}
F(gmN)\overline{\langle\tau(m).v_H, v^\tau_i\rangle_{\tau}}
dm\overline{\overline{\langle v,v^{\tau}_i\rangle_{\tau}}} 
\\& 
=D_{\tilde{r_1}+1-k}(\tau)
\int\nolimits_{M}
F(gmN)\overline{\langle\tau(m).v_H, 
\sum\limits_{1\leq i\leq\dim V_{\tau}}
\langle v,v^\tau_i\rangle_\tau v^\tau_i\rangle_{\tau}} dm 
\\& 
=D_{\tilde{r_1}+1-k}(\tau)\int\nolimits_{M}
F(gmN)\overline{\langle\tau(m).v_H, v\rangle_{\tau}} dm 
=D_{\tilde{r_1}+1-k}(\tau)F^\tau_v(gN).
\end{align*}
Thus we have proved the assertion. 
\qed 
\end{Proof}
\par 
Assume that $m-n\geq 2$ and $n\geq 1$.  
Suppose that $p\geq 0,p+2>n$ and $2p+4\leq m+n$. 
\par 
In the setting of Proposition 3-3(1) and Proposition 2-1, 
we define for $\lambda\in\mathfrak a^*_{\mathbb C}$,  
\begin{equation*}
P(\lambda)
=\prod\limits_{1\leq l\leq n}
(-((Ad(u_M^{-1})\lambda)(Y_{l,m+n+1-l})+(p+2-l))^2).
\tag{4.6.1}
\end{equation*}
\par 
\textbf{Proposition 4-6}
\par 
Assume that $n\geq 1; 0\leq r_1,r_2 \leq m-n ; r_1\not= r_2$. 
Assume that $\dim G/((M\cap S)N))<\dim G/((M\cap T)N))$.
Let $\tilde r_i=\min(m-r_i,r_i+n)$ for $i=1,2$. 
Suppose $f_T\in C^\infty_\lambda(G/((M\cap T)N))$.
Then we have for $\tau\in\hat{M}_T$, 
\begin{align*} & 
((Ad(g)Ad(u_M)
\sum\limits_{I\in T_{2(\tilde{r_1}+1-n)}(m+n),I\subset(n+1,\ldots,m)}
W_{I\cup\{1,\ldots,n\}\cup\{m+1,\ldots,m+n\}}^2)
.(f_T))^{\tau}_v(g N).
\\& 
=P(\lambda)D_{\tilde{r_1}+1-n}(\tau)(f_T)^{\tau}_v(g N) 
\text{ for }g\in G,v\in V_\tau. 
\end{align*} 
\par 
\begin{Proof}
We have for $g\in G$ and $m\in M$,
\begin{align*} & f_T(g m (M\cap T)N)
=f_T(g u_M((u_M)^{-1} m u_M)(u_M)^{-1} (M\cap T)N)
=f_T(g_0 m_0 (u_M)^{-1}(M\cap T)N),
\end{align*} 
where we put $g_0=g u_M\in G$ and 
$m_0=(u_M)^{-1} m u_M\in M(0)$.
By Proposition 2-1 we have  
$$Ad(u_M^{-1})\lambda
\in(\mathfrak a(0))^*_{\mathbb C} 
\text{ and }
f_T(*  u_M^{-1}(M\cap T)N)
\in C^\infty_{Ad(u_M^{-1})\lambda}(G/((M(0)\cap T)N(0))).$$ 
Let $I\in T_{2(\tilde{r_1}+1-n)}(m+n)$ such that 
$I\subset\{n+1,\ldots,\dots,m+n\}$. 
By Proposition 4-3(1),(2.2) and (4.6.1),  
we have for $g\in G$ and $m\in M$,
\begin{align*} & 
((Ad(g)Ad(u_M)
\sum\limits_{I\in T_{2(\tilde{r_1}+1-n)}(m+n),I\subset(n+1,\ldots,m) } 
W_{I\cup\{1,\dots,n\}\cup\{n,\ldots,m+n\}}^2)).
f_T)(g m(M\cap T)N)
\\& 
=((Ad(g_0)
\sum\limits_{I\in T_{2(\tilde{r_1}+1-n)}(m+n),I\subset(n+1,\ldots,m)} 
W_{I\cup\{1,\dots,n\}\cup\{n,\ldots,m+n\}}^2).
f_T)(g_0 m_0 u_M^{-1}(M\cap T)N)
\\& 
=P(\lambda)((Ad(g_0)\sum_{I\in T_{2(\tilde{r_1}+1-n)}(m+n),I\subset(n+1,\ldots,m)}
W_I^2).f_T)(g_0 m_0 u_M^{-1}(M\cap T)N)
\\& 
=P(\lambda)
((Ad(g)Ad(u_M)
\sum\limits_{I\in T_{2(\tilde{r_1}+1-n)}(m+n),I\subset(n+1,\ldots,m)} 
W_I^2).f_T)(gm(M\cap T)N) 
\\& 
=P(\lambda)
((Ad(g)Ad(u_M)W^n_{\tilde{r_1}+1-n}).f_T)(gm(M\cap T)N).
\end{align*} 
By Proposition 4-5(3) we have for $\tau\in\hat{M}_T$ and $v\in V_\tau$,  
\begin{align*} & 
(((Ad(g)Ad(u_M)
\sum\limits_{I\in T_{2(\tilde{r_1}+1-n)}(m+n),I\subset(n+1,\ldots,m) } 
W_{I\cup\{1,\dots,n\}\cup\{n,\ldots,m+n\}}^2).
f_T))^\tau_{v}(g N))
\\& 
=P(\lambda)D_{\tilde{r_1}+1-n}(\tau)(f_T)^\tau_v(g N).
\end{align*} 
Thus we have proved the assertion. \qed 
\end{Proof}
In the setting of  Proposition 4-3(2) and Proposition 2-1, 
we define a polynomial $P_k(\lambda)$ 
of $\lambda\in\mathfrak a^*_{\mathbb C}$ 
for $0\leq k\leq n$ by  
\begin{align*} & P_k(\lambda) 
=\sum\limits_{J\in T_{k}(n)}
(\prod\limits_{1\leq l\leq k}
(-((Ad(u_M^{-1})\lambda)Y_{j_l,m+n+1-j_l}+(p+2-j_l))
\\& 
\times(Ad(u_M^{-1})\lambda)(Y_{l,m+n+1-j_l})+(m+n-4-p+j_l)))
\text{ for }1\leq k\leq n 
\tag{4.7.0.1}
\end{align*}
and 
\begin{equation*} 
P_0(\lambda)=1.
\tag{4.7.0.2 }
\end{equation*}
Then we have 
\begin{align*} & 
P_n(\lambda)
=\prod\limits_{1\leq l\leq n}
(-((Ad(u_M^{-1})\lambda)(Y_{l,m+n+1-l})+(p+2-l))
\\& 
\times((Ad(u_M^{-1})\lambda)(Y_{l,m+n+1-l})+(m+n-4-p+l)))).
\tag{4.7.0.3 }
\end{align*}
Here the degree of $P_n$ is $2n$ 
and the degree of $P_k$ 
for $0\leq k<n$ is less than $2n$. 
\par 
Let $H=H^{(r)}$ for $0\leq r\leq m-n$. 
Let $D$ be a subset in $\mathfrak a^*_{\mathbb C}$ 
such that 
\begin{equation*}
D=H_0+\bigoplus\limits_{1\leq i\leq n}c_i\mathbb R(Ad(u_M)e_i),
\tag{4.7.0.4}
\end{equation*}
where $H_0\in\mathfrak a^*_{\mathbb C}$ and 
$c_i\in\mathbb C\setminus\{0\}$ for $1\leq i\leq n$. 
\par 
\textbf{Proposition 4-7(The eigenvalue of the central elements 
of $\mathfrak s\mathfrak o(m,n)$)}
\par 
Assume that $n\geq 1; 0\leq r_1,r_2 \leq m-n ; r_1\not= r_2$. 
Assume that $\dim G/((M\cap S)N))<\dim G/((M\cap T)N))$.
Let $\tilde r_i=\min(m-r_i,r_i+n)$ for $i=1,2$. 
\par 
$(1)$ Let $\lambda\in\mathfrak a^*_{\mathbb C}$ 
and $F\in C^\infty_\lambda(G/((M\cap T)N))$. 
Then we have for $\tau\in\hat{M}_T$ and $v\in V_\tau$,  
$$(W_{\tilde r_1+1}.F)^\tau_v(gN)
=Q^\tau_{r_1}(\lambda)(F^\tau_v)(gN)
\text{ for }g\in G,m\in M,$$ 
where we put 
$$Q^\tau_{r_1}(\lambda)
=\sum\limits_{0\leq k\leq n}
P_k(\lambda) D_{\tilde{r_1}+1-k}(\tau).$$ 
\par 
$(2)$ If $\tau\in\hat{M}_T\setminus\hat{M}_S$, then we have 
$Q^\tau_{r_1}\not=0$ on a open dense subset of $D$. 
\par 
\begin{Proof}
$(1)$ We have for $g\in G$ and $m\in M$,
\begin{align*} & F(g m (M\cap T)N)
=F(g u_M((u_M)^{-1} m u_M)(u_M)^{-1} (M\cap T)N)
=F(g_0 m_0 (u_M)^{-1} (M\cap T)N),
\end{align*} 
where we put $g_0=g u_M\in G$ and $m_0=(u_M)^{-1} m u_M\in M(0)$.
By Proposition 2-1 we have for $g\in G$ and $m\in M$, 
$$Ad(u_M^{-1})\lambda\in(\mathfrak a(0))^*_{\mathbb C} 
\text{ and }
F(* u_M^{-1}(M\cap T)N)
\in C^\infty_{Ad(u_M^{-1})\lambda}(G/(M(0)\cap T)N(0)).$$ 
By Proposition 4-1(2) and Proposition 4-3(2) with (2.1) and (2.2) and by (4.7.0.1),(4.7.0.2) and (4.7.0.3),  we have for $g\in G$ and $m\in M$, 
\begin{align*} & (W_{\tilde{r_1}+1}.(F))(g m (M\cap T)N)
=(Ad(g)Ad(u_M)W_{\tilde{r_1}+1}.(F))(g m (M\cap T)N)
\\& 
=(Ad(g_0)W_{\tilde{r_1}+1}.(F))(g_0 m_0 u_M^{-1} (M\cap T)N)
\\& 
=((Ad(g_0)(\sum\limits_{0\leq k\leq n}
P_k(\lambda)W^{n}_{\tilde{r_1}+1-k})).F)(g_0 m_0 u_M^{-1} (M\cap T)N)
\\& 
=((\sum\limits_{0\leq k\leq n}
P_k(\lambda)Ad(g)Ad(u_M)W^{n}_{\tilde{r_1}+1-k})).F)(g m (M\cap T)N).
\tag{4.7.1}
\end{align*} 
By Proposition 4-5(3), 
we have for $\tau\in\hat{M}_T$ and $v\in V_\tau$,  
\begin{align*} & (W_{\tilde{r_1}+1}.F)^\tau_v(g N)
=\sum\limits_{0\leq k\leq n}
P_k(\lambda)((Ad(g)(Ad(u_M)W^{n}_{\tilde{r_1}+1-k}))).F)^\tau_v(g N)
\\& 
=\sum\limits_{0\leq k\leq n}
P_k(\lambda)D_{\tilde{r_1}+1-k}(\tau)F^\tau_v(g N)
=Q^\tau_{r_1}(\lambda)F^\tau_v(gN).
\tag{4.7.2}
\end{align*} 
Thus we have proved the assertion. 
\par 
$(2)$ By (4.7.0.1),(4.7.0.2),(4.7.0.3) and the definition of $Q^\tau_{r_1}$ 
we have for $\tau\in\hat{M}_T$,  
\begin{align*} & 
(\prod\limits_{1\leq i\leq n}\frac{\partial^2}{\partial t_i^2})
\bigm\vert_{t_i=0(1\leq i\leq n)}
Q^\tau_{r_1}(\lambda+\sum\limits_{1\leq i\leq n}c_i\cdot t_i(Ad(u_M)e_i))
\\& 
=(-2)^n(\prod\limits_{1\leq i\leq n}c_i^2)D_{\tilde{r_1}+1-n}(\tau)
\text{ for }\lambda\in D.
\tag{4.7.3}
\end{align*} 
We assume that $\tau\in\hat{M}_T\setminus\hat{M}_S$. 
By Proposition 4-4 we have $D_{\tilde{r_1}+1-n}(\tau)\not=0$. 
Then the function of $\lambda\in D$ in (4.7.3)  
is a nonzero constant function on $D$. 
Suppose there exists an open subset $V$ in $D$ such that 
$Q^\tau_{r_1}(\lambda)=0\text{ for }\lambda\in V$.  
Then the function of $\lambda\in D$ in (4.7.3) is zero on $V$. 
This yields a contradiction. 
Then $Q^\tau_{r_1}\not=0$ on a open dense subset of $D$. 
Thus we have proved the assertion.  \qed 
\end{Proof}
\par 
\textbf{Remark}
The following diagram is commutative. 
\[\begin{CD} 
  C^\infty_\lambda(G/((M\cap T)N)) @>  W_{\tilde r_1+1}  >> C^\infty_\lambda(G/((M\cap T)N)) \\
@V\  Fourier\text{ }Transform\text{ }on\text{ }M  VV   @VV\ Fourier\text{ }Transform\text{ }on\text{ }M  V\\
    C^\infty_\lambda((G/N))  @>>  Multiplication\text{ }of\text{ }  of\text{ }Q^\tau_{r_1}(\lambda)   >  C^\infty_\lambda((G/N))
\end{CD}\] 
\par 
\section{ The Injectivity of The Radon Transform on the Horocycle Spaces.} 
\par 
In this section we investigate the injectivity properties 
of the Radon transform $R$ and its dual $R^*$
by use of the Fourier transform on $M$. 
\par 
Assume that $n\geq 1; 0\leq r_1,r_2 \leq m-n ; r_1\not= r_2$. 
We put 
$$\tilde{\mathfrak h}_{\mathfrak m}
=\bigoplus\limits_{n+1\leq i\leq [(m+n)/2]}
\mathbb C Ad(u_M)X_{i,m+n+1-i}.$$ 
Then $\tilde{\mathfrak h}_{\mathfrak m}$ 
is a maximal abelian subalgebra of $\mathfrak m_{\mathbb C}$. 
Let $d(\tau)\in(\tilde{\mathfrak h}_{\mathfrak m})^*$ 
be the highest weight of the finite dimensional irreducible 
unitary representation $(\tau, V_\tau)$ of $M$. 
We put 
$$d_j(\tau)
=d(\tau)(\tilde Ad(u_M)X_{n+j, m+1-j})\text{ for }1\leq j\leq [(m+n)/2]-n.$$
By (2.10) we have 
$$d_j(\tau)\in 2\mathbb Z\text{ for }1\leq j\leq  [(m+n)/2]-n.$$
Let $\vert\tau\vert$ be the norm of the highest weight $d(\tau)$ 
with respect to the Killing form. 
\par 
\textbf{Proposition 5-1}
\par 
Assume that $n\geq 1; 0\leq r_1,r_2 \leq m-n ; r_1\not= r_2$. 
\par 
$(1)$ The homogeneous spaces $G/((M\cap S)N)$ and $G/((M\cap T)N)$ 
are homogeneous spaces in duality 
in the definition in Ch II,\S 1 in $[20]$. 
\par 
$(2)$ The homogeneous spaces $K/(M\cap S)$ and $K/(M\cap T)$ 
are homogeneous spaces in duality 
in the definition in ChII \S 1 in $[20]$. 
\par 
\begin{Proof}
$(1)$ By Proposition 2-1 we have only to prove that 
$G/((M(0)\cap S)N(0))$ and 
\par 
$G/((M(0)\cap T)N(0))$ are homogeneous spaces in duality. 
By ChII,\S 1,Lemma1-3 in [20] we have only to prove that 
$(M(0)\cap S)/(M(0)\cap S\cap T)$ 
and $(M(0)\cap T)/(M(0)\cap S\cap T)$ are homogeneous spaces in duality. 
This follows from Lemma 2-4(3) in [26], (1.2) 
and the fact that $M(0)_0/(M(0)_0\cap S)$ and $M(0)_0/(M(0)_0\cap T)$ 
are homogeneous spaces in duality in the setting of p135 in [26]. 
Thus we have proved the assertion. 
\par 
$(2)$ The assertion follows in a similar way to Proposition 5-1(1). 
\qed 
\end{Proof}
\par 
We put for $\tau\in\hat{M}_T\cap\hat{M}_S$,  
\begin{equation*}
C(\tau)=\langle v_T,v_S\rangle_\tau.
\tag{5.2.0.1}
\end{equation*} 
We put 
\begin{equation*}\gamma(r_1,r_2;\tau)
=\prod\limits_{1\leq j\leq \min(r_1,m-n-r_1)}
\frac{((r_1-j+1)/2)_{d_j(\tau)/2}(m-n-r_2-j+1)/2)_{d_j(\tau)/2}}
{((r_2-j+1)/2)_{d_j(\tau)/2}(m-n-r_1-j+1)/2)_{d_j(\tau)/2}}
\tag{5.2.0.2}
\end{equation*}
where we denote the Pochammer symbol by $(c)_{l}=c(c+1)\cdots(c+l-1)\text{ for }c\in\mathbb C,l\in\mathbb Z_{\geq 0}$ (See Theorem 6.1. in [40]). 
\par 
\textbf{Proposition 5-2(The Properties of the Radon Transform $R_M$)}
\par 
Assume that $n\geq 1; 0\leq r_1,r_2 \leq m-n ; r_1\not= r_2$. 
and that $\dim G/((M\cap S)N)\leq\dim G/((M\cap T)N)$. 
\par 
$(1)$ We have $\hat{M}_S\subset\hat{M}_T$. 
Moreover we have 
\begin{equation*}
\int\nolimits_{(M\cap S)/(M\cap S\cap T)}
(\tau(s)v_T) ds_{(M\cap S\cap T)}
=\begin{cases}  C(\tau)v_S  &  \text{ for }\tau\in\hat{M}_S \\ 
         0  & \text{ for }\tau\in\hat{M}_T\setminus\hat{M}_S, 
\end{cases}
\tag{5.2.1.1}
\end{equation*}
and for $\tau\in\hat{M}_S$ 
\begin{equation*}\int\nolimits_{M\cap T/(M\cap S\cap T)}
(\tau(t)v_S) ds_{(M\cap S\cap T)}
=\overline{C(\tau)}v_T \text{ for }\tau\in\hat{M}_S.
\tag{5.2.1.2}
\end{equation*}
\par 
$(2)$ We have $\tau\in\hat{M}_S$, 
\begin{equation*}\vert C(\tau)\vert=
\begin{cases} (\gamma(r_1,r_2;\tau))^{1/2} & \text{ if }r_1<r_2 \\ 
(\gamma(r_1,r_2;\tau))^{-1/2} & \text{ if }r_1>r_2. 
\end{cases}\tag{5.2.2.1}
\end{equation*}
Moreover, we have $C(\tau)\not=0$ for $\tau\in\hat{M}_S$ and  
\begin{equation*}
\sup\limits_{\tau\in\hat{M}_S}
\vert C(\tau)\vert^{\pm1}
(1+\vert\tau\vert)^{-l_\pm}<\infty
\text{ for some }l_{\pm}\in\mathbb Z^+_{\geq 0}.
\tag{5.2.2.2}
\end{equation*}
\par 
$(3)$ Let $\tau\in\hat{M}_S$. 
For any unit $M\cap S$-fixed element $v_S\in V_\tau$, 
we can choose 
some unit $M\cap T$-fixed element $v_T\in V_\tau$ 
so that $C(\tau)\in \mathbb R$. 
\par 
\begin{Proof}
$(1)$ It follows from the assumption that 
$\dim (M\cap S)/(M\cap S\cap T)
\leq\dim (M\cap T)/(M\cap S\cap T)$. 
In view of Proposition 2-1 
we have only to prove the assertion when $M=M(0)$ and $N=N(0)$. 
Since $\hat{M}_S$ or $\hat{M}_T$ can be identified 
with the set of equivalence classes 
of all irreducible unitary representations of $M_0$ 
with a nonzero $M_0\cap S$ or $M_0\cap T$-fixed element, 
it follows from Proposition 3-10 with $N=m-n$ 
that $\hat{M}_S\subset\hat{M}_T$.   
It follows from Lemma 2-7 in [26] that 
there exists some nonzero $C_1\in\mathbb C$ 
for $\tau\in\hat{M}_S$ such that 
$$\int\nolimits_{(M\cap S)/(M\cap S\cap T)}
(\tau(s)v_T) ds_{(M\cap S\cap T)}
=C_1\cdot v_S.$$ 
Moreover, we have for $\tau\in\hat{M}_T\setminus\hat{M}_S$, 
$$\int\nolimits_{(M\cap S)/(M\cap S\cap T)}
(\tau(s)v_T) ds_{(M\cap S\cap T)}=0.$$ 
This yields the latter equation of (5.2.1.1). 
For any fixed $\tau\in\hat{M}_S$ we have for any $s\in M\cap S$, 
$$C(\tau)=\langle v_T,v_S\rangle_\tau=\langle v_T,\tau(s^{-1}).v_S\rangle_\tau
=\langle \tau(s).v_T,v_s\rangle_\tau.$$ 
By the normalization of the invariant measure, 
we have for $\tau\in\hat{M}_S$, 
\begin{align*} &  C(\tau)=\langle v_T,v_S\rangle_\tau 
=\int\nolimits_{(M\cap S)/(M\cap S\cap T)}
\langle \tau(s).v_T,v_S\rangle_\tau ds_{(M\cap S\cap T)} 
\\& 
=\langle\int\nolimits_{(M\cap S)/(M\cap S\cap T)}
(\tau(s)v_T) ds_{(M\cap S\cap T)}, v_S \rangle_\tau 
=\langle C_1\cdot v_S, v_S \rangle_\tau 
=C_1\cdot \langle v_S,v_S \rangle_\tau=C_1. 
\end{align*} 
This yields the former equation of (5.2.1.1). 
\par 
It follows in a similar way to Lemma 2-7 in [26] that 
there exists some nonzero $C_2\in\mathbb C$ 
for $\tau\in\hat{M}_S$ such that 
$$\int\nolimits_{(M\cap T)/(M\cap S\cap T)}
(\tau(t)v_S) dt_{(M\cap S\cap T)}=C_2\cdot v_T.$$ 
Then in a similar way we have 
$$C_2=\langle v_S,v_T \rangle_\tau=\overline{\langle v_T,v_S \rangle_\tau}=\overline{C(\tau)}.$$  
This yields (5.2.1.2). 
\par 
$(2)$ By Lemma 2-7(2.7.4) in [26] we have 
\begin{align*} & \int\nolimits_{(M\cap S)/(M\cap T\cap S)}
\int\nolimits_{(M\cap T)/(M\cap T\cap S)}
(\tau(s).(\tau(t).v_S)) dt_{(M\cap T\cap S)}ds_{(M\cap T\cap S)}
=C_1\cdot C_2\cdot v_S=\vert C(\tau)\vert^2 v_S.
\tag{5.2.2.3}
\end{align*}  
Suppose $r_1<r_2$. 
By the assumption that $\dim(G/((M\cap S)N))<\dim(G/((M\cap T)N))$, 
we have $r_1<r_2\leq m-n-r_1$. 
It follows from (5.2.2.3), \S 5 in [14], Theorem 6.1 in [40] 
and (5.2.0.2) that we have for $\tau\in\hat{M}_S$, 
$$\vert C(\tau)\vert^2=\gamma(r_1,r_2;\tau).$$ 
This yields the assertion (5.2.2.1) when $r_1<r_2$. 
\par  
Suppose $r_2<r_1$. 
By the assumption that $\dim(G/((M\cap S)N))\leq \dim(G/((M\cap T)N))$, 
we have $m-n-r_1\leq r_2<r_1$ and $m-n-r_1\leq m-n-r_2<r_1$. 
In view of Lemma 2-7 in [26] and Lemma 6.2 in [29] it follows that 
the Radon transform $R_M$ for double fibrations 
$$M/(M\cap H^{(r_1)})\leftarrow M/(M\cap H^{(r_1)}\cap H^{(r_2)})
 \rightarrow M/(M\cap H^{(r_2)})$$  
can be identified with the Radon transform $R_M$ for double fibrations 
$$M/(M\cap H^{(m-n-r_1)})\leftarrow 
M/(M\cap H^{(m-n-r_1)}\cap H^{(m-n-r_2)})
 \rightarrow M/(M\cap H^{(m-n-r_2)}).$$ 
Then we replace $r_i$ by $m-n-r_i$ for $i=1,2$ 
in the assertion (5.2.2.1) for $r_1<r_2$ 
Then we have for $\tau\in\hat{M}_S$, 
$$\vert C(\tau)\vert^2=\gamma(r_1,r_2;\tau)^{-1}.$$ 
This yields the assertion (5.2.2.1) when $r_2<r_1$. 
\par 
The assertion (5.2.2.2) follows directly from (5.2.2.1). 
Thus we have proved all the assertions. 
\par 
$(3)$ Let $v_T$ be any unit $M\cap T$-fixed element in $V\tau$.
Then we put 
$$(v_T)_0
=\frac{\vert\langle v_T,v_S\rangle_\tau\vert}
{\langle v_T,v_S\rangle_\tau}v_T.$$ 
Then we have 
$$\langle (v_T)_0,v_S\rangle_\tau
=\vert\langle v_T,v_S\rangle_\tau\vert\in\mathbb R.$$ 
We replace $v_T$ by $(v_T)_0$. Then we put 
$$C(\tau)=\langle (v_T)_0, v_S \rangle_\tau\in\mathbb R.$$ 
This yields the assertion. \qed 
\end{Proof}
\par 
\textbf{Proposition 5-3(The Radon transform $R_M$)}
\par 
Assume that $n\geq 1; 0\leq r_1,r_2 \leq m-n ; r_1\not= r_2$.
\par 
$(1)$ We have for $f_S\in C^\infty(G/((M\cap S)N))$, 
$$(R(f_S))(gm((M\cap T)N))=R_M(f_S)(g*N)(m(M\cap T))
\text{ for }g\in G,m\in M.$$ 
\par 
In addition, we have for $f_T\in C^\infty(G/((M\cap T)N))$, 
$$(R^* f_T)(gm((M\cap S)N))=R_M^*(f_T)(g*N)(m(M\cap S))
\text{ for }g\in G,m\in M.$$ 
\par 
$(2)$ The Radon transform $R$ maps $\Gamma(G/((M\cap S)N))$ 
into $\Gamma(G/((M\cap T)N))$ 
for $\Gamma=C^\infty,\mathcal S,C^\infty_c,C^\infty_d$ where $d>0$ 
or $C^\infty_\lambda$ where $\lambda\in\mathfrak a_{\mathbb C}$. 
\begin{Proof}
$(1)$ We have 
$$(M\cap T)N/((M\cap S)N\cap(M\cap T)N)=(M\cap T) (M\cap S\cap T)N.$$ 
where $(M\cap T)(M\cap S\cap T)N$ is diffeomorphic to  
$(M\cap T)/(M\cap S\cap T)$. 
Then $((M\cap T)N)$-invariant measure 
on $(M\cap T)N/((M\cap S)N\cap(M\cap T)N)$ 
is also $(M\cap T)$-invariant measure, 
since $(M\cap T)$ is a subgroup of $(M\cap T)N$. 
By the normalization of the invariant measure in \S 2, we have 
$$\int\nolimits_{M\cap T}dt
=\int\nolimits_{M\cap T\cap S}ds
=\int\nolimits_{(M\cap T)/(M\cap T\cap S)}dt_{M\cap T\cap S}=1.$$
It follows from (2.0) and the uniqueness of the invariant measure 
on $((M\cap T)N)/((M\cap S)\cap(M\cap T))$ 
and $(M\cap T)/(M\cap T\cap S)$ that we have 
for $f\in C^\infty(G/((M\cap S)N))$,
\begin{align*} & (R f)(gm((M\cap T)N))
=\int\nolimits_{((M\cap T)N)/((M\cap S)N\cap(M\cap T)N)}
f(gmt((M\cap S)N))dt_{((M\cap T)N)\cap((M\cap S)N)}
\\& 
=\int\nolimits_{(M\cap T)/(M\cap T\cap S)}
f(gmt((M\cap S)N))dt_{(M\cap T\cap S)}
=R_M.(f)(g*n)(m(M\cap T)).
\end{align*} 
This yields the former assertion. 
The latter assertion for $R^*$ follows in a similar way. 
\par 
(2) The assertion directly follows from Proposition 5-3(1). 
\qed 
\end{Proof}
\par 
\textbf{Proposition 5-4(The Projection Slice Theorem of $R$)}
\par 
Assume that $n\geq 1; 0\leq r_1,r_2 \leq m-n ; r_1\not= r_2$. 
and that $\dim G/((M\cap S)N)\leq\dim G/((M\cap T)N)$. 
\par 
$(1)$ We have for $f_S\in C^\infty(G/((M\cap S)N))$ and $g\in G$, 
\begin{equation*}(R f_S)^\tau_v(gN)
=\begin{cases} 
\overline{C(\tau)}(f_S)^\tau_v(gN) & \text{ for }\tau\in\hat{M}_S,v\in V_\tau \\ 
         0  & \text{ for }\tau\in\hat{M}_T\setminus\hat{M}_S,v\in V_\tau. 
\end{cases}
\end{equation*}
\par 
$(2)$ We have for $f_T\in C^\infty(G/((M\cap T)N))$ and $g\in G$, 
$$(R^* f_T)^\tau_v(gN)
=C(\tau)(f_T)^\tau_v(gN)\text{ for }\tau\in\hat{M}_S,v\in V_\tau.$$ 
\par 
\begin{Proof}
$(1)$ By Ch.II,\S 2,Proposition 2-2 in [20] 
and Proposition 5-2(1),  
we have for $\tau\in\hat{M}_T$, 
\begin{align*} & 
(R f)^\tau_v(g N)
=\int\nolimits_{M}(R f)(gm((M\cap T)N))\overline{\langle\tau(m).v_T,v\rangle_\tau} dm
\\& 
=\int\nolimits_{M}(R_{M}.f)(gm(M\cap T))\overline{\langle\tau(m).v_T,v\rangle_\tau} dm
\\& 
=\int\nolimits_{M}f(gm(M\cap S))
\int\nolimits_{(M\cap S)/(M\cap S\cap T)}\overline{\langle\tau(m).\tau(s).v_T,v\rangle_\tau} 
ds_{(M\cap S\cap T)}
\\& 
=\int\nolimits_{M}f(gm(M\cap S))
\overline{\langle \tau(m).
\int\nolimits_{(M\cap S)/(M\cap S\cap T)}(\tau(s)v_T)
ds_{(M\cap S\cap T)},v\rangle_\tau} ds_{(M\cap S\cap T)}.  
\end{align*} 
By Proposition 5-2(1)(5.2.1.1) we have for $g\in G$, 
$$(R f)^\tau_v(gN)
=\begin{cases}  \overline{C(\tau)}
f^\tau_v(gN)  &  \text{ for }\tau\in\hat{M}_S,v\in V_\tau \\ 
         0  & \text{ for }\tau\in\hat{M}_T\setminus\hat{M}_S,v\in V_\tau. 
\end{cases}$$ 
By Ch.II,\S 2,Proposition 2-2 in [20] and Proposition 5-3(1),  
we have for $\tau\in\hat{M}_S$, 
\begin{align*} & 
(R^* f_T)^\tau_v(g N)
=\int\nolimits_{M}(R^* f_T)(gm((M\cap S)N))
\overline{\langle\tau(m).v_S,v\rangle_\tau} dm
\\& 
=\int\nolimits_{M}(R_{M}^*.f_T)(gm((M\cap S)N))
\overline{\langle\tau(m).v_S,v\rangle_\tau} dm
\\& 
=\int\nolimits_{M}f(gm((M\cap T)N))
\int\nolimits_{(M\cap T)/(M\cap S\cap T)}
\overline{\langle\tau(m).\tau(t).v_S,v\rangle_\tau} 
dt_{(M\cap S\cap T)}
\\& 
=\int\nolimits_{M}f_T(gm((M\cap S)N))
\overline{\langle \tau(m).
\int\nolimits_{(M\cap T)/(M\cap S\cap T)}(\tau(t)v_S)
dt_{(M\cap S\cap T)},v\rangle_\tau} ds_{(M\cap S\cap T)}.  
\end{align*} 
By Proposition 5-2(1)(5.2.1.2) we have for $g\in G$, 
$$(R^* f_T)^\tau_v(gN)
=C(\tau)(f_T)^\tau_v(gN)\text{ for }\tau\in\hat{M}_S,v\in V_\tau.$$ 
Thus we have proved all the assertions. 
\qed
\end{Proof}
\par 
\textbf{Remark}
The following diagram is commutative. 
\[\begin{CD} 
  C^\infty_\lambda(G/((M\cap S)N)) @>  R  >> C^\infty_\lambda(G/((M\cap T)N)) \\
@V\  Fourier\text{ }Transform\text{ }on\text{ }M  VV   @VV\ Fourier\text{ }Transform\text{ }on\text{ }M  V\\
    C^\infty_\lambda((G/N))  @>>  Multiplication\text{ }of\text{ }  of\text{ }\overline{C(\tau)}   >  C^\infty_\lambda((G/N))
\end{CD}\] 
\par 
\textbf{Proposition 5-5(The Injectivity of the Radon transform $R$) }
\par 
Assume that $n\geq 1; 0\leq r_1,r_2 \leq m-n ; r_1\not= r_2$.
\par 
$(1)$ Assume that $\dim G/((M\cap S)N<\dim G/((M\cap T)N)$. 
Then the Radon transform $R$ maps $\Gamma(G/((M\cap S)N))$ 
injectively into $\Gamma(G/((M\cap T)N))$ 
for $\Gamma=C^\infty,\mathcal S,C^\infty_c,C^\infty_d$ with $d>0$ 
or $C^\infty_\lambda$ with $\lambda\in\mathfrak a_{\mathbb C}$. 
The dual Radon transform $R^*$ is not injective 
on $\Gamma(G/((M\cap T)N))$ 
for $\Gamma=C^\infty,\mathcal S,C^\infty_c, C^\infty_d$ 
with $d>0$ or $C^\infty_\lambda$ 
with $\lambda\in\mathfrak a_{\mathbb C}$.
\par 
$(2)$ Assume that $\dim G/((M\cap S)N))=\dim G/((M\cap T)N)$. 
Then the Radon transform $R$ maps $\Gamma(G/((M\cap S)N))$ 
injectively into $\Gamma(G/((M\cap T)N))$ 
for $\Gamma=C^\infty,\mathcal S,C^\infty_c$,
$C^\infty_d$ with $d>0$ 
or $C^\infty_\lambda$ with $\lambda\in\mathfrak a_{\mathbb C}$. 
The dual Radon transform $R^*$ 
maps $\Gamma(G/((M\cap T)N))$ injectively 
into $\Gamma(G/((M\cap S)N))$ 
for $\Gamma=C^\infty,\mathcal S,C^\infty_c,C^\infty_d$ 
with $d>0$ or $C^\infty_\lambda$ 
with $\lambda\in\mathfrak a_{\mathbb C}$. 
\par 
$(3)$ Assume that $\dim G/((M\cap S)N))>\dim G/((M\cap T)N)$. 
Then the Radon transform $R$ is not injective 
on $\Gamma(G/((M\cap S)N))$ 
for $\Gamma=C^\infty,\mathcal S,C^\infty_c, C^\infty_d$ with $d>0$ 
or $C^\infty_\lambda$ with $\lambda\in\mathfrak a_{\mathbb C}$.  
The dual Radon transform $R^*$ maps $\Gamma(G/((M\cap T)N))$ 
injectively into $\Gamma(G/((M\cap S)N))$ 
for $\Gamma=C^\infty,\mathcal S,C^\infty_c,C^\infty_d$ with $d>0$ 
or $C^\infty_\lambda$ with $\lambda\in\mathfrak a_{\mathbb C}$. 
\par 
\begin{Proof} 
$(1)$ The assertion for $R$ follows directly from Proposition 5-4(1), since we have $C(\tau)\not=0$ for $\tau\in\hat{M}_S$ 
by Proposition 5-2(1).    
In view of Lemma 2-4(3) in [26] 
and Lemma 2-5(3)(4) in [26] it follows from 
Theorem 3.5 on p533 in [17] that 
there exists some $f_T\in \Gamma(G/((M\cap T)N))$ 
such that we have for any $\tau\in\hat{M}_S$, 
$$(f_T)^\tau_v(gN)=0\text{ for any }g\in G,v\in V_\tau$$ 
and such that we have for some $\tau\in\hat{M}_T\setminus\hat{M}_S$, 
$v\in V_\tau$ and $g\in G$,  
$$(f_T)^\tau_v(gN)\not=0.$$ 
Here we have $f_T\not\equiv 0$. 
By Proposition 5-4(2) we have for any $\tau\in\hat{M}_S,v\in V_\tau$, 
$$(R^*(f_T))^\tau_v(gN)=0\text{ for }g\in G.$$ 
It follows from Proposition 4-5(1) that 
$$R^*(f_T)(g(M\cap S)N)=0\text{ for any }g\in G.$$ 
This yields the assertion. 
\par 
$(2)$ The assertion for $R$ follows directly from Proposition 5-4(1), 
since we have $C(\tau)\not=0$ for $\tau\in\hat{M}_S$ 
by Proposition 5-2(1).   
Suppose that $R^*(f_T)\equiv 0$ for $f_T\in \Gamma(G/((M\cap T)N))$. 
Since we have $C(\tau)\not=0$ 
for $\tau\in\hat{M}_S$ by Proposition 5-2(1), 
it follows from Proposition 5-4(1) that 
$$(f_T)^\tau_v(gN)=0\text{ for }\tau\in\hat{M}_S,v\in V_\tau,g\in G.$$ 
Since it follows from the assumption that 
$\dim (M\cap S)/(M\cap S\cap T)=\dim (M\cap T)/((M\cap S\cap T))$, 
we have $\hat{M}_S=\hat{M}_T$ 
by Lemma 2-7 in [26] and Proposition 3-10 with $N=m-n$. 
Then we have $f_T\equiv 0$ by Proposition 4-5(1). 
This yields the assertion. 
\par 
$(3)$ The assertion follows from Proposition 5-3(1), 
since it is equivalent to the assertion of Proposition 5-3(1). \qed 
\end{Proof}
\par 
\textbf{Proposition 5-6(The definition of $R^+$ and $(R^*)^+$)}
\par 
Assume that $n\geq 1; 0\leq r_1,r_2 \leq m-n ; r_1\not= r_2$ 
and $\dim G/((M\cap S)N))\leq\dim G/((M\cap T)N))$. 
\par 
$(1)$ For $f_T\in C^\infty(G/((M\cap T)N))$, we put for $g\in G$,  
\begin{equation*}
(R^+(f_T))(g N)
=\sum\limits_{\tau\in\hat{M}_S} 
\overline{C(\tau)}^{-1}(\dim V_\tau)
(f_T)^\tau_{v_S}(gN).
\tag{5.6.1}
\end{equation*} 
\par 
Then the right-hand side converges  
and $R^+(f_T)\in C^\infty(G/((M\cap S)N))$. 
Let $\Gamma=\mathcal S,C^\infty_c,C^\infty_d$ with $d>0$ 
or $C^\infty_\lambda$ with $\lambda\in\mathfrak a_{\mathbb C}$.  
If $f_T\in \Gamma(G/((M\cap T)N))$ then 
$(R^+(f_T))\in\Gamma(G/((M\cap S)N))$. 
\par 
$(2)$ For $f_S\in C^\infty(G/((M\cap S)N))$, we put  
\begin{equation*}((R^*)^+(f_S))(g N)
=\sum\limits_{\tau\in\hat{M}_S} 
C(\tau)^{-1}(\dim V_\tau)
(f_S)^\tau_{v_T}(gN)
\text{ for }g\in G.
\tag{5.6.2}
\end{equation*}
\par 
Then the right-hand side converges  
and $(R^*)^+(f_S)\in C^\infty(G/((M\cap T)N))$. 
Let $\Gamma=\mathcal S,C^\infty_c,C^\infty_d$ with $d>0$ 
or $C^\infty_\lambda$ 
with $\lambda\in\mathfrak a_{\mathbb C}$.   
If $f_S\in \Gamma(G/((M\cap S)N))$ then 
$((R^*)^+(f_S))\in\Gamma(G/((M\cap T)N))$. 
\par 
\begin{Proof}
$(1)$ By Theorem 4 in [43] we have for any fixed $g\in G$, 
$$\sup\limits_{\tau\in\hat{M}_T,v\in V_\tau,\Vert v\Vert_\tau=1}
\vert (f_T)^\tau_v)(gN)\vert 
(1+\vert\tau\vert^l)\vert
<\infty\text{ for any }l\in\mathcal Z^+_{\geq 0}.$$ 
Since $\Vert v_S\Vert_\tau=1$ for $\tau\in\hat{M}_S$, 
the convergence follows from (5.2.2.2) 
and the Weyl dimension formula for $\dim V_\tau$
(See Theorem 1.8 on p502 in [17]). 
Moreover we have $R^+(f_T)\in C^\infty(G/N)$ 
by the dominated convergence theorem. 
By the definition we have 
$$(f_T)^\tau_v(g m N)=(f_T)^\tau_{\tau(m).v}(g N)
\text{ for }g\in G,m\in M,\tau\in\hat{M}_T,v\in V_\tau.$$ 
Then we have 
$$R^+(f_T)(g m N)=R^+(f_T)(g N)\text{ for }g\in G,m\in M\cap S.$$ 
Therefore we have $R^+(f_T)\in C^\infty(G/((M\cap S)N))$. 
\par 
Moreover the latter assertion 
for $\Gamma=C^\infty_c,C^\infty_d$ where $d>0$ 
or $C^\infty_\lambda$ where 
$\lambda\in\mathfrak a_{\mathbb C}$ 
follows from (4.5). 
\par 
Let $\Gamma=\mathcal S$. 
Suppose $f_T\in\mathcal S(G/((M\cap T)N))$. 
By Theorem 4 in [43] we have 
\begin{align*} & 
\sup\limits_{u\in K,a\in A,\tau\in\hat{M}_T,v\in V_\tau,\Vert v\Vert_\tau=1}
\vert (X.f_T)^\tau_v)(uaN)\vert 
(1+\vert\tau\vert^{l_1})(1+\vert \log a\vert^{l_2})
<\infty
\text{ for any }X\in U(\mathfrak g),l_1,l_2\in\mathbb Z^+_{\geq 0}
\end{align*} 
This yields $R^+(f_T)\in\mathcal S(G/((M\cap S)N))$. 
Thus we have proved the assertions. 
\par 
$(2)$ The assertion follows in a similar way to Proposition 5-6(1). 
\qed 
\end{Proof}
\par 
\textbf{Proposition 5-7(The Inversion Formula of the Radon transform $R$)}
\par 
Assume that $n\geq 1; 0\leq r_1,r_2 \leq m-n ; r_1\not= r_2$ 
and $\dim G/((M\cap S)N))\leq\dim G/((M\cap T)N))$. 
\par 
$(1)$ Let $f_T\in C^\infty(G/((M\cap T)N))$. 
We have for $\tau\in\hat{M}_S$, 
$$(R(R^+(f_T)))^\tau_v(gN)
=(f_T)^\tau_v(gN)\text{ for }g\in G,v\in V_\tau.$$ 
\par 
$(2)$ Let $f_S\in C^\infty(G/((M\cap S)N)$. 
We have for $\tau\in\hat{M}_S$, 
$$(R^*((R^*)^+(f_S)))^\tau_v(gN)
=(f_S)^\tau_v(gN)\text{ for }g\in G,v\in V_\tau.$$ 
\par 
$(3)$(Inversion Formula) 
We have for $f_S\in C^\infty(G/((M\cap S)N))$, 
$$(R^+(R(f_S)))(g((M\cap S)N))=f_S(g((M\cap S)N))
\text{ for }g\in G.$$ 
\par 
$(4)$(Support Theorem) 
Suppose $f\in C^\infty(G/((M\cap S)N))$. 
Then we have $f\in C^\infty_c(G/((M\cap S)N))$ 
if and only if $R(f)\in  C^\infty_c(G/((M\cap T)N))$. 
Moreover, we have $f\in C^\infty_d(G/((M\cap S)N))$ 
if and only if $R(f)\in  C^\infty_d(G/((M\cap T)N))$ 
for $d>0$.
\par 
\begin{Proof} 
$(1)$ Let $\tau\in\hat{M}_S$ and $v\in V_\tau$, 
By Proposition 5-4(1) we have  
$$(R(R^+(f_T)))^\tau_{v}(gN)
=\overline{C(\tau)}\int\nolimits_{M}
(R^+(f_T))(gm((M\cap S)N))
\overline{\langle\tau(m).v_S,v\rangle_\tau} dm. $$ 
By Proposition 5-6(1) we have 
\begin{align*} & 
(R(R^+(f_T)))^\tau_{v}(gN)
\\& 
=\overline{C(\tau)}\int\nolimits_{M}
(\sum\limits_{\tau_1\in\hat{M}_S}
(\dim V_{\tau_1})\overline{C(\tau_1)}^{-1}
(f_T)^{\tau_1}_{v_S}(gm((M\cap T)N)))
\overline{\langle\tau(m).v_S,v\rangle_\tau} dm
\\& 
=\overline{C(\tau)}\int\nolimits_{M}
\sum\limits_{\tau_1\in\hat{M}_S}
(\dim V_{\tau_1})\overline{C(\tau_1)}^{-1}
\\&
\times\int\nolimits_{M}f_T(gmm_1((M\cap T)N))
\overline{\langle\tau_1(m_1).v_T,v_S\rangle_{\tau_1}} dm_1
\overline{\langle\tau(m).v_S,v\rangle_\tau} dm.
\\& 
=\overline{C(\tau)}\int\nolimits_{M}
\sum\limits_{\tau_1\in\hat{M}_S}
(\dim V_{\tau_1})\overline{C(\tau_1)}^{-1}
\\&
\times\int\nolimits_{M}f_T(gm_0((M\cap T)N))
\overline{\langle\tau_1(m_0).v_T,\tau_1(m).v_S\rangle_{\tau_1}}
dm_0
\overline{\langle\tau(m).v_S,v\rangle_\tau} dm, 
\end{align*}
where we put $m_0=m\cdot m_1$ in the last equality. 
By the Schur orthogonality relations in (8) in p391 in [17], 
we have 
\begin{align*} & 
(R(R^+(f_T)))^\tau_{v}(gN)
\\&
=\overline{C(\tau)}\int\nolimits_{M}
(\dim V_{\tau})\overline{C(\tau)}^{-1}
\int\nolimits_{M}f_T(gm_0((M\cap T)N))
\overline{\langle \tau(m_0).v_T,\tau(m).v_S\rangle_{\tau}}
dm_0
\overline{\langle\tau(m).v_S,v\rangle_\tau} dm 
\\&
=\int\nolimits_{M}f_T(gm_0((M\cap T)N))
(\dim V_\tau)\int\nolimits_{M} 
\langle\tau(m).v_S,\tau(m_0).v_T \rangle_\tau 
\overline{\langle\tau(m).v_S,v\rangle_\tau } dm dm_0.
\end{align*}
By the Schur orthogonality relations in (8) in p391 in [17], 
we have 
\begin{align*} & 
(R(R^+(f_T)))^\tau_{v}(gN)
\\& 
=\overline{C(\tau)}\int\nolimits_{M}
(\dim V_{\tau})\overline{C(\tau)}^{-1}
\int\nolimits_{M}f_T(gm_0((M\cap T)N))
\overline{\langle \tau(m_0).v_T,\tau(m).v_S\rangle_{\tau}}
dm_0
\overline{\langle\tau(m).v_S,v\rangle_\tau} dm 
\\&
=\int\nolimits_{M}f_T(gm_0((M\cap T)N))
(\dim V_\tau)\int\nolimits_{M} 
\langle\tau(m).v_S,\tau(m_0).v_T \rangle_\tau 
\overline{\langle\tau(m).v_S,v\rangle\tau } dm dm_0.
\end{align*}
By the Schur orthogonality relations in (8) in p391 in [17], 
we have 
\begin{align*} & 
(\dim V_\tau)\int\nolimits_{M} 
\langle\tau(m).v_S,\tau(m_0).v_T\rangle_{\tau}
\overline{\langle\tau(m).v_S,v\rangle_\tau} dm
=\langle v_S,v_S\rangle_\tau
\overline{\langle \tau(m_0).v_T,v\rangle_\tau}
=\overline{\langle \tau(m_0).v_T,v\rangle_\tau}.
\end{align*} 
Therefore we have 
\begin{align*} & (R(R^+(f_T)))^\tau_{v}(gN)
=\int\nolimits_{M}f_T(gm_0((M\cap T)N))
\overline{\langle\tau(m_0).v_T,v\rangle_\tau} 
dm_0=(f_T)^\tau_v(gN)
\end{align*} 
Thus we have proved the assertion. 
\par 
$(2)$ The assertion follows in a similar way to Proposition 5-7(1).
\par 
$(3)$ By Proposition 5-7(1) we have for $\tau\in\hat{M}_S$, 
$$(R(R^+(R(f_S))))^\tau_v(gN)=(R(f_S))^\tau_v(gN)
\text{ for }g\in G,v\in V_\tau.$$
By Proposition 5-4(1) we have for $\tau\in\hat{M}_T\setminus\hat{M}_S$, 
$$(R(R^+(R(f_S))))^\tau_v(gN)=(R(f_S))^\tau_v(gN)=0
\text{ for }g\in G,v\in V_\tau.$$ 
Hence we have for $\tau\in\hat{M}_T$, 
$$(R(R^+(R(f_S))))^\tau_v(gN)=(R(f_S))^\tau_v(gN)
\text{ for }g\in G,v\in V_\tau.$$
By Proposition 4-5(1) we have $R(R^+(R(f_S)))\equiv R(f_S)$. 
Therefore it follows from the injectivity of $R$ in Proposition 5-5(1) 
that $R^+(R(f_S))=f_S$. Thus we have proved the assertion. 
\par 
$(4)$ The assertion follows from Proposition 5-3(2), 
Proposition 5-6(1) and Proposition 5-7(3) 
for $\Gamma=C^\infty_c$ or $C^\infty_d$ with $d>0$. 
\qed 
\end{Proof}
\par 
\textbf{Remark}
The assertions of 
Proposition 5-1,5-2,5-3,5-4 and 5-5 hold 
in the following setting(See (1.3) in [26]). 
\begin{align*} (G,S,T) &
=(SU(m,n;\mathbb F),
S(U(m-r_1;\mathbb F)\times U(r_1,n;\mathbb F)), 
S(U(m-r_2;\mathbb F)\times U(r_2,n;\mathbb F)))
\\& 
\text{ where } 
n\geq 1; 0\leq r_1,r_2 \leq m-n ; r_1\not= r_2,
\mathbb F=\mathbb C,\mathbb H.
\end{align*} 

\section{ The Range of the Radon Transform on the Horocycle Spaces.} 
\par 
In this section we investigate the range questions 
of the Radon transform $R$ by use of the Fourier transform on $M$ 
and the Fourier transform on $A$. 
\par 
\textbf{Proposition 6-1(The Surjectivity of the Radon Transform $R$)}
\par 
Assume that $n\geq 1; 0\leq r_1,r_2 \leq m-n ; r_1\not= r_2$. 
Let $\Gamma=C^\infty,\mathcal S,C^\infty_c,C^\infty_d$ with $d>0$ 
or $C^\infty_\lambda$ 
with $\lambda\in\mathfrak a_{\mathbb C}$. 
Assume $\dim G/((M\cap S)N)\geq\dim G/((M\cap T)N)$. 
Then the Radon transform $R$ maps $\Gamma(G/((M\cap S)N))$ 
surjectively onto $\Gamma(G/((M\cap T)N))$. 
\par 
\begin{Proof}
We have only to prove that $R^*$ maps $\Gamma(G/((M\cap T)N))$ surjectively onto $\Gamma(G/((M\cap S)N))$ 
under the condition that $\dim G/((M\cap S)N))\leq\dim G/((M\cap T)N)$. 
This follows from Proposition 5-5(1),Proposition 5-6(2) 
and Proposition 5-7(2). This yields the assertion. \qed 
\end{Proof}
Let $L$ be a closed subgroup of $M$. 
We define the Fourier transform 
$\tilde{f}(g (LN);\lambda)\in C^\infty(G/(LN))$ 
of $f\in \mathcal S(G/(LN))$  on $A$ by 
\begin{equation*}
\tilde{f}(g (LN);\lambda)
=\int\nolimits_{A}f(g a LN)e^{-i\lambda(\log a)}da
\text{ for }\lambda\in\mathfrak a^*.
\tag{6.2}
\end{equation*}
We define 
$$\mathcal L((G/(LN))\times\mathfrak a^*)
=\{F\in C^\infty(G/(LN)\times\mathfrak a^*)\bigm\vert 
F\text{ satisfies the following (L1) and (L2)}\}.$$ 
\begin{equation*}
F(*;\lambda)\in C^\infty_{i\lambda}(G/(LN))
\text{ for }\lambda\in\mathfrak a^*.
\tag{L1}
\end{equation*}
\begin{equation*}
\sup\limits_{u\in K,\lambda\in\mathfrak a^*} 
\vert (X.F)(u (LN);\lambda)\vert
(1+\vert \lambda\vert^l)<\infty\text{ for any }
X\in U(\mathfrak g),l\in\mathbb Z^+_{\geq 0}.
\tag{L2}
\end{equation*}
Here we denote by $\vert *\vert$ the norm on 
$\mathfrak a^*$ induced by the Killing form. 
\par 
\textbf{Proposition 6-2(The Fourier Transform on $A$)}
\par 
Assume that $m\geq n\geq 1$. 
\par 
$(1)$ The Fourier transform $f\mapsto\tilde{f}$ 
maps $\mathcal S(G/(LN))$ bijectively 
onto $\mathcal L((G/(LN))\times\mathfrak a^*)$. 
\par 
$(2)$ Assume that $n\geq 1; 0\leq r_1,r_2 \leq m-n ; r_1\not= r_2$. 
Then we have for $f\in\mathcal S(G/((M\cap S)N))$, 
$$(\tilde{R f})(g((M\cap T)N);\lambda)=R(\tilde{f})
(g((M\cap T)N);\lambda)
\text{ for }g\in G,\lambda\in\mathfrak a^*.$$ 
\par 
\begin{Proof}
$(1)$ The assertion follows from the Schwartz Theorem 
in view of the Iwasawa decomposition $G=KAN$. 
\par 
$(2)$ The assertion follows 
from the change of the order of the integration. 
\qed 
\end{Proof}
\par 
\textbf{Remark}
The following diagram is commutative. 
\[\begin{CD} 
  \mathcal S(G/((M\cap S)N))  @>  R  >> \mathcal S(G/((M\cap T)N))      \\
@V\  Fourier\text{ }Transform\text{ }on\text{ }A  VV   @VV\ Fourier\text{ }Transform\text{ }on\text{ }A  V\\
    \mathcal L((G/((M\cap S)N))\times\mathfrak a^*)    @>>  R  >   
    \mathcal L((G/((M\cap T)N))\times\mathfrak a^*) 
\end{CD}\] 
\par 
By (4.6.1) we define 
\begin{equation*}
(\mathfrak a^*_{\mathbb C})_0
=\{\lambda\in\mathfrak a_{\mathbb C}\bigm\vert 
P(\lambda)\not=0\}.
\tag{6.3.0.1}
\end{equation*}
\par 
By (4.7.0.3) we define 
\begin{equation*}
(\mathfrak a^*_{\mathbb C})_1
=\{\lambda\in\mathfrak a_{\mathbb C}\bigm\vert 
P_n(\lambda)\not=0\}.
\tag{6.3.0.2}
\end{equation*}
Then we have 
$(\mathfrak a^*_{\mathbb C})_1\subset(\mathfrak a^*_{\mathbb C})_0$. 
\par 
\textbf{Proposition 6-3}
\par 
Assume that $n\geq 1; 0\leq r_1,r_2 \leq m-n ; r_1\not= r_2$. 
Assume that $\dim G/((M\cap S)N))<\dim G/((M\cap T)N))$.
Let $\tilde r_i=\min(m-r_i,r_i+n)$ for $i=1,2$. 
\par  
$(1)$ Let $\lambda\in(\mathfrak a^*_{\mathbb C})_0$ 
and $f_T\in C^\infty_\lambda(G/((M\cap T)N))$ 
Suppose that $W_I.f_T\equiv 0$ 
for $I\in T_{2(\tilde{r_1}+1)}(m+n)$.     
Then there exists some $f_S\in C^\infty_\lambda(G/((M\cap S)N))$ 
such that $R f_S=f_T$. 
\par 
$(2)$ Assume that $\tilde{r_2}=\tilde{r_1}+1$. 
Let $\lambda\in(\mathfrak a^*_{\mathbb C})_1$ 
and $f_T\in C^\infty_\lambda(G/((M\cap T)N))$ 
Suppose that $W_{\tilde{r_1}+1}.f_T\equiv 0$. 
Then there exists some $f_S\in C^\infty_\lambda(G/((M\cap S)N))$ 
such that $R f_S=f_T$. 
\begin{Proof}
(1) By the assumption and Proposition 4-1(1) 
we have for $J\in T_{2(\tilde{r_1}+1)}(m+n)$, 
\begin{align*} 
& ((Ad(g)Ad(u_M)(W_J)^2).f_T)(gm(M\cap T)N))
\\& 
=(((Ad(g u_M)W_J)^2).f_T)(gm(M\cap T)N))=0
\text{ for }g\in G,m\in M. 
\end{align*} 
Then we have for $g\in G$ and $m\in M$,  
\begin{align*} & 
((Ad(g)Ad(u_M)
\sum\limits_{I\in T_{2(\tilde{r_1}+1-n)}(m+n),I\subset(n+1,\ldots,m)}
(W_{I\cup\{1,\ldots,n\}\cup\{m+1,\ldots,m+n\}})^2).
(f_T))(g m(M\cap T)N)=0. 
\end{align*}

Then we have for $\tau\in\hat{M}_T$,$v\in V_\tau$, and $g\in G$, 
\begin{align*} & 
((Ad(g)Ad(u_M)
\sum\limits_{I\in T_{2(\tilde{r_1}+1-n)}(m+n),I\subset(n+1,\ldots,m)}
(W_{I\cup\{1,\ldots,n\}\cup\{m+1,\ldots,m+n\}})^2).
(f_T))^{\tau}_v(g N)=0. 
\end{align*}

By Proposition 4-6, we have for $\tau\in\hat{M}_T$, 
$$P(\lambda)D_{\tilde{r_1}+1-n}(\tau)((f_T)^{\tau}_v)(g N)=0
\text{ for }g\in G,v\in V_\tau.$$ 
Since $P(\lambda)\not=0$ 
for $\lambda\in(\mathfrak a^*_\mathbb C)_0$ by (6.3.0.1), 
we have for $\tau\in\hat{M}_T$, 
$$D_{\tilde{r_1}+1-n}(\tau)((f_T)^{\tau}_v)(g N)
=0\text{ for }g\in G,v\in V_\tau.$$ 
Since $D_{\tilde{r_1}+1-n}(\tau)\not=0$ for $\tau\in\hat{M}_T\setminus\hat{M}_S$ 
by Proposition 4-4, we have for $\tau\in\hat{M}_T\setminus\hat{M}_S$, 
\begin{equation*}
((f_T)^{\tau}_v)(g N))=0\text{ for }g\in G,v\in V_\tau.
\tag{6.3.1.1}
\end{equation*}
By Proposition 5-6(1) we define $f_S\in C^\infty_\lambda(G/((M\cap S)N))$ 
by $f_S\equiv R^+(f_T)$. 
By Proposition 5-7(1) we have for $\tau\in\hat{M}_S$, 
\begin{equation*}(R(f_S))^\tau_v(gN)
=(f_T)^\tau_v(gN)
\text{ for any }g\in G,v\in V_\tau.
\tag{6.3.1.2}
\end{equation*}
By Proposition 5-4(1) we have for $\tau\in\hat{M}_T\setminus\hat{M}_S$, 
\begin{equation*}
(R(f_S))^\tau_v(g N)=0 
\text{ for any }g\in G,v\in V_\tau.
\tag{ 6.3.1.3}
\end{equation*}
Then the assertion follows from (6.3.1.1),(6.3.1.2) and (6.3.1.3) 
by Proposition 4-5(1). 
\par  
$(2)$ By the assumption we have for $\tau\in\hat{M}_T$, 
$$((W_{\tilde r_1+1}.(f_T))^\tau_v)(g N)\text{ for }g\in G,v\in V_\tau.$$
It follows from $\tilde{r_2}=\tilde{r_1}+1$ and Proposition 4-4 
that $D_{\tilde{r_1}+1-k}(\tau)=0$ 
for $0\leq k<n$ and $\tau\in\hat{M}_T$. 
By Proposition 4-7(1) we have for $\tau\in\hat{M}_T$, 
$$((W_{\tilde r_1+1}.(f_T))^\tau_v)(g N)
=P_n(\lambda)D_{\tilde{r_1}+1-n}(\tau)((f_T)^{\tau}_v)(g N)
=0\text{ for }g\in G,v\in V_\tau.$$ 
Since $P_n(\lambda)\not=0$ 
for $\lambda\in(\mathfrak a_\mathbb C^*)_1$ by (6.3.0.2), 
we have for $\tau\in\hat{M}_T$,
$$D_{\tilde{r_1}+1-n}(\tau)((f_T)^{\tau}_v)(g N)
=0\text{ for }g\in G,v\in V_\tau.$$ 
Then the assertion follows in a similar way to Proposition 6-3(1). 
\qed 
\end{Proof}
\par 
Let $H=H^{(r)}$ for $0\leq r\leq m-n$. 
Let $U_D$ be an open subset of $D$ in (4.7.0.4). 
We define 
\begin{align*} & 
\mathcal D(G/(M\cap H)N\times U_D)
\\& 
=\{F\in C^\infty(G/((M\cap N)N)\times U_D)
\bigm\vert F\in C^\infty_{\lambda}(G/((M\cap H)N))
\text{ for }\lambda\in U_D\}.
\tag{6.4.0}
\end{align*} 
\par
\textbf{Proposition 6-4}
\par 
Assume that $n\geq 1; 0\leq r_1,r_2 \leq m-n ; r_1\not= r_2$. 
Assume that $\dim G/((M\cap S)N))<\dim G/((M\cap T)N))$.
Let $\tilde r_i=\min(m-r_i,r_i+n)$ for $i=1,2$. 
\par 
$(1)$ Let $F_T\in\mathcal D(G/((M\cap T)N)\times U_D)$. 
Suppose 
$$(W_{\tilde r_1+1}.F_T)(g((M\cap T)N);\lambda)=0 
\text{ for any }g\in G,\lambda\in U_D.$$ 
\par 
Then there exists some 
$F_S\in\mathcal D(G/((M\cap T)N)\times U_D)$ 
such that 
\begin{equation*}
(R(F_S))(g((M\cap T)N);\lambda)=F_T(g((M\cap T)N);\lambda)
\text{ for any }g\in G,\lambda\in U_D.
\tag{6.4.1} 
\end{equation*}
\par 
$(2)$ Let $F_T\in\mathcal L(G/((M\cap T)N)\times\mathfrak a^*)$. 
Suppose 
$$(W_{\tilde r_1+1}.F_T)(g((M\cap T)N);i\lambda)=0 
\text{ for any }g\in G,\lambda\in\mathfrak a^*.$$ 
Then there exists some 
$F_S\in\mathcal L(G/((M\cap T)N)\times\mathfrak a^*)$ 
such that 
\begin{equation*} 
(R(F_S))(g((M\cap T)N);i\lambda)=F_T(g((M\cap T)N);i\lambda)
\text{ for any }g\in G,\lambda\in\mathfrak a^*.
\tag{6.4.2}
\end{equation*}
\par 
\begin{Proof}
$(1)$ By the assumption we have 
for $\tau\in\hat{M}_T$ and $v\in V_\tau$.   
$$((W_{\tilde r_1+1}.(F_T))^\tau_v)(g N;\lambda)=0
\text{ for any }g\in G,\lambda\in D.$$ 
By Proposition 4-7(1) we have 
for $\tau\in\hat{M}_T$ and $v\in V_\tau$, 
$$((W_{\tilde r_1+1}.(F_T))^\tau_v)(g N;\lambda)
=Q^\tau_{r_1}(\lambda)((F_T)^{\tau}_v)(g N;\lambda)
=0\text{ for }g\in G,\lambda\in D.$$ 
Assume that $\tau\in\hat{M}_T\setminus\hat{M}_S$. 
Then by Proposition 4-7(2) 
we have $Q^\tau_{r_1}(\lambda)\not=0$ for almost 
all $\lambda\in U_D$. 
Then we have for $\tau\in\hat{M}_T\setminus\hat{M}_S$, 
\begin{equation*}(F_T)^\tau_v(g N;\lambda)=0 
\text{ for any }g\in G,\lambda\in U_D.
\tag{6.4.1.1}
\end{equation*}
By Proposition 5-6(1) 
we define $F_S\in D(G/((M\cap S)N)\times U_D)$ by 
$$F_S(g((M\cap S)N);\lambda)
=(R^+(F_T))(g((M\cap S)N);\lambda)
\text{ for any }g\in G,\lambda\in U_D.$$ 
By Proposition 5-7(1) we have for $\tau\in\hat{M}_S$, 
\begin{equation*}
(R(F_S))^\tau_v(gN;\lambda)
=(F_T)^\tau_v(gN;\lambda)
\text{ for any }g\in G,v\in V_\tau,\lambda\in U_D.
\tag{6.4.1.2} 
\end{equation*}
By Proposition 5-4(1) we have for $\tau\in\hat{M}_T\setminus\hat{M}_S$, 
\begin{equation*}
(R(F_S))^\tau_v(g N;\lambda)=0 
\text{ for any }g\in G,v\in V_\tau,\lambda\in U_D.
\tag{6.4.1.3}
\end{equation*}
The assertion (6.4.1) follows 
from (6.3.4.1),(6.3.4.2) and (6.3.4.3) by Proposition 4-5(1). 
\par 
\par 
$(2)$ It follows from Proposition 6-4(1) 
for $D=U_D=i\mathfrak a^*$  
that there exists some 
$F_S\in \mathcal D(G/((M\cap S)N)\times i\mathfrak a^*)$ 
such that (6.4.2) holds. 
Hence $F_S$ satisfies the condition (L1) for $L=M\cap S$.  
In view of Theorem 4 in [43] 
it follows from the condition (L2) for $F=F_T$ that 

\begin{align*} & \sup\limits_{u\in K,\lambda\in\mathfrak a^*}
\vert ((X.F_T)^\tau_v)(u N;i\lambda)\vert
(1+\vert\lambda\vert^{l_1})(1+\vert\tau\vert^{l_2})
<\infty
\\& 
\text{ for any }X\in U(\mathfrak g),l_1,l_2
\in\mathbb Z^+_{\geq},\tau\in\hat{M}_S,v\in V_\tau.
\tag{6.4.2.1}
\end{align*} 
Then the condition (L2) for $F=F_S$ follows 
from (6.4.2.1),(6.4.2), Proposition 5-7(3) and Proposition 5-6(1). 
Hence we have $F_S\in\mathcal L(G/(M\cap S)N\times\mathfrak a^*)$. 
Thus we have proved the assertion. 
\qed 
\end{Proof}
\par 
\textbf{Remark}
Proposition 6-4(1) in the setting of $n=1$ and $r_1=0$ 
follows from the proof 
of Proposition 2-3 in $[23]$. 
\par 
\textbf{Theorem 6-5(The Range Theorem of the Radon Transform $R$)}
\par 
Assume that $n\geq 1; 0\leq r_1,r_2 \leq m-n ; r_1\not= r_2$ 
and $\dim G/((M\cap S)N)<\dim G/((M\cap T)N)$. 
Let $\tilde r_i=\min(m-r_i,r_i+n)$ for $i=1,2$. 
Then the following $(i)$,$(ii)$ and $(iii)$ hold. 
\par 
$(i)$ The following $(1)$ and $(2)$ hold 
for $\Gamma=\mathcal S,C^\infty_c$ or $C^\infty_d$ with $d>0$.  
\par 
$(ii)$ The following $(1)$ holds  
for $\Gamma=C^\infty_\lambda$ 
with $\lambda\in(\mathfrak a^*_{\mathbb C})_0$. 
\par 
$(iii)$  Assume $\tilde{r_2}=\tilde{r_1}+1$.
Then the following $(1)$ and $(2)$ hold 
for $\Gamma=C^\infty_\lambda$ 
with $\lambda\in(\mathfrak a^*_{\mathbb C})_1$. 
\begin{equation*}  R(\Gamma(G/((M\cap S)N))))
=\{f\in \Gamma(G/((M\cap T)N)))
\bigm\vert W_{I}.f\equiv 0\text{ for }I\in T_{2(\tilde r_1+1)}(m+n)\}.
\tag{1}
\end{equation*}
\begin{equation*}
R(\Gamma(G/((M\cap S)N))))
=\{f\in \Gamma(G/((M\cap T)N))\bigm\vert 
W_{\\tilde r_1+1}.f\equiv 0\}.
\tag{2}
\end{equation*}
\par 
\begin{Proof}
In the assertions of (i),(ii) and (iii), 
the inclusions of the left-hand side 
in the right-hand side in (1) and (2) 
follow from Proposition 4-2. 
Thus we have only to prove the inclusions of the right-hand side 
in the left-hand side in (1) and (2). 
\par 
Then the assertion (ii) follows directly 
from Proposition 6-3(1). 
\par 
Next we prove the assertions (i) and (iii). 
In order to prove the assertions (i) and (iii), 
we have only to prove the inclusion 
of the right-hand side in the left-hand side in (2), 
since the right-hand side in (1) is included 
in the right-hand side in (2). 
\par 
Then the assertion (iii) follows directly 
from Proposition 6-3(2). 
\par 
Next we prove the assertion (i). 
We assume that $\Gamma=\mathcal S$. 
Suppose $f\in \mathcal S(G/((M\cap T)N))$ 
satisfies $W_{\tilde r_1+1}.f\equiv 0$. 
Then we have $W_{\tilde r_1+1}.\tilde{f}\equiv 0$.  
It also follows from Proposition 6-2(1) that 
$\tilde{f}\in\mathcal L(G/((M\cap T)N)\times\mathfrak a^*)$.  
Then it follows from Proposition 6-4(2) for $F_T=\tilde{f}$ 
that there exists some 
$F_S\in\mathcal L(G/((M\cap S)N)\times\mathfrak a^*)$ 
which satisfies (6.4.2). 
It follows from Proposition 6-2(1) 
that there exists some 
$f_S\in \mathcal S(G/((M\cap S)N))$ such that 
$$\tilde{f_S}(g((M\cap S)N);i\lambda)
=F_S(g((M\cap S)N);i\lambda)
\text{ for any }g\in G,\lambda\in\mathfrak a^*.$$ 
Then it follows from Proposition 6-2(2) 
and (6.4.2) for $F_S=\tilde{f_S}$ that 
$$\tilde{R(f_S)}(g((M\cap S)N);i\lambda)=\tilde{f}(g((M\cap S)N);i\lambda)
\text{ for }g\in G,\lambda\in\mathfrak a^*.$$ 
Then we have $R(f_S)=f$ by Proposition 6-2(1). 
This yields the assertion (i) for $\Gamma=\mathcal S$. 
\par 
Next we assume $\Gamma=C^\infty_d$ with $d>0$. 
Suppose $f\in C^\infty_d(G/((M\cap T)N))$ satisfies 
$W_{\tilde{r_1}+1}.f\equiv 0$. Since 
$f\in\mathcal S(G/((M\cap T)N))$, 
it follows from the assertion (i) for $\Gamma=\mathcal S$ 
that there exists some $f_S\in\mathcal S(G/((M\cap S)N))$ 
such that $R(f_S)\equiv f$. By the support theorem of $R$ 
in Proposition 5-7(4), 
we have $f\in C^\infty_d(G/((M\cap S)N))$. 
This yields the assertion (i) for $\Gamma=C^\infty_d$ with $d>0$. 
The assertion (i) for $\Gamma=C^\infty_c$ follows directly 
from the assertion for $\Gamma=C^\infty_d$ with $d>0$. 
\par 
Thus we have proved all the assertions. \qed 
\end{Proof}
\par 
\textbf{Remark}
 In a similar way to the proof 
of Proposition 5-5, 5-7 and Theorem 6-5
 we can prove the followings. 
Assume that $n\geq 1; 0\leq r_1,r_2 \leq m-n ; r_1\not= r_2$.
\par 
$(1)$ Assume that  $\dim G/((M\cap S)N)\leq\dim G/((M\cap T)N)$. 
Then the Radon transform $R_{M}$ maps 
$C^\infty(M/(M\cap S))$ injectively 
into $C^\infty(M/(M\cap T))$. 
Moreover, we have the following inversion formula for $R$.  
For $f_S\in C^\infty(M/(M\cap S))$, we have 
$$f_S(m((M\cap S)))=(R^+\bigm\vert_{M/(M\cap S)}(R_M(f_S)))(m((M\cap S)))
\text{ for }m\in M.$$ 
\par 
$(2)$ Assume that $\dim G/((M\cap S)N)<\dim G/((M\cap T)N)$. 
Then we have  
\begin{align*} 
& R_{M}(C^\infty(M/(M\cap S)))
=\{f\in C^\infty(M/(M\cap T))\bigm\vert 
((Ad(u_M)W^n_{\tilde r_1+1-n})).f\equiv 0\}. 
\\& 
=\{f\in C^\infty(M/(M\cap T))\bigm\vert 
\\& 
(Ad(u_M)W_{I}).f\equiv 0
\text{ for }I\in T_{2(\tilde r_1+1-n)}(m+n)
\text{ such that }I\subset(n+1,\ldots,m)\}.
\end{align*} 

\begin{flushleft}
\textbf{Contributions}.
\end{flushleft}
The author did all the work. 

\begin{flushleft}
\textbf{Data Availability}. 
\end{flushleft}
No datasets were generated or analyzed during the current study.

\begin{flushleft}
\textbf{Declarations}. 
\end{flushleft} 

\begin{flushleft}
\textbf{Ethics Declarations}. 
\end{flushleft} 
Not applicable.

\begin{flushleft}
\textbf{Competing Interests}. 
\end{flushleft} 
The author declares no competing interests.

\end{document}